\documentclass[12pt, leqno]{amsart}

\usepackage{graphicx}
\usepackage{amssymb,amsmath,amsthm,amscd}
\usepackage{mathrsfs}
\usepackage{enumerate}
\usepackage{enumitem}
\usepackage{verbatim}
\usepackage[usenames,dvipsnames]{color}
\usepackage[colorlinks=true, pdfstartview=FitV,
linkcolor=blue,citecolor=blue,urlcolor=blue]{hyperref}
\usepackage[all]{xy}
\usepackage{tikz}
\usetikzlibrary{shapes, arrows}
\usepackage{tikz-cd}
\usetikzlibrary{matrix, arrows}
\usepackage{stmaryrd} 
\usepackage{dsfont}
\usepackage{bm} 
\usepackage{scalerel}[2016/12/29] 
\usepackage{amsthm}
\usepackage{adjustbox}
\usepackage{float}
\usepackage{array}
\usepackage{multirow} 
\usepackage{ytableau} 
\usepackage{mathdots} 
\usepackage{xcolor}

\allowdisplaybreaks[4]

\DeclareMathSymbol{\shortminus}{\mathbin}{AMSa}{"39}

\newcommand{\nc}{\newcommand}

\numberwithin{equation}{section}

\usepackage[colorinlistoftodos]{todonotes}
\usepackage{mathtools}
\mathtoolsset{showonlyrefs,showmanualtags}
\usepackage{empheq}

\makeatletter
\newcommand{\reqnomode}{\tagsleft@false\let\veqno\@@eqno}

\theoremstyle{plain}
\newtheorem{lem}{Lemma}[section]
\newtheorem{pro}[lem]{Proposition}
\newtheorem{thm}[lem]{Theorem}
\newtheorem{cor}[lem]{Corollary}
\newtheorem{defi}[lem]{Definition}

\newcommand{\Pro}{\begin{pro}}
	\newcommand{\enpro}{\end{pro}}
\newcommand{\Lem}{\begin{lem}}
	\newcommand{\enlem}{\end{lem}}
\newcommand{\Thm}{\begin{thm}}
	\newcommand{\enthm}{\end{thm}}
\newcommand{\Cor}{\begin{cor}}
	\newcommand{\encor}{\end{cor}}
\newcommand{\Defi}{\begin{defi}}
	\newcommand{\enDefi}{\end{defi}}
\newcommand{\Proof}{\begin{proof}}
	\newcommand{\enproof}{\end{proof}}

\theoremstyle{definition} 
\newtheorem{rem}[lem]{Remark}

\newtheorem{Convention}[lem]{Convention}

\newcommand{\Conv}{\begin{Convention}}
	\newcommand{\enconv}{\end{Convention}}
\nc{\Rem}{\begin{rem}}
	\nc{\enrem}{\end{rem}}

\newcommand{\arxiv}[1]{\href{http://arxiv.org/abs/#1}{\tt arXiv:\nolinkurl{#1}}}

\newcommand{\monoto}{\hookrightarrow}
\nc{\epito}{\twoheadrightarrow}
\newcommand{\isoto}[1][]{\mathop{\xrightarrow%
		[{\raisebox{.3ex}[0ex][.3ex]{$\scriptstyle{#1}$}}]%
		{{\raisebox{-.6ex}[0ex][-.6ex]{$\mspace{2mu}\sim\mspace{2mu}$}}}}}

\nc{\rmkend}{\hfill$\triangledown$}
\nc{\defend}{\hfill$\triangle$}

\nc{\ccc}{\mathfrak{c}}
\nc{\CCC}{\mathfrak{C}}
\nc{\Ck}{\mathfrak{C}}

\nc{\kor}{\mathbb{C}}
\nc{\indx}{\mathbb{I}}

\nc{\CC}{C}
\nc{\cc}{c}
\nc{\sss}{s}
\nc{\ck}{\mathfrak{c}}

\nc{\Bg}{B}
\nc{\Ag}{A}
\nc{\As}{\pmb{\Ag}(z)}
\nc{\Asmone}{\pmb{\Ag}^{(-1)}(z)}
\nc{\Aps}{\pmb{\Ag}_+(z)}
\nc{\Apsone}{\pmb{\Ag}^{(1)}_+(z)}
\nc{\Apsmone}{\pmb{\Ag}^{(-1)}_+(z)}
\nc{\Ams}{\pmb{\Ag}_-(z)}
\nc{\Hg}{H}
\nc{\Thg}{\Theta}
\nc{\Thgs}{\boldsymbol{\Theta}}
\nc{\Thgsr}{\boldsymbol\grave{\boldsymbol\Theta}}
\nc{\bTh}{\pmb{\grave{\Theta}}(z)}
\nc{\avar}[1]{\mathbf{A}_{#1, +}(z)} 
\nc{\avarm}[1]{\mathbf{A}_{#1, +}^{(-)}(z)} 
\nc{\avarp}[1]{\mathbf{A}_{#1, +}^{(+)}(z)} 
\nc{\ivar}[2]{\pmb{\Ag}_{#1, +}^{#2}(z)} 
\nc{\thvar}[1]{\boldsymbol\grave{\boldsymbol\Theta}_{#1}(z)}

\nc{\smin}{{\shortminus}}

\nc{\Ui}{\widetilde{\mathbf{U}}^\imath}

\nc{\fext}[2]{{#1}[\negthinspace[#2]\negthinspace]}

\nc{\gTh}{\grave{\Theta}}

\nc{\DD}{\pmb{D}}
\nc{\KK}{\mathbb{K}}

\nc{\xgp}{x^+}
\nc{\Xgp}{\pmb{x}^+}
\nc{\xgm}{x^-}
\nc{\Xgm}{\pmb{x}^-}
\nc{\xgpm}{x^{\pm}}
\nc{\psig}{\phi^+}
\nc{\phig}{\phi^-}
\nc{\phipm}{\phi^\pm}
\nc{\Psig}{\pmb{\phi^+}}
\nc{\Phig}{\pmb{\phi^-}}
\nc{\Phipm}{\pmb{\phi^\pm}}
\nc{\hg}{h}

\nc{\bY}{\mathbf{Y}}
\nc{\bA}{\mathbf{A}}
\nc{\bT}{\mathbf{T}}

\nc{\Eg}{e}
\nc{\Kg}{K}
\nc{\ev}{{\rm ev}}

\nc{\Rep}{\on{Rep}} 

\DeclareRobustCommand{\sqbin}{\genfrac{[}{]}{0pt}{}}

\nc{\qq}{(q-q^{-1})^{-1}}
\nc{\factor}{\Omega}

\nc{\chring}{\Z[Y_a^{\pm 1}]_{a \in \C^\times}}
\nc{\chmap}{\chi_q}
\nc{\chmod}{\Z[\mathbf{Y}_a^{\pm 1}]_{a \in \C^\times}}
\nc{\ourchmap}{\boldsymbol{\chi}_q}
\nc{\Yring}{\mathcal{Y}} 
\nc{\Ymod}{\pmb{\mathcal{Y}}}
\nc{\Xmod}{\pmb{\mathcal{X}}}
\nc{\Yg}{\mathbf{Y}_{i,a}}
\nc{\Ys}{\mathbf{Y}_{i,\mathbf{s}}}
\nc{\Kmat}{\mathcal{K}^0} 

\nc{\Tr}{\on{Tr}}
\nc{\id}{\on{id}}

\nc{\ourR}{\mathbf{R}}
\nc{\ourQ}{\mathbf{Q}}
\nc{\ourY}{\mathbf{Y}}
\nc{\ourX}{\mathbf{X}}
\nc{\ourA}{\mathbf{A}}
\nc{\ourB}{\mathbf{B}}

\nc{\Oq}{\widetilde{\mathcal{O}}_q}
\nc{\Oqi}{\widetilde{\mathcal{O}}_q^{[i]}} 
\nc{\Oqii}{\widetilde{\mathcal{O}}_q^{[i,i+1]}} 
\nc{\Oqiii}{\widetilde{\mathcal{O}}_q^{[i,i+1,i+2]}} 
\nc{\oq}{\mathcal{O}_q}
\nc{\Oqc}{\mathcal{O}_q^{\mathbf{c}}(\widehat{\g})}
\nc{\car}{\mathcal{H}}
\nc{\qaa}{U_q(\widehat{\mathfrak{g}})}
\nc{\uqaa}{\widetilde{U}_q(\widehat{\mathfrak{g}})}
\nc{\qla}{U_q(L\mathfrak{g})}
\nc{\uqla}{\widetilde{U}_q(L\mathfrak{g})}
\nc{\drqaa}{{}^{\mathrm{Dr}}\qaa}
\nc{\uqsl}{U_q(\widehat{\mathfrak{sl}}_2)}
\nc{\uLsl}{U_qL\mathfrak{sl}_2}
\nc{\Serre}{\mathsf{Serre}}
\nc{\Sym}{\on{Sym}}
\nc{\UXp}{UX_+} 

\nc{\Tbr}{\mathbf{T}}
\nc{\degdr}{\on{deg}^{\mathrm{Dr}}}
\nc{\Uq}{\mathbf{U}}
\nc{\Uu}{\widetilde{\mathbf{U}}}
\nc{\gaf}{\widehat{\mathfrak{g}}}

\nc{\ada}{\ad_{\bhg_{-1}}}
\nc{\adb}{\ad_{\bhg_{1}}}
\nc{\adc}{\ad_{\bHg_{1}^{(2)}}}

\nc{\tAg}{\widetilde{A}}
\nc{\hAg}{\widehat{A}}

\nc{\bHg}{\overline{\Hg}}
\nc{\bhg}{\overline{\hg}}
\nc{\tHg}{\widetilde{H}}

\nc {\DrOq}{{}^{\mathrm{Dr}}\Oq}
\nc{\gup}[1]{^{(#1)}}
\nc{\gupp}{^{(1),+}}
\nc{\mi}{^{-1}}
\nc{\adh}{\operatorname{ad}_{\bar{\hg}_{-1}}}
\nc{\adhp}{\operatorname{ad}_{\bar{\hg}_{1}}}
\nc{\Omg}{\Omega^{-1}}
\nc{\Htwo}{\overline{H}_1^{(2)}}
\nc{\ad}{\operatorname{ad}}

\nc{\sspan}{\on{span}}

\newcommand{\commentout}[1]{}

\newcommand{\on}{\operatorname}
\nc{\be}{\begin{enumerate}}
	\nc{\ee}{\end{enumerate}}
\newcommand{\eq}{\begin{equation}}
	\newcommand{\eneq}{\end{equation}}
\nc{\bc}{\begin{cases}}
	\nc{\ec}{\end{cases}}
\newcommand{\eqn}{\begin{eqnarray*}}
	\newcommand{\eneqn}{\end{eqnarray*}}
\newcommand{\ba}{\begin{array}}
	\newcommand{\ea}{\end{array}}

\newcommand{\C}{{\mathbb C}}
\newcommand{\Q}{\mathbb {Q}}
\newcommand{\Z}{{\mathbb Z}}

\newcommand{\g}{{\mathfrak{g}}}
\nc{\Ad}{\operatorname{ad}}

\nc{\gr}{\on{gr}}

\newcommand{\End}{\operatorname{End}}
\nc{\Aut}{\operatorname{Aut}}

\nc{\coker}{\operatorname{coker}}

\nc{\Img}{\on{Im}}
\nc{\res}{\on{res}}

\nc{\modv}[1]{{#1}\operatorname{-mod}}

\nc{\bl}{\bigl(}
\nc{\br}{\bigr)}
 
\definecolor{RED}{rgb}{1,0,0}

\newlength{\mylength}
\DeclareRobustCommand{\SkipTocEntry}[5]{}

\AtBeginDocument{%
   \def\MR#1{}
}

\title[Compatibility of Drinfeld presentations for QSP]{Compatibility of Drinfeld presentations and $q$-characters for affine Kac--Moody quantum symmetric pairs: quasi-split case}

\author[J.-R. Li]{Jian-Rong Li}
\address{Faculty of Mathematics, University of Vienna, Oskar Morgenstern Platz 1, 1090 Vienna, Austria, OrciD: 0000-0001-7896-7391}
\email{\href{mailto:lijr07@gmail.com}{lijr07@gmail.com}}

\author[T. Prze\'{z}dziecki]{Tomasz Prze\'{z}dziecki}
\address{Faculty of Mathematics, University of Vienna, Oskar Morgenstern Platz 1, 1090 Vienna, Austria, OrciD: 0000-0001-9700-1007}
\email{\href{mailto:tprzezdz@exseed.ed.ac.uk}{tprzezdz@exseed.ed.ac.uk}}



\thanks{The first author was supported by the Austrian Science Fund (FWF), P-34602, Grant DOI: 10.55776/P34602, and PAT 9039323, Grant-DOI 10.55776/PAT9039323, and the National Natural Science Foundation of China (No. 12471023). The second author was supported by the EPSRC grant No.\ EP/W022834/1 \emph{Kac--Moody quantum symmetric pairs, KLR algebras and generalized Schur--Weyl duality}, and the Austrian Science Fund (FWF), PAT 9039323, Grant-DOI 10.55776/PAT9039323.}

\begin{document}

\begin{abstract}
Let $(\Uq, \Uq^\imath)$ be a quasi-split affine quantum symmetric pair of type $\mathsf{AIII}$. 
This case is of particular interest thanks to the existence of geometric realizations and Schur--Weyl dualities. 
We establish factorization and coproduct formulae for the Drinfeld--Cartan series $\boldsymbol\Theta_i(z)$ in the Lu--Pan--Wang--Zhang `new Drinfeld'-style presentation, generalizing the split type results from \cite{Przez-23, LP25}. 
As an application, we construct a boundary analogue of the $q$-character map, and show that it is compatible with Frenkel and Reshetikhin's original $q$-character homomorphism.  
\end{abstract}

\maketitle

\setcounter{tocdepth}{1}
\tableofcontents


\section{Introduction}

\renewcommand{\thelem}{\Alph{lem}}
\setcounter{lem}{0}

\nc{\aaui}{\mathbf{U}^\imath_{\mathsf{aff}}}
\nc{\ichmap}{{}^b\boldsymbol\chi_q}
\nc{\Uhz}{\fext{U_q(\widetilde{\mathfrak{h}})}{z}} 
\nc{\auiv}[1]{\Ui_{\mathsf{aff},#1}}
\nc{\auivv}[1]{\mathbf{U}^\imath_{\mathsf{aff},#1}}
\nc{\eva}{\on{ev}_a}
\nc{\aui}{\Ui_{\mathsf{aff}}}

\subsection{Loop presentations and $q$-characters} 

It is well known that quantum affine algebras admit three distinct presentations: the `Drinfeld--Jimbo', `new Drinfeld' (or `loop') and `RTT' realizations. The first quantizes the usual Serre-type presentation of a Kac-Moody Lie algebra $\widehat{\g}$, while the second is a quantization of the central extension presentation of $\widehat{\g}$. The interplay between these different realizations was investigated and described precisely in \cite{Ding-Fr, FrMukHopf}. One of the key features of the new Drinfeld presentation is that it exhibits a large, infinitely generated, commutative subalgebra of $\Uq := U_q(L\g)$. The spectra of the generators of this subalgebra, often called Drinfeld--Cartan operators $\phi^\pm_{i,m}$, play a key role in the classification of finite dimensional representation of $\Uq$ via Drinfeld polynomials \cite{chari-pressley-qaa, chari-pressley-94}, and $q$-character theory \cite{Knight, FrenRes, FrenMuk-comb}. 

Recently, Lu and Wang \cite{lu-wang-21}, building on the work of Baseilhac and Kolb \cite{bas-kol-20} in the split rank one case (see also \cite{ZhangDr}), have constructed a Drinfeld-type presentation for split affine quantum symmetric pair coideal subalgebras, also known as $\imath$quantum groups $ \Uq^\imath$. This work was extended to the quasi-split case in \cite{LWZ-I, LWZ-quasi, LPWZ}. 
The Lu--Wang presentation shares many familiar features of Drinfeld's new realization. In particular, it also exhibits a large commutative subalgebra, generated by $\Theta_{i,m}$. 

In \cite{Przez-23, LP25}, we investigated the relationship between the Lu--Wang and new Drinfeld presentations of $\Uq^\imath$ and 
$\Uq$, respectively, in the \emph{split} case. Roughly speaking, we proved that the two presentations are compatible modulo the Drinfeld positive half $\Uq_+$ of $\Uq$. 
More specifically, we proved that the (renormalized) Cartan currents~$\Thgsr_i(z)$ in the Lu--Wang presentation are well-behaved with respect to the inclusion $ \Uq^\imath \hookrightarrow \Uq$ and the coproduct, in the sense that they satisfy 
the \emph{factorization property}: 
\eq \label{intro 1}
\Thgsr_i(z) \equiv \boldsymbol{\phi}_i^-(z\mi)\boldsymbol{\phi}_i^+(C z) \quad \mod \fext{\Uq_+}{z}, 
\eneq
and the \emph{`approximately group-like' property}: 
\eq \label{intro 2}
\Delta(\Thgsr_i(z)) \equiv \Thgsr_i(z) \otimes \Thgsr_i(z) \quad  \mod \fext{\Uq^\imath \otimes \Uq_{+}}{z}. 
\eneq 
These results allowed us to define a notion of \emph{boundary $q$-characters} for affine quantum symmetric pairs, compatible with the usual $q$-character homomorphism due to Frenkel and Reshetikhin \cite{FrenRes}. 
Boundary $q$-characters of evaluation modules in type $\mathsf{AI}$ were investigated in detail in \cite{LP25b}. 

\subsection{Quantum affine symmetric pairs of type $\mathsf{AIII}$} 
\label{sec: QSP AIII}

The goal of the present paper is to generalize the aforementioned results to quasi-split affine quantum symmetric pair of type $\mathsf{AIII}_N^{(\tau)}$. This case is of particular interest, thanks to the existence of geometric realizations, Schur--Weyl dualities and categorification results. 

A geometric interpretation of $\mathbf{U}^\imath$, when $N$ is odd, was recently obtained in \cite{SuWang}. The construction uses equivariant $K$-groups of non-connected Steinberg varieties of type $\mathsf{C}$, endowed with a convolution product. This, in turn, paves the way to the application of standard techniques from geometric representation theory \cite{CGinz}, based on the Beilinson--Bernstein--Deligne--Gabber decomposition theorem, to study problems such as the classification of simple modules and Kazhdan--Lusztig-type multiplicities. 

A Schur--Weyl duality, in the form of a double centralizer theorem, between coideal subalgebras of affine type $\mathsf{AIII}$ and affine Hecke algebras of type $\mathsf{C}$, with three parameters, was established in \cite{FLLLWW}. This relationship was further exploited in \cite{AppelPrz}, where a functor was constructed from the category of modules over \emph{orientifold} Khovanov--Lauda--Rouquier (KLR) algebras to the category of $\mathbf{U}^\imath$-modules. Orientifold KLR algebras can be seen as graded analogues of affine Hecke algebras of type $\mathsf{C}$. They were introduced in \cite{VV-Hecke}, and basic results about their representation theory were established in \cite{PrzPac}. It is expected that the representation theory of $\mathbf{U}^\imath$ can also be approached via this functor, which should, additionally, reveal new features, such as the existence of non-trivial gradings on categories of $\mathbf{U}^\imath$-modules. 

Quantum symmetric pair coideal subalgebras of affine type $\mathsf{AIII}_N^{(\tau)}$ are also connected to, via Schur--Weyl duality and orientifold KLR algebras, to categorification results due to Enomoto and Kashiwara \cite{EK1, EK2, EK3, Enomoto} in their work on symmetric crystals. 
More precisely, inspired by Ariki--Lascoux--Leclerc--Thibon theory, they constructed a certain remarkable module ${}^\theta V(\lambda)$ over a $q$-boson-type algebra. Just as KLR algebras categorify $U_q(\mathfrak{n}_-)$, orientifold KLR algebras were later shown to categorify ${}^\theta V(\lambda)$ \cite{VV-Hecke}. Recently, there has been renewed interest in the work of Enomoto and Kashiwara. For instance, in \cite{SSX} (see also \cite{Wang26}), the Enomoto--Kashiwara module was proposed as the Higgs branch of a certain three-dimensional $\mathcal{N}=4$ supersymmetric gauge theory, mirror-symmetric dual to $\mathbf{U}^\imath$, which, in turn, arises as its Coulomb branch.  
It is plausible that Enomoto--Kashiwara type categorifications can also be constructed independently in terms of quantum symmetric pairs, and extended to other Satake types (e.g., $\mathsf{AI}$).

\subsection{Main results} 

Let $(\Uq, \Uq^\imath)$ be a quasi-split affine quantum symmetric pair of type $\mathsf{AIII}_N^{(\tau)}$, where $N \geq 1$ (both odd and even cases are allowed). 
The following is the main result of the paper. 

\Thm[Theorem \ref{thm: main overall}] \label{thm: A intro}
The `Drinfeld--Cartan' series $\thvar{i}$ admit the following factorization: 
\begin{align*}
\thvar{i} \equiv K_i K_{\tau(i)}\mi   \boldsymbol\phi^-_i(z\mi)\boldsymbol{\phi}^+_{\tau(i)}(\CCC z)  \mod \fext{\Uq_+}{z}. 
\end{align*}
\enthm 

The proof uses the key methods introduced in our previous works \cite{Przez-23, LP25}: reduction to rank one, and analysis of the compatibility of the relative braid group action on $\mathbf{U}^\imath$ with the usual Lusztig braid group action on $\Uq$. Nevertheless, type $\mathsf{AIII}$ presents new challenges and requires substantial refinement of these methods. In particular, we need to deal with two new rank one cases, corresponding to a pair of vertices $i, \tau(i)$ of the Satake diagram in the orbit of the associated involution $\tau$ being either isolated or connected by an edge. The third rank one case, corresponding to $i = \tau(i)$, gives rise to the $q$-Onsager algebra and was already considered in \cite{Przez-23}.   
The case when $i$ and $\tau(i)$ are connected, which yields a quasi-split coideal subalgebra of $U_q(L\mathfrak{sl}_3)$, is the hardest, and is handled separately in \S \ref{sec: new rank 1}. 

The higher rank factorization formula is then proven by reducing to the rank one cases, using results on the compatibility of the braid group operators corresponding to fundamental weights, established in \S \ref{sec: compat of br}. This part of the proof also presents novel features, largely due to the appearance of the \emph{relative} braid group and root systems of type $\mathsf{C}_n^{(1)}$ and $\mathsf{A}_{2n}^{(2)}$. In contrast, in our previous works, dedicated to the split case, the relative braid group always coincided with the full braid group. 

The next result is a crucial application of Theorem \ref{thm: A intro}. 

\Thm[Theorem \ref{thm: coproduct main}] 
The `Drinfeld--Cartan' series $\thvar{i}$ satisfy the following `group-like' coproduct formula: 
\begin{align*}
\Delta(\thvar{i}) \equiv \thvar{i} \otimes \thvar{i} \qquad \mod \fext{(\mathbf{U}^\imath \otimes \Uq_+)}{z}. 
\end{align*}
\enthm 

Theorem \ref{thm: A intro} also allows us to compute the generalized eigenvalues of $\thvar{i}$ on finite-dimensional representations of $\Uq$, restricted to $\mathbf{U}^\imath$. 

\Thm[Theorem \ref{cor: FR thm Oq}] \label{thm: C intro}
Let $W$ be a finite-dimensional representation of $\Uq$. 
Then the generalized eigenvalues of $\thvar{i}$ on $W$ are of the form: 
\eq 
\boldsymbol\gamma_i(z) = q^{\deg \mathbf{Q}_i - \deg \mathbf{Q}_{\tau(i)}}  \kappa_i \kappa_{\tau(i)}\mi \frac{\ourQ_i(q^{-1}z)}{\ourQ_i(qz)} \frac{\ourQ_{\tau(i)}^\dag(q z)}{\ourQ_{\tau(i)}^\dag(q^{-1}z)}, 
\eneq
for some explicit polynomials $\ourQ_i(z)$ and scalars $\kappa_i$. 
\enthm 

The final and most important application of Theorems \ref{thm: A intro}--\ref{thm: C intro} is the construction of a `\emph{boundary $q$-character theory}' for quantum affine symmetric pairs of type $\mathsf{AIII}_N^{(\tau)}$, compatible with Frenkel and Reshetikhin's $q$-character theory for quantum affine algebras. 

We define the boundary $q$-character $\ichmap(W)$ of a finite-dimensional representation~$W$ of $\mathbf{U}^\imath$ to be the trace\footnote{More precisely, we first evaluate the left leg of the tensor $\mathcal{K}^0$ on $W$ and then take the trace over~$W$, i.e., $\ichmap(W) =\Tr_W((\pi_{W(z)} \otimes 1)(\mathcal{K}^0))$.} of the following explicit operator 
\[
\mathcal{K}^0 = \exp \left( - (q-q\mi) \sum_{i \in \indx_0} \sum_{k >0} \frac{k}{[k]} \grave{H}_{i,k} \otimes \tilde{h}_{i,-k} z^k \right) \cdot TT_\tau\mi, 
\] 
reminiscent of the formula for the abelian part of the universal $R$-matrix (see \eqref{eq: K0} for a more detailed explanation of the notation). This definition is shown to be equivalent to the `more naive' construction of the boundary $q$-character as the series counting the dimensions of the generalized eigenspaces with respect to the action of the commutative algebra generated by the elements $\Theta_{i,r}$. 

Since $\mathbf{U}^\imath$ is a coideal of $\Uq$, the Grothendieck group ${[\Rep \mathbf{U}^\imath]}$ of finite-dimensional representations of $\mathbf{U}^\imath$ is a module over the Grothendieck ring ${[\Rep \Uq]}$. The compatibility of boundary $q$-characters with ordinary $q$-characters can be formulated in the following way. 

\Thm[Theorem \ref{cor: comm diagram qchar actions}] \label{thm: D intro}
Boundary $q$-characters are compatible with usual $q$-characters, i.e., the following diagram commutes: 
\[
\begin{tikzcd}[ row sep = 0.2cm]
{[\Rep \Uq]}  \arrow[r, "\chmap"] & \Z[Y_{i,a}^{\pm 1}]_{i \in \indx_0, a \in \C^{\times}}  \\
 \curvearrowright  & \curvearrowright  \\
{[\Rep \mathbf{U}^\imath]} \arrow[r, "\ichmap"] & \Uhz. 
\end{tikzcd}
\]
\enthm 

In more `down to earth' terms, if $V \in \Rep \mathbf{U}^\imath$ and $W \in \Rep \Uq$, then 
Theorem \ref{thm: D intro} means that $\ichmap(V \otimes W) = \ichmap(V) \triangleleft \chi(W)$, for an appropriately defined right action $\triangleleft$ of $\Z[Y_{i,a}^{\pm 1}]_{i \in \indx_0, a \in \C^{\times}}$ (see \eqref{eq: Yring module}). 

\subsection{Future directions} 

In future works, we intend to explore how boundary $q$-characters interact with the salient features of quantum affine symmetric pairs of type $\mathsf{AIII}$ outlined in \S \ref{sec: QSP AIII}, including geometric realizations, Schur--Weyl duality and categorification. 

In particular, we expect that the boundary $q$-characters introduced in this paper admit a geometric interpretation \`{a} la Nakajima \cite{Naka-qt2, Naka-qt}. In the context of twisted Yangians, a similar expectation was recently formulated in \cite{Naka-new}. Roughly speaking, boundary $q$-characters should appear as generating functions of the Poincar\'{e} polynomials of an analogue of graded quiver varieties. Such an analogue was, for instance, proposed in  \cite{YiqiangLi}, and the cotangent bundles to flag varieties of type $\mathsf{C}$, considered in \cite{SuWang}, are an example thereof. The homological grading should also give rise to a natural $t$-deformation of boundary $q$-characters.  
We expect that these boundary $t,q$-characters will correspond, under Schur--Weyl duality, to the graded  characters of modules over orientifold KLR algebras. 

\subsection{Organization} 

In \S \ref{sec: QSP}, we recall the definitions of the main objects considered in the paper, including quantum affine algebras, quantum symmetric pairs, braid group actions and loop presentations. This section does not contain new material. In \S \ref{sec: new rank 1}, we prove a factorization formula (Theorem \ref{thm: rank 1 factorization aij-1}) for the series $\thvar{i}$ in the special case of the rank one quantum symmetric pair coideal subalgebra $\Ui(\widehat{\mathfrak{sl}}_3, \tau)$. This formula is new, and complements the previous rank one result for the $q$-Onsager algebra established in \cite{Przez-23}. In \S \ref{sec: roots n weight}, we recall some combinatorial facts about fundamental weights in root systems of type $\mathsf{C}_n^{(1)}$ and $\mathsf{A}_{2n}^{(2)}$. Section \ref{sec: compat of br} establishes a result about the compatibility of braid group actions, modulo the positive half of $\Uq$. Section \ref{sec: factor formula} is devoted to the proof of the main result of the paper, i.e., a factorization formula for the series $\thvar{i}$ in arbitrary rank (Theorem \ref{thm: main overall}). 
In the last two sections, we propose several applications of our factorization formula. In \S \ref{sec: coproduct}, we derive a formula for the coproduct $\Delta(\thvar{i})$ (Theorem \ref{thm: coproduct main}). The factorization and coproduct formulae are then used in \S \ref{sec: bound q-char} to compute the generalized eigenvalues of the series $\thvar{i}$ on restriction representations (Theorem \ref{cor: FR thm Oq}), and show that boundary $q$-characters are compatible with Frenkel and Reshetikhin's $q$-characters (Theorem \ref{cor: comm diagram qchar actions}). 
Some of the more cumbersome calculations have been moved to the appendices to improve the readability of the paper. 

\addtocontents{toc}{\SkipTocEntry}
	
\section*{Acknowledgements} 
We would like to thank Changjian Su for insightful discussions on this and related projects. We would also like to thank the  International Centre for Mathematical Sciences (ICMS) in Edinburgh and the Erwin Schr\"{o}dinger Institute (ESI) in Vienna, where parts of this project were carried out, for their hospitality. 

\renewcommand{\thelem}{\thesection.\arabic{lem}}
\setcounter{lem}{0}

\section{Quasi-split quantum symmetric pairs} 
\label{sec: QSP}
\nc{\sln}{\mathfrak{sl}_{n+1}}
\nc{\asln}{\widehat{\mathfrak{sl}}_{n+1}}

We work over the field of complex numbers and assume that $q \in \C^\times$ is not a root of unity throughout. 
We will also abbreviate
\[
[k] = \frac{q^k-q^{-k}}{q-q\mi}, \qquad [x,y]_q = xy - qyx, \qquad \rho = q - q\mi. 
\]

\subsection{Quantum affine algebras}
\nc{\Ka}{\widetilde{K}}
\nc{\Kb} {\widetilde{K}'}

Let $\indx_0 = \{ 1,\cdots, N \}$ and $\indx = \indx_0 \cup \{0\}$. Let $\g = \mathfrak{sl}_{N+1}$, with Cartan matrix $(a_{ij})_{i,j \in \indx_0}$, and let $\widehat{\g}$ be the corresponding untwisted affine Lie algebra with affine Cartan matrix $(a_{ij})_{i,j \in \indx}$.
We use standard notations and conventions regarding root systems, weight lattices, Weyl groups, etc. In particular, we let $\mathcal{R}$ denote the affine root system associated to $\widehat{\g}$, with positive roots $\mathcal{R}^+$ and simple roots $\alpha_i$ ($i \in \indx$); let $\theta$ denote the highest root, and $\delta$ the basic imaginary root; let $P$ and $Q$ denote the weight and root  lattices of $\g$, respectively; and let $\omega_i \in P$ ($i \in \indx_0$) be the fundamental weights.

The \emph{Drinfeld double quantum affine algebra} $\Uu' = \Uu'(\widehat{\g})$ is the algebra with generators $\Eg_i^\pm, \Ka_i, \Kb_i$ $(i \in \indx)$, where $\Ka_i$ and $\Kb_i$ are invertible, subject to the relations: 
\begin{align} \label{eq: DDrel 1}
\Ka_i\Eg_j^{\pm} &= q^{\pm a_{ij}} \Eg_j^\pm\Ka_i, & 
\Kb_i\Eg_j^{\pm} &= q^{\mp a_{ij}} \Eg_j^\pm\Kb_i, \\ \label{eq: DDrel 2} \tag{\theequation} \stepcounter{equation}
[\Ka_i,\Ka_j] &= [\Kb_j, \Kb_i] = [\Ka_i, \Kb_j] = 0, & 
[\Eg_i^{+}, \Eg_j^{-}] &= \delta_{ij} \frac{\Ka_i - {\Kb_i}}{q - q^{-1}}, 
\end{align}
\eq \label{eq: DDrel 3}
\Serre_{ij}(\Eg_i^{\pm}, \Eg_j^\pm) = 0, 
\eneq
where
\[
\Serre_{ij}(x,y) = \begin{cases} 
x^3y - [3] x^2yx +[3]xyx^2 - yx^3 & \mbox{if} \ a_{ij} = -2, \\
x^2 y - [2] xyx + yx^2 & \mbox{if} \ a_{ij} = -1, \\
xy - yx & \mbox{if} \ a_{ij} = 0, \\
0 & \mbox{if} \ a_{ij} = 2. 
\end{cases}
\]
We will also abbreviate 
$E_i = e_i^+$, $F_i = e_i^-$. 
The algebra $\Uu'$ is a Hopf algebra, with the coproduct
\begin{align} \label{eq: coprod on U}
\Delta(\Eg_i^+) &= \Eg_i^+ \otimes 1 + \Ka_i \otimes \Eg_i,  &\quad \Delta(\Eg_i^-) &= \Eg_i^- \otimes \Kb_i + 1 \otimes \Eg_i^-, \\
 \Delta(\Ka_i) &= \Ka_i \otimes \Ka_i,  &\quad \Delta(\Kb_i) &= \Kb_i \otimes \Kb_i,
\end{align}
and the counit given by $\varepsilon(\Eg_i^{\pm}) = 0$ and $\varepsilon(\Ka_i) = \varepsilon(\Kb_i) = 1$.

\nc{\uaff}{\mathbf{U}_{\mathsf{aff}}} 

The \emph{quantum affine algebra} $\Uq'$ is the algebra generated by $\Eg_i^\pm, K_i, K_i\mi$ $(i \in \indx)$ with relations \eqref{eq: DDrel 1}--\eqref{eq: DDrel 3}, modified by replacing $\Ka_i$ with $K_i$, and $\Kb_i$ with $K_i\mi$. 
Given $\mathbf{c} = (c_i) \in (\C^\times)^{\indx}$, let 
\[
\Uq'_{\mathbf{c}} = \Uu'  / \langle \Ka_i \Kb_i - c_i \mid i \in \indx  \rangle. 
\footnote{We incorporate the parameters $\mathbf{c}$ already at this stage so that \eqref{eq: Kolb emb} naturally descends to \eqref{eq: Kolb emb 2}.}
\] 
We will use the same notation for elements of $\Uu'$ and their images in the quotient $\Uq'_{\mathbf{c}}$. 
There is an isomorphism 
\eq \label{eq: 2 qaa iso}
\Uq' \mapsto \Uq'_{\mathbf{c}}, \qquad F_i \mapsto F_i, \quad E_i \mapsto c_i^{-1/2} E_i, \quad K_i \mapsto c_i^{-1/2} \Ka_i, \quad  K_i\mi \mapsto c_i^{-1/2} \Kb_i. 
\eneq
We will also consider the \emph{quantum loop algebra}\footnote{The passage to the quantum loop algebra is motivated by the fact that its Drinfeld--Cartan subalgebra is actually commutative.} 
and its Drinfeld double, defined as the quotients 
\[
 \Uq = \Uq' / \langle K_\delta - 1 \rangle, \qquad \Uq_{\mathbf{c}} = \Uq'_{\mathbf{c}} / \langle \Ka_\delta - c_\delta^{1/2} \rangle,  \qquad \Uu = \Uu' / \langle \Ka_\delta - c_\delta^{1/2}, \Kb_\delta - c_\delta^{1/2} \rangle, 
\]
where $c_\delta = \prod_{i \in \indx} c_i$. 
The isomorphism \eqref{eq: 2 qaa iso} descends to an isomorphism $\Uq \cong \Uq_{\mathbf{c}}$. 

\subsection{Lusztig's braid group action} 
\label{subsec: br action Lusztig}

The Weyl group $W_0$ of $\g$ is generated by the simple reflections $s_i$ ($i \in \indx_0$), and acts on $P$ in the usual way. The extended affine Weyl group $\widetilde{W}$ is the semi-direct product $W_0 \ltimes P$. 
It contains the affine Weyl group $W = W_0 \ltimes Q = \langle s_i \mid i \in \indx \rangle$ as a subgroup, and can also be realized as the semi-direct product 
\eq \label{eq: ext aff weyl 2 rel}
\widetilde{W} \cong \Lambda \ltimes W,
\eneq
where $\Lambda = P/Q$ is a finite group of automorphisms of the Dynkin diagram of $\widehat{\g}$. Let $\widetilde{\mathcal{B}}$ be the braid group associated to $\widetilde{W}$. 
When referring to specific types, we will use notation such as $\widetilde{W}_{\mathsf{A}_n^{(1)}}$, $\widetilde{W}_{\mathsf{C}_n^{(1)}}$, etc. 

It is well known \cite[\S 37.1.3]{lusztig-94} that,
for each $i \in \indx$, there exists an automorphism $T_i$ of $\Uu$ such that $T_i(\Ka_j) = \Ka_i^{-a_{ij}}\Ka_j$, $T_i(\Kb_j) = (\Kb_i)^{-a_{ij}}\Kb_j$ and 
\begin{align*}
T_i(E_j) &= \begin{cases}
 -F_i(\Kb_i)\mi & \mbox{if } j=i, \\ 
[E_i, E_j]_{q\mi} & \mbox{if } a_{ij} = -1, \\ 
[2]\mi E_i^2 E_j - q\mi E_i E_j E_i + q^{-2} [2]\mi E_j E_i^2 & \mbox{if } a_{ij} = -2, \\ 
E_j & \mbox{if } a_{ij} = 0, 
\end{cases} \\
T_i(F_j) &= \begin{cases}
 -\Ka_i\mi E_i & \mbox{if } j=i, \\ 
[F_j, F_i]_{q} & \mbox{if } a_{ij} = -1, \\ 
[2]\mi F_j F_i^2 - q F_i F_j F_i + q^2 [2]\mi F_i^2 F_j & \mbox{if } a_{ij} = -2, \\ 
F_j & \mbox{if } a_{ij} = 0. 
\end{cases}
\end{align*}
For each $\lambda \in \Lambda$, there is also an automorphism $T_\lambda$ such that $T_\lambda(E_i) = E_{\lambda(i)}$, $T_\lambda(F_i) = F_{\lambda(i)}$, $T_\lambda(\Ka_i) = \Ka_{\lambda(i)}$ and $T_\lambda(\Kb_i) = \Kb_{\lambda(i)}$.   
These automorphisms satisfy the relations of the braid group $\widetilde{\mathcal{B}}$, and descend to automorphisms of $\Uq$. 
Given a reduced expression $w = \lambda s_{i_1} \cdots s_{i_r} \in \widetilde{W}$, one defines $T_w = T_\lambda T_{i_1} \cdots T_{i_r}$. This is independent of the choice of reduced expression. 

\subsection{Drinfeld's new presentation}

Let us recall the ``new" Drinfeld presentation of~$\Uu$. By (a slight adaptation of) \cite{drinfeld-dp, beck-94}, $\Uu$ is isomorphic to the algebra generated by the pairwise commuting elements $\hg_{i,l}, \Ka_{i}, \Kb_i$ ($i \in \indx_0$, $l \in \Z - \{0\}$) and  $x_{i,k}^{\pm}$ ($i \in \indx_0, \ k \in \Z$), subject to the following relations:
\begin{align} 
\Ka_i \xgpm_{j,k} =& \ q^{\pm a_{ij}} \xgpm_{j,k}\Ka_i,  &
\Kb_i \xgpm_{j,k} =& \ q^{\mp a_{ij}} \xgpm_{j,k}\Kb_i,  \\
[\hg_{i,k},x^-_{j,l}] =& - \textstyle \frac{[k\cdot a_{ij}]}{k}x^-_{j,k+l}, &
[\hg_{i,k},x^+_{j,l}] =& \ c_{\delta}^{-k/2} \textstyle \frac{[k\cdot a_{ij}]}{k}x^+_{j,k+l}
,\\ \label{eq: rel sln xx} \ \ \ 
[\xgpm_{i,k+1}, \xgpm_{j,l}]_{q^{\pm a_{ij}}}  =& \ q^{\pm a_{ij}}[\xgpm_{i,k},\xgpm_{j,l+1}]_{q^{\mp a_{ij}}}, &
[\xgp_{i,k}, \xgm_{j,l}] =& \ \delta_{ij} \textstyle\frac{c_i^{1/2}c_\delta^{k/2}}{q-q^{-1}}(\psig_{i, k+l} - \phig_{i, k+l}), 
\end{align}
\begin{align}
\Sym_{k_1, k_2} \big( \xgpm_{j,l} \xgpm_{i,k_1} \xgpm_{i,k_2} - [2]  \xgpm_{i,k_1} \xgpm_{j,l} \xgpm_{i,k_2} +  \xgpm_{i,k_1} \xgpm_{i,k_2} \xgpm_{j,l} \big) =& \ 0 \quad (a_{ij} = -1), \\ 
[\xgpm_{i,k}, \xgpm_{j,l}] =& \ 0 \quad (a_{ij} = 0), 
\end{align} 
where $\Sym_{k_1, k_2} $ denotes symmetrization with respect to the indices $k_1, k_2$ and 
\begin{align}
\boldsymbol\phi^{+}_i(z) &= \sum_{k=0}^\infty \phipm_{i,\pm k} z^{\pm k} = c_i^{-1/2} \Ka_i \exp\left( (q-q^{-1}) \sum_{k=1}^\infty \hg_{i, k} z^{k} \right), \\ 
\boldsymbol\phi^{-}_i(z) &= \sum_{k=0}^\infty \phipm_{i,\pm k} z^{\pm k} = c_i^{-1/2} \Kb_i \exp\left( - (q-q^{-1}) \sum_{k=1}^\infty \hg_{i,- k} z^{- k} \right). 
\end{align}

Explicitly, the Drinfeld new generators can be described in terms of the braid group action:
\[
x_{i,k}^{\pm} = o(i)^k T_{\omega_i}^{\mp k}(e_i^{\pm}) \qquad (i \in \indx_0, \ k \in \Z), 
\]
for a choice of sign function $o(i)$ such that $o(i)o(j) = -1$ whenever $a_{ij} < 0$. 
We will use the following notation for generating series: 
\[
\mathbf{x}^{\pm}_{i, \leq r}(z) = \sum_{k\leq r} x_{i,k} z^k, \quad 
\mathbf{x}^{\pm}_{i, -}(z) = \mathbf{x}^{\pm}_{i, \leq 0}(z), \quad 
\mathbf{x}^{\pm}_{i, \geq r}(z) = \sum_{k\geq r} x_{i,k} z^k, \quad 
\mathbf{x}^{\pm}_{i, +}(z) = \mathbf{x}^{\pm}_{i, \geq 0}(z). 
\]


We consider $\Uq$ as a $Q$-graded algebra with 
\[
\degdr \xgpm_{i,k} = \pm \alpha_i, \quad \degdr h_{i,k} = 0 \qquad (i \in \indx_0). 
\] 
This grading can be extended to $\Uu$ by setting 
\[
\degdr \Ka_i = \degdr \Kb_i = 0 \qquad (i \in \indx). 
\]
Note that 
\[
\degdr e_0^{\pm} = \mp\theta. 
\]

We will need to consider subalgebras of $\Uu$ spanned by elements of specified degree. 
Write $\degdr_i$ for degree along the simple root $\alpha_i$. Let $Q_i$ be the sublattice of~$Q$ spanned by all the simple roots except $\alpha_i$. 
Let $\Uu_{d_i=r,+}$ be the subalgebra of $\Uu$ spanned by elements which are both: (i) of degree $r$ along $\alpha_i$, and (ii) positive $Q_i$-degree, i.e., non-negative along all the simple roots and positive along at least one simple root in $Q_i$. That is, 
\[
\Uu_{d_i=r,+} = \langle y \in \Uu \mid \degdr_i y = r;\ \forall j \neq i \in \indx_0\colon \degdr_j y \geq 0;\ \exists j \neq i\colon \degdr_j y > 0 \rangle.
\]
Moreover, set 
\begin{align*} 
\Uu_{d_i\geq r} =& \ \langle y \in \Uu \mid \degdr_i y \geq r \rangle, \\
\Uu_{d_i\geq r,+} =& \ \langle y \in \Uu \mid \degdr_i y \geq r;\ \forall j \neq i \in \indx_0\colon \degdr_j y \geq 0;\ \exists j \neq i\colon \degdr_j y > 0 \rangle, \\
\Uu_{d_i \geq r, d_j\geq s,+} =& \ \langle y \in \Uu \mid \degdr_i y \geq r, \degdr_j y \geq s;\ \forall k \neq i,j \in \indx_0\colon \degdr_k y \geq 0; \\
&\quad \exists k \neq i,j \colon \degdr_k y > 0 \rangle, \\
\Uu_+ =& \ \langle y \in \Uu \mid \forall j \in \indx_0\colon \degdr_j y \geq 0,\ \exists j \colon \degdr_j y > 0 \rangle.  
\end{align*} 
Analogous notations are used with `$\geq$' replaced by `$\leq$', or `$+$' replaced by `$-$', etc. 
We also abbreviate $\Uu_{d_i, d_j\geq r,+} = \Uu_{d_i \geq r, d_j\geq r,+}$, and use analogous notations for the corresponding subalgebras of $\Uq$.

\subsection{Quantum symmetric pairs of quasi-split type}

\nc{\Uiv}[1]{\Ui_{#1}}
\nc{\auic}{\mathbf{U}^\imath_{\mathbf{c}}}
\nc{\uic}{\mathbf{U}^\imath_{\mathbf{c}}}
\nc{\DrUi}{{}^{\mathrm{Dr}}\Ui}
\nc{\ui}{\mathbf{U}^\imath} 

\nc{\tindx}{\indx_\tau}

In this paper, we consider quantum symmetric pair coideal subalgebras $\Ui$ of proper quasi-split affine restrictable type $\mathsf{A}$, i.e., corresponding to Satake diagrams with no black nodes and a non-trivial involution $\tau$ which fixes the affine node. They are also known as quantum symmetric pairs (or $\imath$quantum groups) of affine type $\mathsf{AIII}$.  
The resulting types of Satake diagrams are illustrated in Figure~\ref{fig: Satake}. 
As a limit case, we will also consider the $q$-Onsager algebra, which is of type $\mathsf{AI}$.  
Let $\tindx$ be a fixed set of representatives of $\tau$-orbits on $\indx$, and let $\indx_{0,\tau} = \indx_0 \cap \tindx$. 

\newcommand*{\vc}[1]{\vcenter{\hbox{#1}}}

\begin{figure}[h!] \label{fig: Satake} 
\centering
\begin{minipage}{0.48\textwidth}
\centering
\adjustbox{scale=1}{
\begin{tikzcd}[row sep=0.2cm, column sep=0.6cm]
  & {\overset{1}{\circ}} & \cdots & {\overset{n-2}{\circ}} & {\overset{n-1}{\circ}} \\
  \vc{\scriptsize$0$}\ \circ &&&&& \circ\ \vc{\scriptsize$n$} \\
  & {\underset{2n-1}{\circ}} & \cdots & {\underset{n+2}{\circ}} & {\underset{n+1}{\circ}}
  \arrow[no head, from=1-2, to=1-3]
  \arrow["\tau"', color={rgb,255:red,214;green,92;blue,92}, bend left=20, <->, from=3-2, to=1-2]
  \arrow[no head, from=1-3, to=1-4]
  \arrow[no head, from=1-4, to=1-5]
  \arrow["\tau"', color={rgb,255:red,214;green,92;blue,92}, bend left=20, <->, from=3-4, to=1-4]
  \arrow[no head, from=1-5, to=2-6]
  \arrow["\tau"', color={rgb,255:red,214;green,92;blue,92}, bend left=20, <->, from=3-5, to=1-5]
  \arrow[no head, from=2-1, to=1-2]
  \arrow[no head, from=2-1, to=3-2]
  \arrow[no head, from=3-2, to=3-3]
  \arrow[no head, from=3-3, to=3-4]
  \arrow[no head, from=3-4, to=3-5]
  \arrow[no head, from=3-5, to=2-6]
\end{tikzcd}
}
\\[1ex]
(a)
\end{minipage}
\hfill
\begin{minipage}{0.48\textwidth}
\centering
\adjustbox{scale=1}{
\quad \quad
\begin{tikzcd}[row sep=0.2cm, column sep=0.6cm]
  & {\overset{1}{\circ}} & \cdots & {\overset{n-1}{\circ}} & {\overset{n}{\circ}} \\
  {\underset{0}{\circ}} \\
  & {\underset{2n}{\circ}} & \cdots & {\underset{n+2}{\circ}} & {\underset{n+1}{\circ}}
  \arrow[no head, from=1-2, to=1-3]
  \arrow["\tau"', color={rgb,255:red,214;green,92;blue,92}, bend left=20, <->, from=3-2, to=1-2]
  \arrow[no head, from=1-3, to=1-4]
  \arrow[no head, from=1-4, to=1-5]
  \arrow["\tau"', color={rgb,255:red,214;green,92;blue,92}, bend left=20, <->, from=3-4, to=1-4]
  \arrow["\tau"', shift left=3, color={rgb,255:red,214;green,92;blue,92}, bend left=20, <->, from=3-5, to=1-5]
  \arrow[no head, from=1-5, to=3-5]
  \arrow[no head, from=2-1, to=1-2]
  \arrow[no head, from=2-1, to=3-2]
  \arrow[no head, from=3-2, to=3-3]
  \arrow[no head, from=3-3, to=3-4]
  \arrow[no head, from=3-4, to=3-5]
\end{tikzcd}
}
\\[1ex]
(b)
\end{minipage}
\caption{(a) Satake diagram of type ${\mathsf{AIII}}_{2n-1}^{(\tau)}$ ($n \ge 2$); (b) Satake diagram of type ${\mathsf{AIII}}_{2n}^{(\tau)}$ ($n \ge 1$).}
\label{fig:Satake}
\end{figure}
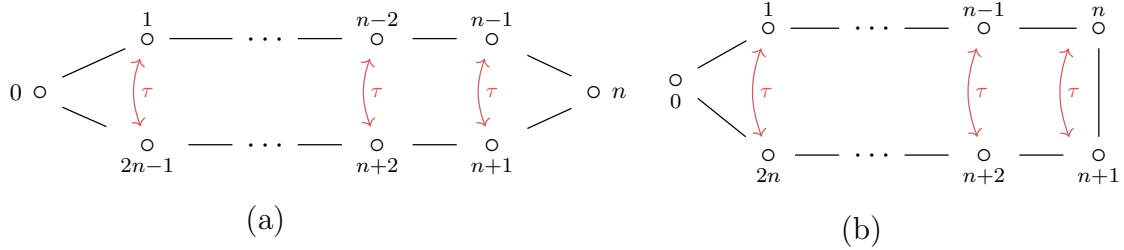

We first recall the definition of $\Ui$ in the Kolb--Letzter presentation. 

\Defi 
Let $\Ui$ be the algebra generated by $\Bg_i$ and commuting invertible elements $\KK_i$ $(i \in \indx)$ subject to relations 
\begin{align*} 
\KK_i B_j &= q^{-a_{ij} +a_{\tau(i),j}}B_j \KK_i, \\ 
\Serre_{ij}(B_i,B_j) &= 
\left\{
\begin{array}{r l l l l} 
0 & \mbox{if} & a_{\tau(i),j} \neq 2 & \& &  \tau(i) \neq i, \\ 
 (q-q\mi)\mi (\KK_{\tau(i)} - \KK_{i}) & \mbox{if} & a_{\tau(i),j} = 2 & \& & a_{i,\tau(i)} = 0, \\ 
- q[2] (\KK_i B_i + B_i \KK_{\tau(i)}) & \mbox{if} & a_{\tau(i),j} = 2 & \& & a_{i,\tau(i)} = -1, \\ 
0 & \mbox{if} & a_{ij} = 0 & \& & \tau(i) = i, \\ 
-q\mi \Bg_j \KK_i & \mbox{if} & a_{ij} = -1 & \& & \tau(i) = i, \\ 
- q\mi \KK_i [2]^2[\Bg_i, \Bg_j] & \mbox{if} & a_{ij} = -2. 
\end{array}
\right. 
\end{align*}
\enDefi

Given $\mu = \sum_{i \in \indx} d_i \alpha_i \in \Z \indx = \bigoplus_{i \in \indx} \Z \alpha_i$, set $\KK_\mu = \prod \KK_i^{d_i}$.

\subsection{The relative braid group action} 
\nc{\balpha}{\boldsymbol\alpha} 

We first recall the notion of a relative root system. Let $\mathcal{R}$ be the affine root system associated to the Dynkin diagram underlying a chosen Satake diagram. Given $\alpha \in \mathcal{R}$, let 
\[
\boldsymbol\alpha = \frac{\alpha + \tau(\alpha)}{2} \in \Q\mathcal{R}. 
\]
Then 
\[
\mathcal{R}^{\rm rel} = \{ \boldsymbol\alpha \mid \alpha \in \mathcal{R} \}
\]
is the \emph{relative root system} associated to the Satake diagram of the affine symmetric pair $(\widehat{\g},\widehat{\g}^{\omega \tau})$. In particular, $\{ \balpha_i \mid i \in \tindx \}$ forms a set of simple roots, and $\{ \balpha \mid \alpha \in \mathcal{R}^+ \}$ a set of positive roots for the relative root system.

The corresponding \emph{relative affine} (resp.\ finite) \emph{Weyl group} 
is the subgroup $W_\tau$ (resp.\ $W_{0,\tau}$) of the usual Weyl group $W$ 
generated by 
\[
r_i = \begin{cases}
s_i & \mbox{if} \ a_{i,\tau(i)} = 2, \\
s_i s_{\tau(i)} & \mbox{if} \ a_{i,\tau(i)} = 0, \\ 
s_i s_{\tau(i)} s_i  & \mbox{if} \ a_{i,\tau(i)} = -1.  
\end{cases}
\] 
for $i \in \tindx$ (resp.\ $i \in \indx_{0,\tau}$). It is a Coxeter group isomorphic to a Weyl group of type $\mathsf{A}^{(1)}_1$, $\mathsf{C}^{(1)}_n$ or $\mathsf{A}^{(2)}_{2n}$, as illustrated in Table \ref{Relative affine Weyl groups of types AIII and DI}.

\begin{table}[h]
\centering
\renewcommand{\arraystretch}{1.6}
\begin{tabular}{|c|c|c|}
\hline
\shortstack{Satake type} & 
\shortstack{Relative root system} & 
Dynkin diagram for relative root system \\
\hline

$\mathsf{AI}_{1}$ & 
$\mathsf{A}_{1}^{(1)}$ & 
\begin{tikzpicture}[scale=0.5, baseline=(current bounding box.center)] 
\draw (0, 0.1 cm) -- +(2 cm,0);
\draw (0, -0.1 cm) -- +(2 cm,0);
\draw[fill=white] (0 cm, 0 cm) circle (.25cm) node[below=4pt]{$0$};
\draw[fill=white] (2 cm, 0 cm) circle (.25cm) node[below=4pt]{$1$};
\end{tikzpicture} \\
\hline

$\mathsf{AIII}^{(\tau)}_{2n-1}\ (n \ge 2)$ & 
$\mathsf{C}_n^{(1)}$ & 

\begin{tikzpicture}[scale=0.5, baseline=(current bounding box.center)]
\draw (0, 0.1 cm) -- +(2 cm,0);
\draw (0, -0.1 cm) -- +(2 cm,0);
\draw[shift={(1.2, 0)}, rotate=0] (135 : 0.45cm) -- (0,0) -- (-135 : 0.45cm);
{
\pgftransformxshift{2 cm}
\draw (0 cm,0) -- (2 cm,0);
\draw (2 cm,0) -- (4 cm,0);
\draw[dashed] (4 cm,0) -- (6 cm,0);
\draw (6 cm,0) -- (8 cm,0);
\draw (8 cm, 0.1 cm) -- +(2 cm,0);
\draw (8 cm, -0.1 cm) -- +(2 cm,0);
\draw[shift={(8.8, 0)}, rotate=180] (135 : 0.45cm) -- (0,0) -- (-135 : 0.45cm);
\draw[fill=white] (0 cm, 0 cm) circle (.25cm) node[below=4pt]{$1$};
\draw[fill=white] (2 cm, 0 cm) circle (.25cm) node[below=4pt]{$2$};
\draw[fill=white] (4 cm, 0 cm) circle (.25cm) node[below=4pt]{};
\draw[fill=white] (6 cm, 0 cm) circle (.25cm) node[below=4pt]{};
\draw[fill=white] (8 cm, 0 cm) circle (.25cm) node[below=4pt]{$n-1$};
\draw[fill=white] (10 cm, 0 cm) circle (.25cm) node[below=4pt]{$n$};
}
\draw[fill=white] (0 cm, 0 cm) circle (.25cm) node[below=4pt]{$0$};
\end{tikzpicture}  \\
\hline

$\mathsf{AIII}^{(\tau)}_{2}$ & 
$\mathsf{A}_{2}^{(2)}$ & 
\begin{tikzpicture}[scale=0.5, baseline=(current bounding box.center)] 
\draw (0, 0.05 cm) -- +(2 cm,0);
\draw (0, -0.05 cm) -- +(2 cm,0);
\draw (0, 0.15 cm) -- +(2 cm,0);
\draw (0, -0.15 cm) -- +(2 cm,0);
\draw[shift={(0.8, 0)}, rotate=180] (135 : 0.45cm) -- (0,0) -- (-135 : 0.45cm);
\draw[fill=white] (0 cm, 0 cm) circle (.25cm) node[below=4pt]{$0$};
\draw[fill=white] (2 cm, 0 cm) circle (.25cm) node[below=4pt]{$1$};
\end{tikzpicture} \\
\hline

$\mathsf{AIII}^{(\tau)}_{2n} \ (n \ge 2)$ & 
$\mathsf{A}_{2n}^{(2)}$ & 
\begin{tikzpicture}[scale=0.5, baseline=(current bounding box.center)]

\draw (0, 0.1 cm) -- +(2 cm,0);
\draw (0, -0.1 cm) -- +(2 cm,0);
\draw[shift={(1.2, 0)}, rotate=180] (135 : 0.45cm) -- (0,0) -- (-135 : 0.45cm);

{
\pgftransformxshift{2 cm}

\draw (0 cm,0) -- (2 cm,0);
\draw (2 cm,0) -- (4 cm,0);
\draw[dashed] (4 cm,0) -- (6 cm,0);
\draw (6 cm,0) -- (8 cm,0);

\draw (8 cm, 0.1 cm) -- +(2 cm,0);
\draw (8 cm, -0.1 cm) -- +(2 cm,0);
\draw[shift={(8.8, 0)}, rotate=180] (135 : 0.45cm) -- (0,0) -- (-135 : 0.45cm);

\draw[fill=white] (0 cm, 0 cm) circle (.25cm) node[below=4pt]{$1$};
\draw[fill=white] (2 cm, 0 cm) circle (.25cm) node[below=4pt]{$2$};
\draw[fill=white] (4 cm, 0 cm) circle (.25cm) node[below=4pt]{};
\draw[fill=white] (6 cm, 0 cm) circle (.25cm) node[below=4pt]{};
\draw[fill=white] (8 cm, 0 cm) circle (.25cm) node[below=4pt]{$n-1$};
\draw[fill=white] (10 cm, 0 cm) circle (.25cm) node[below=4pt]{$n$};
}

\draw[fill=white] (0 cm, 0 cm) circle (.25cm) node[below=4pt]{$0$};
\end{tikzpicture} \\
\hline 
\end{tabular} \\[2ex]
\caption{Correspondence between Satake types and relative root systems.}
\label{Relative affine Weyl groups of types AIII and DI}
\end{table}

The \emph{extended relative affine Weyl group} $\widetilde{W}_\tau$ is the semi-direct product 
\[
\widetilde{W}_\tau = W_{0,\tau} \ltimes P^\tau \cong W_\tau \rtimes \Lambda^\tau, 
\]
where $P^\tau = \{ x \in P \mid \tau(x) = x\}$, and $\Lambda^\tau = P^\tau/Q^\tau$ can be identified with the group of automorphisms of the Dynkin diagram associated to the relative root system. 
 Equivalently, 
\[
\widetilde{W}_\tau = \{ w \in \widetilde{W} \mid w\tau = \tau w\}. 
\]
While discussing particular types, we will use the notation $\widetilde{W}_{\mathsf{AIII}_{N}^{(\tau)}}$ to denote the corresponding extended relative affine Weyl group. Let
\[
\varpi_i = \begin{cases}
\omega_i & \mbox{if} \ a_{i,\tau(i)} = 2, \\
\omega_i + \omega_{\tau(i)} & \mbox{if} \ a_{i,\tau(i)} \neq 2, 
\end{cases}
\] 
for $i \in \indx_{0,\tau}$. 
Let $\widetilde{\mathcal{B}}_\tau$ be the braid group associated to $\widetilde{W}_\tau$. 

\Thm 
There exists an action of $\widetilde{\mathcal{B}}_\tau$ on $\Ui$ by algebra automorphisms 
\[ r_i \mapsto \Tbr_i \quad (i \in \indx_\tau), \qquad \lambda \mapsto \Tbr_\lambda \quad  (\lambda \in \Lambda^\tau) \] 
given below. 
\be
\item 
For each $i \in \tindx$ with $\tau(i) = i$, there exists an automorphism $\Tbr_i$ of $\Ui$ such that $\Tbr_i(\KK_\mu) = \KK_{r_i\mu}$ and 
\[
\Tbr_i(B_j) = 
\left\{ \begin{array}{l l}
\KK_j\mi B_j & \mbox{if } i=j, \\[2pt]
B_j & \mbox{if } a_{ij} = 0, \\[2pt]
[B_j,B_i]_q  & \mbox{if } a_{ij} = -1, \\[2pt]
\end{array} \right. 
\] 
for $\mu \in \Z\indx$ and $j \in \indx$. 
\item For each $i \in \tindx$ with $a_{i,\tau(i)} = 0$, there exists an automorphism $\Tbr_i$ of $\Ui$ such that $\Tbr_i(\KK_j) = (-q)^{-a_{ij}-a_{\tau(i),j}} \KK_{j}\KK_i^{-a_{ij}}\KK_{\tau(i)}^{-a_{\tau(i),j}}$ and 
\[
\Tbr_i(B_j) = 
\left\{ \begin{array}{l l}
-\KK_j\mi B_{\tau(j)} & \mbox{if } i=j \mbox{ or } \tau(i) = j, \\[2pt]
[B_j,B_i]_q   & \mbox{if } a_{ij} = -1 \ \& \ a_{\tau(i),j}=0, \\[2pt] 
[B_j,B_{\tau(i)}]_q   & \mbox{if } a_{ij} = 0 \ \& \ a_{\tau(i),j}=-1, \\[2pt] 
[[B_j,B_{i}]_q,B_{\tau(i)}]_q - qB_j \KK_i   & \mbox{if } a_{ij} = -1 \ \& \ a_{\tau(i),j}=-1, \\[2pt]
B_j & \mbox{otherwise}, \\[2pt]
\end{array} \right. 
\] 
for $\mu \in \Z\indx$ and $j \in \indx$. 
\item For $i \in \tindx$ with $a_{i,\tau(i)} = -1$, there exists an automorphism $\Tbr_i$ of $\Ui$ such that $\Tbr_i(\KK_j) = q^{-a_{ij}-a_{\tau(i),j}} \KK_{j}(\KK_i\KK_{\tau(i)})^{-a_{ij}-a_{\tau(i),j}}$ and 
\[
\Tbr_i(B_j) = 
\left\{ \begin{array}{l l}
-q^{-2}B_j \KK_{\tau(j)}\mi  & \mbox{if } i=j \mbox{ or } \tau(i) = j, \\[2pt]
[[B_j,B_i]_q, B_{\tau(i)}]_q - \KK_i B_j   & \mbox{if } a_{ij} = -1 \ \& \ a_{\tau(i),j}=0, \\[2pt] 
[[B_j,B_{\tau(i)}]_q, B_i]_q - \KK_{\tau(i)}B_j   & \mbox{if } a_{ij} = 0 \ \& \ a_{\tau(i),j}=-1, \\[2pt] 
q \big[ [[B_j,B_{i}]_q,B_{\tau(i)}], [B_{\tau(i)}, B_i]_q \big] & \\
\qquad -[B_j, [B_{\tau(i)}, B_i]_{q^3}] \KK_i - qB_j \KK_i\KK_{\tau(i)}   & \mbox{if } a_{ij} = -1 \ \& \ a_{\tau(i),j}=-1, \\[2pt]
B_j & \mbox{otherwise}, \\[2pt]
\end{array} \right. 
\] 
for $\mu \in \Z\indx$ and $j \in \indx$. 
\item 
For each $\lambda \in \Lambda^\tau$, there is also an automorphism $\Tbr_\lambda$ such that $\Tbr_\lambda(B_i) = B_{\lambda(i)}$ and $\Tbr_\lambda(\KK_i) = \KK_{\lambda(i)}$.   
\ee
\enthm 

\Proof
See \cite[Theorem 2.6]{LWZ-quasi} and \cite[Theorem 2.2]{LPWZ} (based on \cite{ZhangBr}). 
\enproof

Given a reduced expression $w = \lambda r_{i_1} \cdots r_{i_k} \in \widetilde{W}_\tau$, we set $\Tbr_{w} = \Tbr_\lambda \Tbr_{i_1} \cdots \Tbr_{i_k}$. 
We will need the following lemma. 

\Lem \label{lem: br act on wt}
We have: $T_w(e^{\pm}_i) = e^{\pm}_{wi}$ ($w \in \widetilde{W}$), and $\Tbr_w(B_i) = B_{wi}$ ($w \in \widetilde{W}_\tau$), for $i \in \indx$ such that $wi \in \indx$. 
\enlem

\Proof
See \cite{lusztig-94} and \cite[Lemma 2.9]{LWZ-quasi}. 
\enproof

\subsection{The Lu--Wang presentation} 

Next, let us recall the ``Drinfeld-type" presentation of $\Ui$ from \cite{LWZ-quasi, LPWZ}. 
We will use the following shorthand notations: 
\begin{align*}
S(r_1, r_2 \mid s; i,j) =& \ \Ag_{i,r_1} \Ag_{i,r_2} \Ag_{j,s} - [2] \Ag_{i,r_1} \Ag_{j,s} \Ag_{i,r_2} + \Ag_{j,s}\Ag_{i,r_1} \Ag_{i,r_2}, \\
\mathbb{S}(r_1, r_2 \mid s; i,j) =& \ S(r_1, r_2 \mid s; i,j) + S(r_2, r_1 \mid s; i,j), \\ 
R(r_1, r_2 \mid s; i,j) =& \ \KK_i\Ck^{r_1}\big( - \sum_{p \geq 0} q^{2p} [2] [\Thg_{i, r_2-r_1-2p-1}, \Ag_{j,s-1}]_{q^{-2}}  \Ck^{p+1} \\ 
& \ - \sum_{p \geq 1} q^{2p-1} [2] [\Ag_{j,s}, \Thg_{i, r_2-r_1-2p}]_{q^{-2}} \Ck^{p} - [\Ag_{j,s}, \Thg_{i,r_2-r_1}]_{q^{-2}}\big) \\
&\quad \mbox{if}   \quad i = \tau(i), \\ 
R(r_1, r_2 \mid s; i, \tau(i)) =& \ [2]\sum_{p \geq 0} q^{2p} [\Thg_{\tau(i), s-r_2-p}\KK_i - q\Thg_{\tau(i), s-r_2-p-2}\CCC\KK_i, \Ag_{i,r_1-p}]_{q^{-4p-1}}  \Ck^{r_2+p} \\ 
& \ + q[2] \sum_{p \geq 1} q^{2p} [\Ag_{i,r_1+p+1}, \Thg_{i, r_2-s-p+1}\KK_{\tau(i)} - q\Theta_{i,r_2-s-p-1}\CCC\KK_{\tau(i)}]_{q^{-4p-3}} \Ck^{s-1}  \\
&\quad \mbox{if}   \quad a_{i,\tau(i)} = -1, \\ 
\mathbb{R}(r_1, r_2 \mid s; i,j) =& \ R(r_1, r_2 \mid s; i,j) + R(r_2, r_1 \mid s; i,j). 
\end{align*}

\Defi \label{defi: LW pres gen}
Let $\DrUi$ be the algebra generated by $\KK_i^{\pm1}$, $\Hg_{i,m}$ and $\Ag_{i,r}$, where $m\geq1$, $r\in\Z$ and $i \in \indx_0$, as well as a central element $\CCC$, subject to the following relations:
\begin{align} 
\KK_i A_{j,r} &= q^{a_{\tau(i),j} - a_{ij}} A_{j,r} \KK_i, \\
[\Hg_{i,m},\Hg_{j,l}] = [\Hg_{i,m},\KK_{j}] &= [\KK_{i}, \KK_{j}] = 0, \label{eq: grel1} \\
[\Hg_{i,m}, \Ag_{j,r}] &= \textstyle \frac{[m \cdot a_{ij}]}{m} \Ag_{j,r+m} - \frac{[m \cdot a_{\tau(i),j}]}{m} \Ag_{j,r-m}\CCC^m, \label{eq: grel2} \\ 
[\Ag_{i,r}, \Ag_{\tau(i),s}] &= \KK_{\tau(i)} \CCC^s\Theta_{i,r-s} - \KK_{i} \CCC^r\Theta_{\tau(i),s-r} \qquad \mbox{if} \quad a_{i,\tau(i)} = 0, \label{eq: grel3} \\ 
[A_{i,r}, A_{j,s}] &= 0 \qquad \mbox{if} \quad a_{ij} = 0 \ \& \ j \neq \tau(i), \\
[\Ag_{i,r}, \Ag_{j,s+1}]_{q^{-a_{ij}}}   &= q^{-a_{ij}} [\Ag_{i,r+1}, \Ag_{j,s}]_{q^{a_{ij}}}  \quad \mbox{if} \quad j \neq \tau(i), \label{eq: grel4} \\ 
\label{eq: grel5}
[\Ag_{i,r}, \Ag_{i,s+1}]_{q^{-2}}  -q^{-2} [\Ag_{i,r+1}, \Ag_{i,s}]_{q^{2}}
&= q^{-2}\KK_i\CCC^r\Thg_{i,s-r+1} - q^{-4}\KK_i\CCC^{r+1}\Thg_{i,s-r-1}  \\ \notag
&\quad  + q^{-2}\KK_i\CCC^s\Thg_{i,r-s+1}  -q^{-4}\KK_i \CCC^{s+1}\Thg_{i,r-s-1} \\ 
&\quad \mbox{if} \quad i = \tau(i),\\ 
\label{eq: grel5a}
[\Ag_{i,r}, \Ag_{\tau(i),s+1}]_{q}  -q [\Ag_{i,r+1}, \Ag_{\tau(i),s}]_{q\mi}
&= -\KK_i\CCC^r\Thg_{\tau(i),s-r+1} + q\KK_i\CCC^{r+1}\Thg_{\tau(i),s-r-1}  \\ \notag
&\quad  - \KK_{\tau(i)}\CCC^s\Thg_{i,r-s+1}  +q\KK_{\tau(i)} \CCC^{s+1}\Thg_{i,r-s-1} \\ 
&\quad \mbox{if} \quad a_{i,\tau(i)} = -1,\\
\mathbb{S}(r_1, r_2 \mid s; i,j) &= \mathbb{R}(r_1, r_2 \mid s; i,j) \qquad \mbox{if} \quad a_{ij} = -1 \ \& \ \tau(i) = i, \label{eq: grel6} \\ 
\mathbb{S}(r_1, r_2 \mid s; i, \tau(i)) &= \mathbb{R}(r_1, r_2 \mid s; i, \tau(i)) \qquad \mbox{if} \quad a_{i,\tau(i)} = -1, \\ 
\mathbb{S}(r_1, r_2 \mid s; i,j) &= 0  \qquad \mbox{if} \quad a_{ij} = -1 \ \& \ j \neq \tau(i) \neq i, 
\end{align}
where $l,m\geq1$; $r,s, r_1, r_2\in \Z$ and 
\[
1+ \sum_{m=1}^\infty (q-q^{-1})\Thg_{i,m} z^m  =  \exp \left( (q-q^{-1}) \sum_{m=1}^\infty \Hg_{i,m} z^m \right).
\]
By convention, $\Thg_{i,0} = (q-q^{-1})^{-1}$ and $\Thg_{i,m} = 0$ for $m\leq -1$. Moreover, we call the commutative subalgebra
\begin{equation} \label{eq: HT Cartans}
\mathcal{H} = \langle \Theta_{i,r} \mid i \in \indx_0, r \in \Z_{\geq0} \rangle. 
\end{equation}
of $\Ui$ the \emph{Lu--Wang Cartan subalgebra}. 
\enDefi 

We choose a sign function 
$
o(\cdot) \colon \indx_0 \to \{\pm 1\}$ such that $o(i) = -o(j)$ if $a_{ij} = -1$. 
By \cite[Theorem 4.5]{LWZ-quasi} and \cite[Theorem 4.7]{LPWZ}, there is an algebra isomorphism 
\eq
\label{eq: droq->oq}
\DrUi \isoto{} \Ui, \qquad A_{i,r} \mapsto o(i)^r \Tbr_{\varpi_i}^{-r}(B_i), \quad
\KK_i \mapsto \KK_i, \quad  \Ck \mapsto (-q)^{N-1} \KK_\delta, 
\eneq
for $i \in \indx_0$. 

We will use the following notation for generating series: 
\[
\avar{i} = \sum_{r \geq 0} A_{i,r} z^r, \qquad \Thgs_i(z) = \Theta_{i,r} z^r. 
\]
Abbreviate $\overline{H}_{i,1} = [2]\mi H_{i,1}$ and $\rho = (q-q\mi)$. 
We will also use the following renormalized series: 
\[
\thvar{i} = \frac{\rho(1 - q^{-a_{i,\tau(i)}} \CCC z^2)}{1 - \CCC z^2} \boldsymbol\Theta_i(z) = \begin{cases}
\frac{\rho(1-q^{-2}\CCC z^2)}{1 - \CCC z^2} \Thgs_i(z) & \mbox{if } a_{i,\tau(i)} = 2, \\[3pt]
\rho \Thgs_i(z) & \mbox{if } a_{i, \tau(i)} = 0, \\[3pt]
\frac{\rho(1-q\CCC z^2)}{1 - \CCC z^2} \Thgs_i(z) & \mbox{if } a_{i,\tau(i)} = -1.
\end{cases} 
\]

\Lem \label{lem: A Th as series}
The series $\avar{i}$ and $\thvar{i}$ are expansions of the following rational functions: 
\[
\avar{i} = 
\begin{cases} 
(1 - \Ad_{\overline{H}_{i,1}} z)\mi (A_{i,0}) & \mbox{if } a_{i,\tau(i)} = 0, \\ 
(1 - \Ad_{\overline{H}_{i,1}} z - \CCC z^2)\mi (A_{i,0} + \CCC z A_{i,-1}) & \mbox{if } a_{i,\tau(i)} = 2, \\
([2] - \Ad_{H_{i,1}} z + \CCC z^2)\mi ([2]A_{i,0} - \CCC z A_{i,-1}) & \mbox{if } a_{i,\tau(i)} = -1, 
\end{cases}
\]
\[
\thvar{i} = 
\begin{cases} 
\KK_i \KK_{\tau(i)}\mi + \rho\KK_{\tau(i)}\mi [\avar{i}, A_{\tau(i),0}] & \mbox{if } a_{i,\tau(i)} = 0, \\[3pt]
1 + \frac{q^2\rho\CCC \KK_i\mi z^2}{1 - \CCC z^2}\left(z\mi[A_{i,-1}, \avar{i}]_{q^{-2}} - q^{-2}[A_{i,0}, \avar{i}]_{q^2} \right) & \mbox{if } a_{i,\tau(i)} = 2, \\[3pt]
\frac{-\rho\CCC \KK_{\tau(i)}\mi z^2}{1 - \CCC z^2} \Big( z\mi[A_{\tau(i),-1}, \avar{i}]_{q} - q[A_{\tau(i),0}, \avar{i}]_{q\mi}  \\[3pt]
\quad - o(i) q \CCC\mi \KK_{\tau(i)} \Tbr_{\boldsymbol\theta_n}\mi(B_0) \KK_{i} z\mi 
- \rho\mi (\CCC\mi \KK_{\tau(i)} z^{-2} - \KK_{i}) \Big)   & \mbox{if } a_{i,\tau(i)} = -1, 
\end{cases}
\]
with $\boldsymbol\theta_n$ as in \eqref{eq: thetar}, and $\Ad_y = [y, - ]$. 
\enlem

\Proof
The formulae follow rather straightforwardly from the presentation in Definition~\ref{defi: LW pres gen} (c.f. \cite[Proposition 3.1]{Przez-23}). We will limit ourselves to explicitly deriving the formula for $\thvar{i}$ in the case $a_{i,\tau(i)} = -1$. The Lu--Wang presentation implies that: 
\[
[\Ag_{i,-1}, \Ag_{\tau(i),s+1}]_{q}  -q [\Ag_{i,0}, \Ag_{\tau(i),s}]_{q\mi} =
\begin{cases}
\KK_i (- \CCC\mi \Theta_{\tau(i), s+2} + q \Theta_{\tau(i),s}) & \mbox{if } s \geq 1, \\ 
-\KK_i\CCC\mi \Theta_{\tau(i), 2} + \rho\mi(q\KK_i -\KK_{\tau(i)}) & \mbox{if } s =0, \\ 
-\KK_i\CCC\mi \Theta_{\tau(i), 1} - \KK_{\tau(i)}\CCC\mi \Theta_{i, 1} & \mbox{if } s =-1. 
\end{cases}
\]
Therefore, 
\begin{align*}
Z:= z\mi[A_{\tau(i),-1}, \ &\avar{i}]_{q} - q[A_{\tau(i),0}, \avar{i}]_{q\mi} = 
[A_{i,-1}, A_{\tau(i),0} ]_q z\mi \\ 
&\quad + \sum_{s \geq 0} \left( [\Ag_{i,r}, \Ag_{\tau(i),s+1}]_{q}  -q [\Ag_{i,r+1}, \Ag_{\tau(i),s}]_{q\mi} \right) z^s \\ 
&= \sum_{s \geq 1} \KK_i (- \CCC\mi \Theta_{\tau(i), s+2} + q \Theta_{\tau(i),s}) z^s \\ 
&\quad -\KK_i\CCC\mi \Theta_{\tau(i), 2} + \rho\mi(q\KK_i -\KK_{\tau(i)}) + [A_{i,-1}, A_{\tau(i),0} ]_q z\mi \\
&= q \KK_i \Theta_{\tau(i),s} z + \sum_{s \geq 2} \KK_i (q - \CCC\mi z^{-2}) \Theta_{\tau(i),s} z^s \\
&\quad + \rho\mi(q\KK_i -\KK_{\tau(i)}) + [A_{i,-1}, A_{\tau(i),0} ]_q z\mi \\ 
&= \KK_i (q - \CCC\mi z^{-2}) \boldsymbol\Theta_{\tau(i)}(z) + \CCC\mi \rho\mi \KK_i z^{-2} - \rho\mi \KK_{\tau(i)} \\
&\quad + \CCC\mi z\mi \KK_i \Theta_{\tau(i),1} + [A_{i,-1}, A_{\tau(i),0} ]_q z\mi. 
\end{align*}
By \cite[(4.7)]{LPWZ}, 
\[
\Theta_{\tau(i), 1} = - [A_{i,-1}, A_{\tau(i)}]_q \CCC \KK_i\mi + o(\tau(i)) q \Tbr_{\boldsymbol\theta_n}\mi(B_0)\KK_{\tau(i)}. 
\]
It follows that 
\begin{align} \label{eq: case3 theta formula}
\KK_i (q - \CCC\mi z^{-2}) \boldsymbol\Theta_{\tau(i)}(z) &= Z - \rho\mi (\CCC\mi \KK_i z^{-2} - \KK_{\tau(i)}) \\
&\quad - o(\tau(i)) q \CCC\mi \KK_i \Tbr_{\boldsymbol\theta_n}\mi(B_0) \KK_{\tau(i)} z\mi. 
\end{align} 
The desired formula now follows by interchanging the indices $i \leftrightarrow \tau(i)$. 
\enproof

\subsection{Coideal structures} 

Assume that $\mathbf{c} = (c_i) \in (\C^\times)^{\indx}$ satisfies $c_i = c_{\tau(i)}$
There exists an injective algebra homomorphism 
\eq
\label{eq: Kolb emb}
\eta \colon \Ui \monoto \Uu', \quad \Bg_i \mapsto F_i +  E_{\tau(i)} \Kb_i, \quad 
\KK_i \mapsto
\begin{cases}
\Ka_i \Kb_{\tau(i)} & \mbox{if} \ \tau(i) \neq i, \\
- q^2 \Ka_i \Kb_i & \mbox{if} \ \tau(i) = i, 
\end{cases}
\eneq 
for $i \in \indx$. 
This homomorphism makes $\Ui$ into a coideal subalgebra of $\Uu'$, with the coideal structure explicitly given by 
\eq \label{eq: coproduct explicit} 
\Delta(B_i) = 1 \otimes \eta(B_i) + \eta(B_i) \otimes \Kb_i, \qquad \Delta(\KK_i) = \KK_i \otimes \KK_i. 
\eneq  

Let $\auic$ be the quotient of $\Ui$ by the two-sided ideal generated by 
\[
\KK_i + q^2 c_i \quad (\tau(i) = i), \qquad \KK_i \KK_{\tau(i)} - c_i^2 \quad (\tau(i) \neq i). 
\] 
The monomorphism \eqref{eq: Kolb emb} descends to the embedding 
\eq \label{eq: Kolb emb 2}
\eta_{\mathbf{c}} \colon \auic \hookrightarrow \Uq_{\mathbf{c}} \cong \Uq, \qquad 
B_i \mapsto F_i + c_i E_{\tau(i)} K_i\mi, \quad \KK_i \mapsto 
\begin{cases}
c_i K_i K_{\tau(i)}\mi & \mbox{if} \ \tau(i) \neq i, \\
- q^2 c_i & \mbox{if} \ \tau(i) = i, 
\end{cases} 
\eneq
making $\auic$ into a quantum symmetric pair coideal subalgebra of $\Uq$ (c.f. \cite{Kolb-14}). 
In particular,
\[
\CCC \mapsto 
\begin{cases}
q^{N+3} c_\delta & \mbox{if } N \mbox{ is odd}, \\
q^{N+1} c_\delta & \mbox{if } N \mbox{ is even}. 
\end{cases}
\]

\subsection{Small rank subalgebras of $\Uu'$}

\nc{\Uii}{\Ui_{[i]}} 

\nc{\Uqi}{\Uu'_{[i]}} 
\nc{\Uui}{\Uu_{[i]}} 
\nc{\Usl}{U_q(L\mathfrak{sl}_2)}
\nc{\Uusl}{\widetilde{U}_q(L\mathfrak{sl}_2)}

For $i \in \indx_0$, set $\omega'_i = \omega_i s_i$, and let $\Uqi$ be the subalgebra of $\Uu'$ generated by 
\eq \label{eq: rank 1 U gen}
E_i, \ F_i, \ \Ka^{\pm 1}_i, \ (\Kb_i)^{\pm 1}, \ T_{\omega'_i}(E_i), \ T_{\omega'_i}(F_i), \ T_{\omega'_i}(\Ka_i^{\pm 1}), \ T_{\omega'_i}((\Kb_i)^{\pm 1}). 
\eneq
By (a slight modification of) \cite[Proposition 3.8]{beck-94}, there is an algebra isomorphism $\iota_i \colon \Uu'(\widehat{\mathfrak{sl}}_2) \to \Uqi$ sending 
\begin{align*}
E_1 &\mapsto E_i, \ F_1 \mapsto F_i, \ \Ka_1^{\pm1} \mapsto \Ka_i^{\pm1}, \ (\Kb_1)^{\pm1} \mapsto (\Kb_i)^{\pm1}, \\ 
E_0 &\mapsto T_{\omega'_i}(E_i), \ F_0 \mapsto T_{\omega'_i}(F_i), \ \Ka_0^{\pm1} \mapsto T_{\omega'_i}(\Ka_i^{\pm1}), \ (\Kb_0)^{\pm1} \mapsto T_{\omega'_i}((\Kb_i)^{\pm1}). 
\end{align*}
In terms of Drinfeld's new presentation, $\iota_i$ sends $x^{\pm}_{r} \mapsto x^{\pm}_{i,r}$ and $\phi^{\pm}_s \mapsto \phi^{\pm}_{i,s}$.  

Next, assume that $N \geq 2$ and $1 \leq j \leq N-1$. Let $\Uu'_{[j,j+1]}$ be the subalgebra of $\Uu'$ generated by \eqref{eq: rank 1 U gen} for $i \in \{j, j+1\}$. The isomorphism $\iota_i$ generalizes to an isomorphism $\iota_{j,j+1} \colon  \Uu'(\widehat{\mathfrak{sl}}_3) \to \Uu'_{[j,j+1]}$. We will only needs its description in terms of the Drinfeld presentation, namely, $\iota_{j,j+1}$ sends $x^{\pm}_{k,r} \mapsto x^{\pm}_{j+k-1,r}$ and $\phi^{\pm}_{k,s} \mapsto \phi^{\pm}_{j+k-1,s}$ for $k \in \{1,2\}$.

\subsection{Rank one subalgebras of $\Ui$} 

For $i \in \indx_0$, let $\varpi'_i = \varpi_i r_i$. 
The affine rank one subalgebra of $\Ui$ associated with $i \in \tindx$ is the subalgebra 
\[
\Uii =
\begin{cases}
\langle B_j, \KK_j^{\pm1}, \Tbr_{\varpi'_i}(B_i), \CCC^{\pm1} \mid j \in \{i, \tau(i) \} \rangle & \mbox{if } a_{i,\tau(i)} \in \{0, 2\}, \\[3pt] 
\langle B_i, B_{\tau(i)}, \KK_i^{\pm1}, \KK_{\tau(i)}^{\pm1}, \Tbr_{\boldsymbol\theta_n}\mi(B_0), \Tbr_{\boldsymbol\theta_n}(\KK_0)^{\pm1} \rangle & \mbox{if } a_{i,\tau(i)} = -1. 
\end{cases}
\]
The second case occurs only if $N= 2n$, and we may then take $i = n$, $\tau(i) = n+1$ as $\Ui_{[n]} = \Ui_{[n+1]}$. 
Here $\boldsymbol\theta_n$ is the longest element in $\langle r_i \mid 1 \leq i < n \rangle$, with a reduced expression given by 
\eq \label{eq: thetar}
\boldsymbol\theta_n = 
r_{n-1} (r_{n-2} r_{n-1}) \cdots (r_1 \cdots r_{n-1}). 
\eneq
By \cite[(3.2)]{LPWZ}, one has 
\eq \label{eq: T theta b0}
\Tbr_{\boldsymbol\theta_n}\mi(B_0) = \Tbr\mi_{n-1} \cdots \Tbr\mi_2 \Tbr\mi_1(B_0). 
\eneq

Altogether, there are three types of rank one subalgebras of $\Ui$, depending on whether $a_{i, \tau(i)}$ equals $2,0$ or $-1$. 

\subsubsection{Case $a_{i,\tau(i)} = 2$: the $q$-Onsager algebra} 

The $q$-Onsager algebra is a split quantum symmetric pair coideal subalgebra of $U_q(\widehat{\mathfrak{sl}}_2)$. It was first introduced by Terwilliger in \cite{Ter01}.  
As in \cite{lu-wang-21}, we consider a central extension of the $q$-Onsager algebra called the \emph{universal} $q$-Onsager algebra. It has the following definition by generators and relations. 

\Defi
The universal $q$-Onsager algebra $\Oq$ is the algebra generated by $B_0, B_1$ and invertible central elements $\KK_0, \KK_1$ subject to relations:
\eq
\label{eq: Ons rel}
\sum_{r=0}^3 (-1)^r \sqbin{3}{r} \Bg_i^{3-r}\Bg_j\Bg_i^r = - q\mi \KK_i [2]^2[\Bg_i, \Bg_j] \qquad (0 \leq i \neq j \leq 1).
\eneq 
\enDefi

By \cite[Proposition 3.3]{LWZ-quasi}, there is an algebra isomorphism $\iota_i \colon \Oq \to \Uii$ sending 
\eq \label{eq: rank 1 q onsager}
B_1 \mapsto B_i, \quad B_0 \mapsto \Tbr_{\varpi'_i}(B_i), \quad \KK_1 \mapsto \KK_i, \quad \KK_0 \mapsto \KK_\delta \KK_i^{-1}. 
\eneq

\subsubsection{Case $a_{i,\tau(i)} = 0$}

Let $\overline{U}_q(\widehat{\mathfrak{sl}}_2)$ be the quotient of $\Uu'(\widehat{\mathfrak{sl}}_2)$ by the ideal generated by the central element $\widetilde{K}_\delta - \widetilde{K}'_\delta$. It is a variant of the Drinfeld double of quantum $\widehat{\mathfrak{sl}}_2$. 
By \cite[Proposition 3.7]{LWZ-quasi}, there is an algebra isomorphism $\iota_i \colon \overline{U}_q(\widehat{\mathfrak{sl}}_2) \to \Uii$ sending 
\begin{align} \label{eq: rank 1 sl2}
F_1 &\mapsto B_i, &\quad F_0 &\mapsto \Tbr_{\varpi'_i}(B_i), &\quad E_1 &\mapsto B_{\tau(i)}, &\quad E_0 &\mapsto \Tbr_{\varpi'_i}(B_{\tau(i)}), \\ 
\Ka_1 &\mapsto \KK_i, &\quad \Ka_0 &\mapsto \KK_\delta \KK_i^{-1}, &\quad
\Kb_1 &\mapsto \KK_{\tau(i)}, &\quad \Kb_0 &\mapsto \KK_\delta \KK_{\tau(i)}^{-1}. 
\end{align}

\subsubsection{Case $a_{i,\tau(i)} = -1$}

Let $\Ui(\widehat{\mathfrak{sl}}_3, \tau)$ be the quasi-split affine $\imath$-quantum group of type $\mathsf{AIII}_2^{(\tau)}$. It has the following presentation by generators and relations. 

\Defi
The algebra $\Ui(\widehat{\mathfrak{sl}}_3, \tau)$ is generated by $B_0, B_1, B_2$, and invertible commuting elements $\KK_0, \KK_1, \KK_2$, subject to the relations: 
\begin{align*}
\KK_i B_j &= q^{a_{\tau(i),j} - a_{ij}} B_j \KK_i, \\
\Serre_{ij}(B_i, B_j) &= 
\begin{cases}
-q[2] (\KK_i B_i + B_i \KK_j) & \mbox{if } (i,j) \in \{ (1,2), (2,1) \}, \\ 
0 & \mbox{if } j=0 \ \& \ i \in \{1,2\}, \\ 
-q\mi \KK_0 B_j & \mbox{if } i=0 \ \& \ j \in \{1,2\}. 
\end{cases}
\end{align*}
\enDefi

By \cite[Proposition 3.4]{LPWZ}, there is an algebra isomorphism $\iota_n \colon \Ui(\widehat{\mathfrak{sl}}_3, \tau) \to \Ui_{[n]}$ sending 
\eq \label{eq: rank 1 qs}
B_1 \mapsto B_n, \quad B_2 \mapsto B_{n+1}, \quad B_0 \mapsto \Tbr\mi_{\boldsymbol\theta_n}(B_0), \quad \KK_1 \mapsto \KK_n, \quad \KK_2 \mapsto \KK_{n+1}, \quad \KK_0 \mapsto \Tbr_{\boldsymbol\theta_n}(\KK_0). 
\eneq 
The isomorphism is compatible with the loop presentations, in the sense that it sends 
$A_{i,r} \mapsto A_{n+i-1,r}$ and $\Theta_{i,s} \mapsto \Theta_{n+i-1,s}$ for $i \in \{1,2\}$.  


\section{Factorization in rank one} 
\label{sec: new rank 1}

A key ingredient of our argument in a detailed analysis of rank one cases. More precisely, we require an expression for the series $\thvar{i}$ in terms of the usual Drinfeld--Cartan series $\boldsymbol\phi^{\pm}(z)$, modulo the Drinfeld positive half of a quantum loop algebra. We refer to such results as `factorization formulae'. In the case of the $q$-Onsager algebra, a  factorization formula for the series $\Thgsr(z)$ was established in \cite{Przez-23}. We recall it below. Since $\left| \indx_0 \right| = 1$, we can omit the subscripts on all the generators. 

\Thm[{\cite[Corollary 4.11]{Przez-23}}] \label{thm: rank 1 factorization}
The series $\Thgsr(z)$ admits the following factorization 
\begin{align*} 
\Thgsr(z) \equiv \boldsymbol\phi^-(z\mi)\boldsymbol{\phi}^+(\CCC z)  \qquad \mod \fext{\Uu_+}{z}. 
\end{align*} 
\enthm

In this paper, we additionally require a factorization formula for the (real rank one) quantum symmetric pair coideal subalgebra of type $\mathsf{AIII}_2^{(\tau)}$. The remainder of the section is devoted to its proof. The argument is inspired by \cite[\S 4]{Przez-23}, but differs from it in some important aspects.

\subsection{Type $\mathsf{AIII}_2^{(\tau)}$: loop presentation}

If $N=2$, Definition \ref{defi: LW pres gen} specializes to the following. 

\Defi 
Let ${}^{\mathrm{Dr}}\Ui(\widehat{\mathfrak{sl}}_3, \tau)$ be the algebra generated by $\KK_i^{\pm1}$, $\Hg_{i,m}$ and $\Ag_{i,r}$, where $m\geq1$, $r\in\Z$ and $i \in \{1,2\}$, as well as a central element $\CCC$, subject to the following relations:
\begin{align} 
\KK_i A_{j,r} &= q^{a_{\tau(i),j} - a_{ij}} A_{j,r} \KK_i, \\
[\Hg_{i,m},\Hg_{j,l}] = [\Hg_{i,m},\KK_{j}] &= [\KK_{i}, \KK_{j}] = 0,  \\
[\Hg_{i,m}, \Ag_{j,r}] &= \textstyle \frac{[m \cdot a_{ij}]}{m} \Ag_{j,r+m} - \frac{[m \cdot a_{\tau(i),j}]}{m} \Ag_{j,r-m}\CCC^m, \label{rank1 HA rel} \\ 
[\Ag_{i,r}, \Ag_{\tau(i),s}] &= \KK_{\tau(i)} \CCC^s\Theta_{i,r-s} - \KK_{i} \CCC^r\Theta_{\tau(i),s-r} \qquad \mbox{if} \quad a_{i,\tau(i)} = 0,  \\ 
[\Ag_{i,r}, \Ag_{i,s+1}]_{q^{-2}}   &= q^{-2} [\Ag_{i,r+1}, \Ag_{j,s}]_{q^{2}},  \\ 
[\Ag_{i,r}, \Ag_{\tau(i),s+1}]_{q}  -q [\Ag_{i,r+1}, \Ag_{\tau(i),s}]_{q\mi}
&= -\KK_i\CCC^r\Thg_{\tau(i),s-r+1} + q\KK_i\CCC^{r+1}\Thg_{\tau(i),s-r-1}  \\ \notag
&\quad  - \KK_{\tau(i)}\CCC^s\Thg_{i,r-s+1}  +q\KK_{\tau(i)} \CCC^{s+1}\Thg_{i,r-s-1}, \\ 
\mathbb{S}(r_1, r_2 \mid s; i, \tau(i)) &= \mathbb{R}(r_1, r_2 \mid s; i, \tau(i)) \qquad \mbox{if} \quad a_{i,\tau(i)} = -1, 
\end{align}
for $m\geq1$; $r,s, r_1, r_2\in \Z$. 
\enDefi 

By \cite[Theorem 5.5]{LWZ-I}, there is an algebra isomorphism ${}^{\mathrm{Dr}}\Ui(\widehat{\mathfrak{sl}}_3, \tau) \isoto{} \Ui(\widehat{\mathfrak{sl}}_3, \tau)$, which is a special case of \eqref{eq: droq->oq}. Its inverse sends
\begin{align*}
\KK_0 &\mapsto -q\mi \CCC (\KK_1 \KK_2)\mi, \quad \KK_i \mapsto \KK_i, \quad B_i \mapsto B_{i,0} \qquad (i \in \{1,2\}), \\
B_0 &\mapsto o(1) q\mi \big( \Theta_{1,1} - q[B_1, B_{2,-1}]_{q\mi} \CCC \KK_2\mi)\KK_1\mi. 
\end{align*}
Note that $\CCC = q^3 c_\delta$. 
We will also need information about some particular root vectors:
\begin{align*}
A_{i,-1} &= - o(i)q^{-2}[B_i,B_0]_q \KK_i \KK_\delta^{-1}, \\
\Theta_{i,1} &= - o(i) \big( [B_i, [B_{\tau(i)}, B_0]_q]_{q^2} - qB_0 \KK_i  \big). 
\end{align*} 
We view $\Ui(\widehat{\mathfrak{sl}}_3, \tau)$ as a subalgebra of $\Uu'$ via \eqref{eq: Kolb emb} (we will in fact work with its image in the quotients $\Uu$ and $\Uq$). 

We will need more notations for various graded subalgebras of $\Uu$.  
Let
\begin{align*}
\Uu_{>k} &= \langle y \in \Uu \mid \forall j \in \{ 1, 2 \}\colon \degdr_j y \geq k;\ \exists j \in \{ 1, 2 \} \colon \degdr_j y > k \rangle, \\ 
\Uu_{\geq k} &= \langle y \in \Uu \mid \forall j \in \{ 1, 2 \}\colon \degdr_j y \geq k \rangle. 
\end{align*}

\subsection{Type $\mathsf{AIII}_2^{(\tau)}$: estimation of $\avar{i}$}

We will need the following lemma. 

\Lem \label{lem: Hernandez rel}
We have 
\[
\boldsymbol\phi^-_i(z\mi) x_{\tau(i),k}^+ = q x_{\tau(i),k}^+ \boldsymbol\phi^-_i(z\mi) + \mathbf{x}^+_{\leq k-1} (q\mi c_\delta^{-1/2} z\mi) \boldsymbol\phi^-_i(z\mi). 
\]
\enlem 

\Proof
The proof is similar to that of \cite[Lemma 4.2]{Przez-23}. 
\enproof

The following estimates are crucial. 

\Lem \label{lem: rank1 H1}
We have 
\begin{align*}
H_{1,i} = \Theta_{i,1} &\equiv - h_{i,-1} + \CCC h_{\tau(i),1} + \mathcal{Y}_i && \mod \Uu_{\geq 2} , \\
A_{i,-1} &\equiv q^{-4} c_\delta^{-1/2} x_{\tau(i),-1}^+ \Ka_i +  x^-_{i,1}  && \mod \Uu_{>1}, 
\end{align*}
where $\degdr \mathcal{Y}_i = \alpha_1 + \alpha_2$ and, explicitly, 
\[
\mathcal{Y}_i = c_\delta^{1/2} \left( (1-q^{-2}) x^+_{\tau(i), -1} x_{i,0}^+ + (q^{-2} - q)  x_{\tau(i),0}^+x^+_{i, -1} +  (1-q) x^+_{i, -1} x_{\tau(i),0}^+ \right). 
\]
\enlem 

\Proof
A direct computation shows that 
\begin{align*}
\Theta_{i,1} &= -o(i)
\Big(  [F_i, [F_{\tau(i)}, F_0]_q]_{q^2} + c_\delta^{1/2} [E_{\tau(i)},[E_i,E_0]_q]_{q^2} + \frac{1-q^2}{q^4} F_0 [E_{\tau(i)},E_i]_{q^3} \widetilde{K}'_{\tau(i)} \widetilde{K}'_i \\
& \quad + 
(q^{-1} - q) [F_i, F_0]_q E_i \widetilde{K}'_{\tau(i)} + (q^{-2}-q) F_0 \widetilde{K}'_{\tau(i)} \widetilde{K}'_i +(q^{-2}-q^2) [F_{\tau(i)},F_0]_q E_{\tau(i)} \widetilde{K}'_i \Big). 
\end{align*}
Let us denote the six summands (including the factor $-o(i)$) in the formula above on the RHS by $S_1, \cdots, S_6$, in the natural order of appearance from top left to bottom right. The summand $S_3$ clearly lies in $\Uu_{\geq 2}$ and can therefore be ignored. 
The other summands need to be converted into Drinfeld generators. One has 
\begin{align*}
x^+_{i,-1} &= o(i) T_{\omega_i} (E_i) = - o(i) c_\delta^{-1/2} [F_{\tau(i)}, F_0]_q \Kb_i , \\ 
\phi_{i,-1}^- &= - \rho c_\delta^{1/2} c_i^{-1/2} [x^+_{i,-1}, x^-_{i,0}], \\ 
h_{i,-1} &= - \rho\mi c_i^{1/2} (\Kb_i)\mi  \phi_{i,-1}^- = c_\delta^{1/2} [x^+_{i,-1}, x^-_{i,0}](\Kb_i)\mi, \\
[x^+_{i,-1}, x^-_{i,0}] &= o(i) c_\delta^{-1/2} [F_i, [F_{\tau(i)}, F_0]_q \Kb_i] = o(i) c_\delta^{-1/2} [F_i, [F_{\tau(i)}, F_0]_q]_{q^2} \Kb_i. 
\end{align*}
Hence 
\[
h_{i,-1} = o(i) [F_i, [F_{\tau(i)}, F_0]_q]_{q^2}. 
\]
Similarly, 
\begin{align*}
x^-_{\tau(i),1} &= o(\tau(i)) T_{\omega_{\tau(i)}} (F_{\tau(i)}) =  - o(i)q\mi c_\delta^{-1/2} \Ka_{\tau(i)} [E_i, E_0]_q  , \\ 
\phi_{\tau(i),1}^+ &= \rho c_{\tau(i)}^{-1/2} [x^+_{\tau(i),0}, x^-_{\tau(i),1}], \\ 
h_{\tau(i),1} &= \rho\mi c_{\tau(i)}^{1/2} \Ka_{\tau(i)}\mi  \phi_{\tau(i),1}^+ = \Ka_{\tau(i)}\mi [x^+_{\tau(i),0}, x^-_{\tau(i),1}], \\
[x^+_{\tau(i),0}, x^-_{\tau(i),1}] &= - o(i)q\mi c_\delta^{-1/2} [E_{\tau(i)}, \Ka_{\tau(i)} [E_i, E_0]_q] \\
 &=- o(i)q^{-3} c_\delta^{-1/2} \Ka_{\tau(i)} [E_{\tau(i)},  [E_i, E_0]_q]_{q^2}. 
\end{align*}
Hence 
\[
h_{\tau(i),1} = -o(i) q^{-3} c_\delta^{-1/2}  [E_{\tau(i)},  [E_i, E_0]_q]_{q^2}. 
\]
These calculations allow us to deduce that 
\[
S_1  = - h_{i,-1}, \quad S_2 = \CCC h_{\tau(i),1}, \quad S_4 = c_\delta^{1/2} (1-q^{-2}) x^+_{\tau(i), -1} x_{i,0}^+. 
\]
Next, we check that 
\begin{align*}
x^+_{\tau(i), 0} x_{i,-1}^+ &= -o(i) c_\delta^{-1/2} E_{\tau(i)} [F_{\tau(i)}, F_0]_q \Kb_i  \\
&= -o(i) c_\delta^{-1/2}   [E_{\tau(i)} F_{\tau(i)}, F_0]_q \Kb_i \\ 
&= -o(i) c_\delta^{-1/2}   [F_{\tau(i)}E_{\tau(i)} + \rho\mi (\Ka_{\tau(i)} - \Kb_{\tau(i)}), F_0]_q \Kb_i \\ 
&= -o(i) c_\delta^{-1/2} \left(  [F_{\tau(i)}, F_0]_q E_{\tau(i)} \Kb_i + F_0 \Kb_{\tau(i)} \Kb_i  \right), \\ 
x_{i,-1}^+ x^+_{\tau(i), 0}  &= -o(i) q c_\delta^{-1/2} [F_{\tau(i)}, F_0]_q E_{\tau(i)} \Kb_i, 
\end{align*}
which implies that 
\[
S_5 + S_6 =  (q^{-2} - q)  x_{\tau(i),0}^+x^+_{i, -1} +  (1-q) x^+_{i, -1} x_{\tau(i),0}^+. 
\]

A direct calculation also shows that 
\begin{align}
A_{i,-1} &= 
 o(i) q^{-4} c_\delta\mi ( [F_i, F_0]_q + q [E_{\tau(i)}, E_0]_q \widetilde{K}_i'  \widetilde{K}_0'  + ( q^{-1} - q ) F_0 E_{\tau(i)}  \widetilde{K}_i' ) \Ka_i \Kb_{\tau(i)}.  
\end{align}
We can ignore the third summand, while the others are easily seen to convert into Drinfeld generators as claimed. 
\enproof

We will abbreviate $\bar{h}_{i,-1} = [2]\mi h_{i,-1}$, \ $\bar{h}_{\tau(i),1} = [2]\mi h_{\tau(i),1}$ and 
\[
\Omega_i = [2] + (\Ad_{h_{i,-1}} - \CCC \Ad_{h_{\tau(i),1}})  z + \CCC z^2 = [2] (1 + \Ad_{\bar{h}_{i,-1}} z) (1 - \Ad_{\bar{h}_{\tau(i),1}} \CCC z).
\]

Lemma \ref{lem: rank1 H1} and relation \eqref{rank1 HA rel} imply that, modulo $\fext{\Uu_{>1}}{z}$, the series $\avar{i}$ consists of terms of degree $- \alpha_i$ and $\alpha_{\tau(i)}$. Let us denote them as $\avarm{i}$ and $\avarp{i}$, respectively. We will use similar notation for modes of these series, e.g., $A_{i,0}^{(-)}$ and $A_{i,0}^{(+)}$. 

\Pro \label{N=2 A est}
We have 
\[
\avar{i} \equiv \avarm{i} + \avarp{i} \quad \mod \fext{\Uu_{>1}}{z},
\]
with
\begin{align*}
\avarm{i} &= \mathbf{x}^-_{i,-}(z\mi), \\
\avarp{i} &=  c_i^{1/2} \Omega_i\mi \Big( -q\mi c_\delta^{1/2} x^+_{\tau(i), -1}z + [2] x^+_{\tau(i), 0} \\
&\quad + (q-1) \mathbf{x}^+_{\tau(i), \leq -1}(q\mi c_\delta^{-1/2} z\mi) \Big) \boldsymbol\phi^-_i(z\mi). 
\end{align*} 
\enpro 

\Proof
Let us first calculate $\avarm{i}$. Lemmas \ref{lem: A Th as series} and \ref{lem: rank1 H1} imply that 
\begin{align*}
\avarm{i}  &= \Omega_i\mi ([2]x_{i,0}^- - \CCC x_{i,1}^- z) = [2]\Omega_i\mi  (1 - [2]\mi \Ad_{h_{\tau(i),1}} \CCC z)(x_{i,0}^-)  \\
&= (1 + [2]\mi \Ad_{h_{i,-1}} z)\mi (x_{i,0}^-) = \mathbf{x}^-_{i,-}(z\mi). 
\end{align*}

Let $Z$ be the term of degree $\alpha_{\tau(i)}$ in $([2] - \Ad_{H_{i,1}} z + \CCC z^2) (\avar{i})$. It follows from Lemmas \ref{lem: A Th as series} and \ref{lem: rank1 H1} that 
\begin{align} \label{eq: Zomega form}
[2] x^+_{\tau(i),0} \Kb_i - c_\delta^{1/2} \Ka_i x_{\tau(i),-1}^+    z = Z &= \Omega_i(\avarp{i}) - \Ad_{\mathcal{Y}_i} (\avarm{i}) z. 
\end{align} 
We calculate
\begin{align*}
[x^+_{\tau(i), -1} x_{i,0}^+, \mathbf{x}^-_{i,-}(z\mi)] &= x^+_{\tau(i), -1} [ x_{i,0}^+, \mathbf{x}^-_{i,-}(z\mi)] 
+ [x^+_{\tau(i), -1}, \mathbf{x}^-_{i,-}(z\mi)] x_{i,0}^+ \\
&= -\rho\mi c_i^{1/2} x^+_{\tau(i), -1} \boldsymbol\phi^-_i(z\mi) + \rho\mi x^+_{\tau(i), -1} \Ka_i, \\
[x_{\tau(i),0}^+x^+_{i, -1}, \mathbf{x}^-_{i,-}(z\mi)]  &= x^+_{\tau(i), 0} [ x_{i,-1}^+, \mathbf{x}^-_{i,-}(z\mi)] 
+ [x^+_{\tau(i), 0}, \mathbf{x}^-_{i,-}(z\mi)] x_{i,-1}^+ \\ 
&= -\rho\mi c_i^{1/2} c_\delta^{-1/2 }x^+_{\tau(i), 0} \boldsymbol\phi^-_i(z\mi) z\mi + \rho\mi c_\delta^{-1/2}x_{\tau(i),0}^+ \Kb_i z\mi \\ 
[x^+_{i, -1}x_{\tau(i),0}^+, \mathbf{x}^-_{i,-}(z\mi)]  &= x_{i,-1}^+ [ x_{i,-1}^+, \mathbf{x}^-_{i,-}(z\mi)] 
+ [x^+_{\tau(i), 0}, \mathbf{x}^-_{i,-}(z\mi)] x^+_{\tau(i), 0} \\ 
&= -\rho\mi c_i^{1/2} c_\delta^{-1/2 } \boldsymbol\phi^-_i(z\mi)x^+_{\tau(i), 0} z\mi + \rho\mi c_\delta^{-1/2}\Kb_i x_{\tau(i),0}^+ z\mi. 
\end{align*} 
Therefore, using Lemma \ref{lem: Hernandez rel}, we get 
\begin{align*}
\Ad_{\mathcal{Y}_i} (\avarm{i}) &= 
 -[2] x^+_{\tau(i),0} \Kb_i z\mi + c_\delta^{1/2} \Ka_i   x^+_{\tau(i),-1} 
+ c_i^{1/2}  \Big( -q\mi c_\delta^{1/2} x^+_{\tau(i), -1} + [2] x^+_{\tau(i), 0}z\mi  \\
&\quad + (q-1) \mathbf{x}_{\tau(i), \leq -1}(q\mi c_\delta^{-1/2} z\mi) z\mi  \Big) \boldsymbol\phi^-_i(z\mi). 
\end{align*} 
Substituting the above result into \eqref{eq: Zomega form} and rearranging yields the proposition. 
\enproof

\subsection{Type $\mathsf{AIII}_2^{(\tau)}$: estimation of $\thvar{i}$} 

Let us abbreviate
\[
\mathbb{B}_i = -q\mi c_\delta^{1/2} x^+_{\tau(i), -1}z + [2] x^+_{\tau(i), 0} + (q-1) \mathbf{x}^+_{\tau(i), \leq -1}(q\mi c_\delta^{-1/2} z\mi), \quad \mathbb{A}_i = c_i^{1/2} \Omega_i\mi (\mathbb{B}_i),
\]
so that $\avarp{i} = \mathbb{A}_i \boldsymbol\phi^-_i(z\mi)$. 

\Lem \label{lem: AA into x}
We have 
\[
(1 - q \Ad_{\bar{h}_{\tau(i),-1}} z) \mathbb{A}_i = c_i^{1/2} (1 - \Ad_{\bar{h}_{\tau(i),-1}} z) \mathbf{x}_{\tau(i),+}^+(c_\delta^{-1/2}\CCC z). 
\]
\enlem

\Proof
We calculate
\begin{align*}
(q-1) (1 - q \Ad_{\bar{h}_{\tau(i),-1}} z) \mathbf{x}^+_{\tau(i), \leq -1}(q\mi c_\delta^{-1/2} z\mi) &= (q^2-q) c_\delta^{1/2} x_{\tau(i),-1}^+ z, \\ 
[2] (1 - q \Ad_{\bar{h}_{\tau(i),-1}} z) x_{\tau(i),0}^+ &= [2] x_{\tau(i),0}^+ - (q^2+1)c_\delta^{1/2} x_{\tau(i),-1}^+ z, \\ 
-q\mi c_\delta^{1/2} z (1 - q \Ad_{\bar{h}_{\tau(i),-1}} z) x_{\tau(i),-1}^+ &= -q\mi c_\delta^{1/2} x_{\tau(i),-1}^+ z + c_\delta x_{\tau(i),-2}^+ z^2. 
\end{align*}
Therefore, 
\begin{align*}
(1 - q \Ad_{\bar{h}_{\tau(i),-1}} z) \mathbb{B}_i &= [2] x_{\tau(i),0}^+ - (q\mi + 1 + q) c_\delta^{1/2} x_{\tau(i),-1}^+ z + c_\delta x_{\tau(i),-2}^+ z^2 \\ 
&= [2](1 + \Ad_{\bar{h}_{i,-1}} z) (1 - \Ad_{\bar{h}_{\tau(i),-1}} z) x_{\tau(i),0}^+. 
\end{align*}
Applying $c_i^{1/2} \Omega_i\mi$ yields 
\begin{align*}
(1 - q \Ad_{\bar{h}_{\tau(i),-1}} z) \mathbb{A}_i &= c_i^{1/2} \frac{(1 - \Ad_{\bar{h}_{\tau(i),-1}} z)}{(1 - \Ad_{\bar{h}_{\tau(i),1}} \CCC z)} x_{\tau(i),0}^+ \\
&= c_i^{1/2} (1 - \Ad_{\bar{h}_{\tau(i),-1}} z) \mathbf{x}_{\tau(i),+}^+(c_\delta^{-1/2}\CCC z), 
\end{align*}
completing the proof. 
\enproof

\Thm \label{thm: rank 1 factorization aij-1} 
The `Drinfeld--Cartan' series $\thvar{i}$ admits, for each $i \in \{1,2\}$, the following factorization: 
\begin{align*}
\thvar{i} \equiv K_i K_{\tau(i)}\mi   \boldsymbol\phi^-_i(z\mi)\boldsymbol{\phi}^+_{\tau(i)}(\CCC z)  \mod \fext{\Uq_+}{z}. 
\end{align*}
\enthm

\Proof
Consider the expression $z\mi[A_{\tau(i),-1}, \avar{i}]_{q} - q[A_{\tau(i),0}, \avar{i}]_{q\mi}$. 
It follows from the first relation in \eqref{eq: rel sln xx} that 
\begin{align*}
z\mi[A_{\tau(i),-1}^{(-)}, \avarm{i}]_{q} &- q[A_{\tau(i),0}^{(-)}, \avarm{i}]_{q\mi} = \\ 
&= z\mi[x_{\tau(i),1}^-, \mathbf{x}^-_{i,-}(z\mi)]_{q} - q[x_{\tau(i),0}^-, \mathbf{x}^-_{i,-}(z\mi)]_{q\mi} \\ 
&= z\mi [x_{\tau(i),1}^-, {x}^-_{i,0}]_{q} = o(i) c_0^{-1/2} q^{-2} E_0 K_{\tau(i)} K_i z\mi \\
&= o(i) c_0^{-1/2} q^{-2} E_0 K_0\mi z\mi. 
\end{align*}

Next, one has 
\begin{align*}
z\mi[A_{\tau(i),-1}^{(+)},& \ \avarm{i}]_{q} - q[A_{\tau(i),0}^{(+)}, \avarm{i}]_{q\mi} = \\
&= z\mi q^{-4}c_\delta^{-1/2}[ x^+_{i,-1} \Ka_{\tau(i)}, \mathbf{x}^-_{i,-}(z\mi)]_{q} - q[x^+_{i,0} \Kb_{\tau(i)}, \mathbf{x}^-_{i,-}(z\mi)]_{q\mi} \\
&= z\mi q^{-3}c_\delta^{-1/2}[ x^+_{i,-1}, \mathbf{x}^-_{i,-}(z\mi)] \Ka_{\tau(i)} - [x^+_{i,0}, \mathbf{x}^-_{i,-}(z\mi)] \Kb_{\tau(i)} \\
&=  -(\rho z^2 \CCC)\mi c_i^{1/2} \left( c_i^{1/2} \boldsymbol\phi^-_{i}(z\mi) - \Kb_i  \right) \Ka_{\tau(i)} +  \rho\mi \left( c_i^{1/2} \boldsymbol\phi^-_{i}(z\mi) - \Ka_i  \right) \Kb_{\tau(i)} \\
&= c_i \rho\mi \left(K_{\tau(i)}\mi - \CCC\mi K_{\tau(i)} z^{-2} \right) \boldsymbol\phi^-_{i}(z\mi) - \rho\mi \left( \KK_i - \CCC\mi z^{-2} \KK_{\tau(i)} \right). 
\end{align*} 

Moreover, 
\begin{align*}
z\mi[A_{\tau(i),-1}^{(-)}, \avarp{i}]_{q} &- q[A_{\tau(i),0}^{(-)}, \avarp{i}]_{q\mi} = \\
&= z\mi [ x_{\tau(i),1}^-, \mathbb{A}_i \boldsymbol\phi_i^-(z\mi)]_q -q[x_{\tau(i),0}^-, \mathbb{A}_i \boldsymbol\phi_i^-(z\mi)]_{q\mi} \\
&= \left( z\mi [ x_{\tau(i),1}^-, \mathbb{A}_i] -q[x_{\tau(i),0}^-, \mathbb{A}_i] \right) \boldsymbol\phi_i^-(z\mi) \\
&= z\mi \left[ (1 + q \Ad_{\bar{h}_{\tau(i),-1}} z)x^-_{\tau(i),1}, \ \mathbb{A}_i \right] \boldsymbol\phi_i^-(z\mi), 
\end{align*} 
where, in the second equality, we used the fact that 
\[
\boldsymbol\phi_i^-(z\mi) \left( x_{\tau(i),0}^- -q z\mi x_{\tau(i),1}^- \right) =\left( qx_{\tau(i),0}^- -z\mi x_{\tau(i),1}^- \right) \boldsymbol\phi_i^-(z\mi). 
\]
Next, we compute 
\begin{align*}
\left[ (1 + q \Ad_{\bar{h}_{\tau(i),-1}} z)x^-_{\tau(i),1}, \ \mathbb{A}_i \right] &= 
\left[ x^-_{\tau(i),1}, \ (1 - q \Ad_{\bar{h}_{\tau(i),-1}} z) \mathbb{A}_i \right] \\ 
&=  c_i^{1/2} \left[ x^-_{\tau(i),1}, \  (1 - \Ad_{\bar{h}_{\tau(i),-1}} z) \mathbf{x}_{\tau(i),+}^+(c_\delta^{-1/2}\CCC z) \right] \\ 
&=  c_i^{1/2} \left[ (1 + \Ad_{\bar{h}_{\tau(i),-1}} z) x^-_{\tau(i),1}, \   \mathbf{x}_{\tau(i),+}^+(c_\delta^{-1/2}\CCC z) \right] \\ 
&=  c_i^{1/2} \left[  x^-_{\tau(i),1} - x^-_{\tau(i),0}z , \   \mathbf{x}_{\tau(i),+}^+(c_\delta^{-1/2}\CCC z) \right] \\ 
&= - c_i \rho\mi \CCC\mi z\mi \left( (1 - \CCC z^2) \boldsymbol\phi^+_{\tau(i)}(\CCC z) - K_{\tau(i)} + \CCC z^2 K_{\tau(i)}\mi \right). 
\end{align*} 
In the first and third equalities, we used the fact that $ \Ad_{\bar{h}_{\tau(i),-1}}$ is a derivation acting trivially on elements of degree zero. For the second equality, we used Lemma \ref{lem: AA into x}. 
As a consequence, 
\begin{align*}
z\mi[A_{\tau(i),-1}^{(+)},& \ \avarm{i}]_{q} - q[A_{\tau(i),0}^{(+)}, \avarm{i}]_{q\mi} = \\
&= - c_i \rho\mi \CCC\mi z^{-2} \left( (1 - \CCC z^2) \boldsymbol\phi^+_{\tau(i)}(\CCC z) - K_{\tau(i)} + \CCC z^2 K_{\tau(i)}\mi \right) \boldsymbol\phi_i^-(z\mi). 
\end{align*} 

Moreover, $\boldsymbol\theta_n$ is trivial, and so 
\begin{align*}
o(i) q \CCC\mi \KK_{\tau(i)} \Tbr_{\boldsymbol\theta_n}\mi(B_0) \KK_{i} z\mi &\equiv 
o(i) q \CCC\mi \KK_{\tau(i)} E_0 \Kb_0 \KK_{i} z\mi \\
&= -o(i) E_0 \Kb_0 \KK_0\mi z\mi =  o(i) q^{-2} c_0^{-1/2} E_0 K_0\mi z\mi. 
\end{align*}
Using \eqref{eq: case3 theta formula}, it follows that 
\begin{align*} 
\KK_{\tau(i)} (q - \CCC\mi z^{-2}) \boldsymbol\Theta_{i}(z) &= 
z\mi[A_{\tau(i),-1}, \avar{i}]_{q} - q[A_{\tau(i),0}, \avar{i}]_{q\mi} \\
&\quad - o(i) q \CCC\mi \KK_{\tau(i)} \Tbr_{\boldsymbol\theta_n}\mi(B_0) \KK_{i} z\mi 
- \rho\mi (\CCC\mi \KK_{\tau(i)} z^{-2} - \KK_{i}) \\
&\equiv - c_i \rho\mi \CCC\mi z^{-2} (1 - \CCC z^2) \boldsymbol\phi^+_{\tau(i)}(\CCC z) \boldsymbol\phi_i^-(z\mi). 
\end{align*} 
Rearranging, we get
\[
\thvar{i} = \frac{\rho(1- q \CCC z^2)}{1-\CCC z^2} \boldsymbol\Theta_i(z) \equiv K_i K_{\tau(i)}\mi   \boldsymbol\phi^-_i(z\mi)\boldsymbol{\phi}^+_{\tau(i)}(\CCC z), 
\]
completing the proof. 
\enproof

\section{Root and weight combinatorics} 
\label{sec: roots n weight}

The relative root systems associated to affine symmetric pairs of types $\mathsf{AIII}_{N}^{(\tau)}$ are the affine root systems of type $\mathsf{C}_n^{(1)}$ and $\mathsf{A}_{2n}^{(2)}$. 
In this section, we first recall some combinatorial results concerning the root system of type $\mathsf{C}_n^{(1)}$ from \cite{LP25}, and then state corresponding results for type $\mathsf{A}_{2n}^{(2)}$.

\subsection{Reduced expressions for fundamental weights}
 
We denote $[a,b] = r_a r_{a+1} \cdots r_b$ for $a \le b$ and $[a,b] = r_{a} r_{a-1} \cdots r_b$ for $a>b$. 

\subsubsection{Type $\mathsf{AIII}_{2n-1}^{(\tau)}$} 

The relative root system is of type $\mathsf{C}_n^{(1)}$ (for $n \geq 2$). 
The fundamental group is $\Lambda^\tau = \{ 1, \pi_n\}$, where $\pi_n$ interchanges the roots $\balpha_k$ and $\balpha_{n-k}$. If $\Lambda = \{1, \pi, \cdots, \pi^{2n-1} \}$ is the fundamental group for type $\mathsf{A}_{2n-1}$, with $\pi$ shifting the affine Dynkin diagram by one, then $\pi_n = \pi^{n}$ under the canonical inclusion $\widetilde{W}_{{\mathsf{AIII}}_{2n-1}^{(\tau)}} \hookrightarrow \widetilde{W}_{\mathsf{A}_{2n-1}^{(1)}}$. The fundamental weights can be expressed using the following formulae.

\Pro \label{prop:reduced expression of omegai in extended affine Weyl groups type AIIItau}
Fundamental weights in $\widetilde{W}_{{\mathsf{AIII}}_{2n-1}^{(\tau)}} \cong \widetilde{W}_{\mathsf{C}_n^{(1)}}$ are given by the following.  If $\tau(i) \ne i$, then
\begin{align*}
\varpi_i = (r_0 [1,n])^i [n-i, n-1] \cdots [2,i+1][1,i]. 
\end{align*} 
If $i=n$, then 
\begin{align*}
\varpi_n = \pi_n r_n [n-1,n] \cdots [1,n] = \pi_n [n,1][n,2] \cdots [n,n-1]r_n. 
\end{align*} 
Moreover, the formulae above yield reduced expressions. 
\enpro 

\Proof
This follows from \cite[Proposition 3.6]{LP25}. 
\enproof

\subsubsection{Type $\mathsf{AIII}_{2n}^{(\tau)}$}

The relative root system is of type $\mathsf{A}_n^{(2)}$. The fundamental weights can be expressed using the following formulae. 

\Pro \label{prop:reduced expression of omegai in extended affine Weyl groups type AIII even}
Fundamental weights in $\widetilde{W}_{{\mathsf{AIII}}_{2n}^{(\tau)}} \cong \widetilde{W}_{\mathsf{A}_{2n}^{(2)}}$  are given by 
\begin{align*}
\varpi_{i} = (r_0 [1,n])^i [n-i, n-1] [n-i-1,n-2] \cdots [1,i], 
\end{align*}  
for $1 \le i \le n$. 
Moreover, the formulae above yield reduced expressions. 
\enpro 

\Proof
These formulas can be obtained from those in \cite[Corollary 2.10]{LPWZ} by easy algebraic manipulations. 
\enproof

\subsection{Root combinatorics} 

We will also need information about the action of some special elements of the relative braid group. 

\subsubsection{Type ${\mathsf{AIII}}_{2n-1}^{(\tau)}$} 

Let $\zeta_i = (r_0r_1\cdots r_n)^i \in \widetilde{W}_{{\mathsf{AIII}}_{2n-1}^{\tau}} \subset \widetilde{W}_{\mathsf{A}_{2n-1}^{(1)}}$, for $1 \leq i \leq n{-}1$. Set 
\begin{align*}
\tilde{\alpha}_k = \begin{cases}
\alpha_k & 1 \leq k \leq n-1, \\
\alpha_0 + \sum_{j=0}^{n-1} \alpha_{n+j} & k=0, 
\end{cases}
\qquad 
\tilde{\alpha}'_k = \begin{cases}
\alpha_{2n-k} & 1 \leq k \leq n-1, \\
\sum_{j=0}^{n} \alpha_j & k=0. 
\end{cases}
\end{align*} 
Note that we are dealing with roots in the \emph{non-relative} root system. 

\Lem \label{lem: action on sroots 1}
Let $1 \leq k \leq n$ and $1 \leq i \leq n-1$. Then 
\[
\zeta_i \cdot \tilde{\alpha}_k = \tilde{\alpha}_{k+i}, \qquad \zeta_i \cdot \tilde{\alpha}'_k = \tilde{\alpha}'_{k+i}, 
\]
with indices taken modulo $n$. 
\enlem 

\Proof
A straightforward calculation shows that $(r_0r_1\cdots r_n) \cdot \tilde{\alpha}_k = \tilde{\alpha}_{k+1}$ and \linebreak $(r_0r_1\cdots r_n) \cdot \tilde{\alpha}'_k = \tilde{\alpha}'_{k+1}$. Repeating this $i$-times yields the lemma. 
\enproof

\subsubsection{Type $\mathsf{AIII}_{2n}^{(\tau)}$}  

Let $\zeta_i = (r_0r_1\cdots r_n)^i \in \widetilde{W}_{{\mathsf{AIII}}_{2n}^{\tau}} \subset \widetilde{W}_{\mathsf{A}_{2n}^{(1)}}$, for $1 \leq i \leq n$. Set 
\begin{align*}
\tilde{\alpha}_k = \begin{cases}
\alpha_k & 1 \leq k \leq n-1, \\
\alpha_0 + \sum_{j=0}^{n} \alpha_{n+j} & k=0, 
\end{cases}
\qquad 
\tilde{\alpha}'_k = \begin{cases}
\alpha_{2n-k+1} & 1 \leq k \leq n-1, \\
\sum_{j=0}^{n+1} \alpha_j & k=0. 
\end{cases}
\end{align*} 

\Lem \label{lem: action on sroots AIII even}
Let $1 \leq k \leq n$ and $1 \le i \le n$. Then 
\[
\zeta_i \cdot \tilde{\alpha}_k = \tilde{\alpha}_{k+i}, \qquad \zeta_i \cdot \tilde{\alpha}'_k = \tilde{\alpha}'_{k+i}, 
\]
with indices taken modulo $n$.
\enlem 

\Proof
This is a straightforward calculation, as above. 
\enproof

\section{Compatibility of braid group actions} 
\label{sec: compat of br}

The goal of this section is to express $\Tbr_{\varpi'_i}(B_i)$ in terms of elements defined using the usual Lusztig braid group action, modulo an appropriate subalgebra of the Drinfeld positive half of $\Uu$. More precisely, we will prove that 
\[
\Tbr_{\varpi'_i}(B_i) \equiv T_{\omega_i'}(F_i) + \gamma  T_{\omega'_{\tau(i)}}(E_{\tau(i)}) T_{\omega_i'}(\Kb_i) \quad \mod \Uu_{d_{i} \geq 1,+}, 
\]
if $a_{i,\tau(i)} \in \{0,2\}$, and 
\[
\Tbr_{\varpi'_i}(B_i) \equiv T_{\omega'_{\tau(i)}}(F_{\tau(i)}) + \gamma T_{\omega'_{i}}(E_{i})T_{\omega'_{\tau(i)}}(\widetilde{K}'_{\tau(i)}) 
\quad \mod \Uu_{d_{\tau(i)} \geq 1,+}, 
\]
if $a_{i,\tau(i)} = -1$, for some scalar $\gamma$. 
We do this by expressing $\Tbr_{\varpi'_i}(B_i)$ as an explicit polynomial in the generators $B_j$ and showing that appropriate terms of negative degree vanish.

\nc{\tE}{\widetilde{E}}

\subsection{Iterated $q$-brackets} 

Define non-commutative polynomials $P(y_1, \cdots, y_k)$ and $P'(y_1, \cdots, y_k)$ over $\C$ by induction in the following way: 
\begin{alignat}{3}
P(y_1) =& \ y_1, \qquad& P(y_1, \cdots, y_{k+1}) =& \ P(y_1, \cdots, y_{k-1}, [y_k,y_{k+1}]_q), \\
P'(y_1) =& \ y_1, \qquad& P'(y_1, \cdots, y_{k+1}) =& \ [P'(y_1, \cdots, y_k), y_{k+1}]_q.
\end{alignat}
Clearly, we have 
\begin{align}
P(y_1, \cdots, y_k) =& \ P(y_1, \cdots, y_l, P(y_{l+1}, \cdots, y_k)), \\
P'(y_1, \cdots, y_k) =& \ P'(P'_l(y_1, \cdots, y_l), y_{l+1}, \cdots, y_k). 
\end{align}
for any $1 \leq l \leq k-1$. 

We say that a tuple $(y_1, \cdots, y_k)$ is \emph{almost commuting} if $y_my_n = y_ny_m$, for $|m-n| > 1$. Recall the following two lemmas from \cite[\S5]{LP25}. 

\Lem \label{lem: PnP'}
The polynomials $P$ and $P'$ are equal if $(y_1, \cdots, y_k)$ is almost commuting. In particular, in that case, 
\[
P(y_1, P(y_2, y_3)) = P'(P'(y_1, y_2), y_3).
\]
\enlem

\Lem \label{lem: exchange}
If $y_my_{m+1} = y_{m+1}y_m$ then
\[
P(\cdots, y_m, y_{m+1}, \cdots) = P(\cdots, y_{m+1}, y_{m}, \cdots). 
\]
\enlem

Let us abbreviate 
\[ \tE_i = E_{\tau(i)} \Kb_i. 
\] 
We will also need the following polynomials: 
\begin{align*}
{\bf P}(X_{i},Y,Z) = [X_{i}, [Y, Z]_q]_q - q \KK_{i} Z, 
\end{align*}
for $X \in \{B, F, \widetilde{E}, E \}$, and $Y,Z$ arbitrary.

\subsection{Type ${\sf AIII}_{2n-1}^{(\tau)}$} 

The following lemma expresses $\Tbr_{\varpi'_i}(B_i)$ as an explicit polynomial.

\Lem \label{lem: T B pol C} 
If $n \geq 2$, then 
\eq \label{eq: C good pol}
\Tbr_{\varpi'_i}(B_i) = 
\begin{cases}
P(B_{i-1}, \cdots B_1, P(B_{i+1}, \cdots, B_{2n-1},B_0)) & \mbox{if } 1 \leq i < n, \\[3pt]
P(B_{i+1}, \cdots B_{2n-1}, P(B_{i-1}, \cdots, B_{1},B_0)) & \mbox{if } n < i < 2n, \\[3pt]
{\bf P}( B_{n+1}, B_{n-1}, {\bf P}( \cdots {\bf P}( B_{2n-2}, B_{2}, {\bf P}( B_{2n-1}, B_{1}, B_{0} )  ) )) & \mbox{if } i=n. 
\end{cases}
\eneq 
\enlem

\begin{proof}
Let $1 \leq i<n$. By Proposition \ref{prop:reduced expression of omegai in extended affine Weyl groups type AIIItau}, one has $\varpi_i' = \zeta_i \tau_i$, where $\zeta_i=(r_0 [1,n])^i$ and $\tau_i=[n{-}i, n{-}1] \cdots [2,i{+}1][1,i{-}1]$. 
By \cite[Lemma 5.3]{LP25}, 
\begin{align*}
\Tbr_{\tau_i}(B_i) = P(B_{n-1}, \ldots, B_{n-i+1}, P(B_1, \ldots, B_{n-i})),
\end{align*}
and, by Lemma \ref{lem: action on sroots 1}, $\Tbr_{\zeta_i}(B_j) = B_{(j+i) \: {\rm mod} \: n}$ if $j \ne n-i$. Therefore
\begin{align*}
\Tbr_{\varpi'_i}(B_i) = P(B_{i-1}, \ldots, B_{1}, P(B_{i+1}, \ldots, B_{n-1}, \Tbr_{\zeta_i}( B_{n-i} ))).
\end{align*}
Since, by Lemma \ref{lem: TjTiBj formula} (second case), 
\begin{align*}
\Tbr_{\zeta_i}( B_{n-i} ) = \Tbr_{r_0[1,n]}( B_{n-1} ) = P(B_n, \ldots, B_{2n-1},B_0),
\end{align*}
the result follows. The case $ n < i < 2n$ is proven similarly. 

Now let $i=n$. By Proposition \ref{prop:reduced expression of omegai in extended affine Weyl groups type AIIItau}, one has  $\varpi_n' = \pi_n [n,1] \cdots [n, n-1]$. By Lemma~\ref{lem: TjTiBj formula}, we get $\Tbr_{[n, n-1]}(B_n) = {\bf P}(B_{n-1}, B_{n+1}, B_n)$. Therefore
\begin{align*}
\Tbr_{[n,n-2] [n, n-1]}(B_n) = {\bf P}(  B_{n-2}, B_{n+2}, {\bf P}( B_{n-1}, B_{n+1}, B_n) ).
\end{align*} 
Repeating this procedure, we deduce that, for $1 \le k < n$, $\Tbr_{[n,k] \cdots [n, n-1]}(B_n)$ is equal to 
\begin{align*}
{\bf P}( B_{k}, B_{2n-k}, {\bf P}( \cdots {\bf P}( B_{n-2}, B_{n+2}, {\bf P}( B_{n-1}, B_{n+1}, B_{n} )  ) )). 
\end{align*} 
In particular, $\Tbr_{[n,1] \cdots [n, n-1]}(B_n)$ is 
\begin{align*}
{\bf P}( B_{1}, B_{2n-1}, {\bf P}( \cdots {\bf P}( B_{n-2}, B_{n+2}, {\bf P}( B_{n-1}, B_{n+1}, B_{n} )  ) )). 
\end{align*} 
Applying ${\bf T}_{\pi_n}$, we obtain the desired result. 
\end{proof}

\Pro \label{pro: Aodd igood}
For $1 \leq i \leq 2n-1$, one has  
\[
\Tbr_{\varpi'_i}(B_i) \equiv T_{\omega_i'}(F_i) + q^{2n -2\delta_{i,n}}  T_{\omega'_{\tau(i)}}(E_{\tau(i)}) T_{\omega_i'}(\Kb_i) \quad \mod \Uu_{d_{i} \geq 1,+}. 
\]
\enpro 

\Proof
Let $1 \leq i < n$ (the argument for $n < i < 2n$ is analogous, so we will omit it).  
By, e.g., \cite[Lemma 5.3]{LP25}, 
\[
P(F_{i-1}, \cdots F_1, P(F_{i+1}, \cdots, F_{2n-1}, F_0))
= T_{\omega_i'}(F_i). 
\]
One then sees immediately from Lemma \ref{lem: T B pol C} and \eqref{eq: Kolb emb} that 
\[
\Tbr_{\varpi'_i}(B_i) \equiv T_{\omega_i'}(F_i) + P(B_{i-1}, \cdots B_1, P(B_{i+1}, \cdots, B_{2n-1}, \tE_0)) \quad \mod \Uu_{d_{i} \geq 1,+}. 
\]
We will show that 
\eq
P(B_{i-1}, \cdots B_1, P(B_{i+1}, \cdots, B_{2n-1}, \tE_0)) = q^{2n}  T_{\omega'_{\tau(i)}}(E_{\tau(i)}) T_{\omega_i'}(\Kb_i). 
\eneq

Write $\Tbr_{\varpi'_i}(B_i)$ as 
\[
\Tbr_{\varpi'_i}(B_i) = P((B_{i{-}1}, \cdots, B_1), P((B_{i{+}1}, \cdots, B_{n{-}1}), {P}(B_n, \cdots, B_{2n{-}1}, B_0))). 
\]
First consider the innermost nested polynomial ${P}(B_n, \cdots, B_{2n-1}, B_0)$. 
We will show by induction that 
\[
{P}(B_n, \cdots, B_{2n-1}, \tE_0) = {P}(\tE_n, \cdots, \tE_{2n-1}, \tE_0). 
\]
Define 
\[
R_k = \left\{
\begin{array}{r l}
\tE_0 & \quad \mbox{if} \ k=1, \\
P(B_{2n+1-k}, \cdots, B_{2n-2}, B_{2n-1}, \tE_0) & \quad \mbox{if} \ 2 \leq k \leq n+1. 
\end{array}
\right.
\]
Let $S_k$ be defined analogously, with each $B_l$ replaced by $\tE_l = E_{\tau(l)} \Kb_l$. 
Observe that $R_1 = S_1$ trivially, and $R_2 = S_2$ by  Lemma \ref{lem: FE vanishing}. For $k \geq 2$, 
we have 
\begin{align*}
R_{k+1} - S_{k+1} = P(F_{2n-k}, S_k) = P(F_{2n-k}, \tE_k, S_{k-1}) = [[F_{2n-k}, \tE_{2n+1-k}]_q, S_{k-1}]_q = 0, 
\end{align*}
where the first equality follows from induction, the third from Lemma \ref{lem: PnP'} (since $F_{2n-k}$ and $S_{k-1}$ commute), and the fourth from the vanishing of $[F_{2n-k}, \tE_{2n+1-k}]_q = 0$ due to Lemma~\ref{lem: FE vanishing}.  
Continuing the induction, we conclude that $R_{n+1} = S_{n+1}$. 

Next, we consider the middle nested polynomial and show that
\[
P(B_{i+1}, \cdots, B_{n{-}1}, R_{n+1}) = P(\tE_{i+1}, \cdots, \tE_{n{-}1}, S_{n+1}).
\]
Set
\[
R_k^{\mathsf{mid}} = \left\{
\begin{array}{r l}
R_{n+1} & \quad \mbox{if} \ k=0, \\
P(B_{n{-}k}, \cdots, B_{n{-}1}, R_{n+1}) & \quad \mbox{if} \ 1 \leq k \leq n{-}i{-}1,
\end{array}
\right.
\]
and define $S_k^{\mathsf{mid}}$ analogously, with each $B_l$ replaced by $\tE_l$, and $R_{n+1}$ by $S_{n+1}$. 
We already know that $R_0^{\mathsf{mid}} = S_0^{\mathsf{mid}}$. Moreover, $R_1^{\mathsf{mid}} = R_1^{\mathsf{mid}}$ as well since, using Lemma~\ref{lem: PnP'} repeatedly, we get 
\begin{align*}
P(F_{n-1}, \tE_{n}, \tE_{n+1}, \tE_{n+2}, S_{n-3}) &= 
[F_{n-1}, P(\tE_{n}, \tE_{n+1}, \tE_{n+2}, S_{n-3})]_q = \\ 
= [F_{n-1}, P'(\tE_{n}, \tE_{n+1}, \tE_{n+2}, S_{n-3})]_q 
&= P(F_{n-1}, P'(\tE_{n}, \tE_{n+1}, \tE_{n+2}), S_{n-3}) = \\ 
= P'(F_{n-1}, P'(\tE_{n}, \tE_{n+1}, \tE_{n+2}), S_{n-3}) 
&= [[F_{n-1}, P'(\tE_{n}, \tE_{n+1}, \tE_{n+2})]_q, S_{n-3}]_q = \\
= [[F_{n-1}, P(\tE_{n}, \tE_{n+1}, \tE_{n+2})]_q, S_{n-3}]_q 
&= [P(F_{n-1}, \tE_{n}, \tE_{n+1}, \tE_{n+2}), S_{n-3}]_q. 
\end{align*}
The last expression vanishes because, by Lemma \ref{lem: P4 van}, $P(F_{n-1}, \tE_{n}, \tE_{n+1}, \tE_{n+2}) = 0$.

Next, we claim that $[F_{n{-}k}, S_{k{-}2}^{\mathsf{mid}}] = 0$ for each $2 \leq k \leq n{-}i{-}1$. 
Since $S_{k{-}2}^{\mathsf{mid}} = P(\tE_{n{-}k+2}, \cdots, \tE_{n{-}1}, S_{n{+}1})$, and $F_{n-k}$ commutes with $\tE_{n{-}k+2}, \cdots, \tE_{n{-}1}$, 
it suffices to show that $[F_{n{-}k}, S_{n{+}1}] = 0$. The latter is indeed true because, by Lemma \ref{lem: action on sroots 1}, 
\[
[F_{n{-}k}, S_{n{+}1}] = T_{\zeta_i} ( [F_{n{-}k{-}i}, \tE_{n{-}i}] ) = 0. 
\] 
Therefore, 
\[ 
R_k^{\mathsf{mid}} - S_k^{\mathsf{mid}} = 
P(F_{n{-}k}, S_{k{-}1}^{\mathsf{mid}}) = [[F_{n{-}k}, \tE_{n{-}k+1}]_q, S_{k{-}2}^{\mathsf{mid}}]_q = 0. 
\]
By induction, we conclude that $R_{n{-}i{-}1}^{\mathsf{mid}} = S_{n{-}i{-}1}^{\mathsf{mid}}$. 

Finally, consider the outer polynomial $P(B_{i{-}1}, \cdots, B_{1}, R_{n{-}i{-}1}^{\mathsf{mid}})$ and show that 
\[
P(B_{i{-}1}, \cdots, B_{1}, R_{n{-}i{-}1}^{\mathsf{mid}}) = P(\tE_{i{-}1}, \cdots, \tE_{1}, S_{n{-}i{-}1}^{\mathsf{mid}}). 
\]
Set
\[
R_k^{\mathsf{out}} = \left\{
\begin{array}{r l}
R_{n{-}i{-}1}^{\mathsf{mid}} & \quad \mbox{if} \ k=0, \\
P(B_{k}, \cdots, B_{1}, R_{n{-}i{-}1}^{\mathsf{mid}}) & \quad \mbox{if} \ 1 \leq k \leq i{-}1,
\end{array}
\right.
\]
and define $S_k^{\mathsf{out}}$ analogously, with each $B_l$ replaced by $\tE_l$, and $R_{n{-}i{-}1}^{\mathsf{mid}}$ by $S_{n{-}i{-}1}^{\mathsf{mid}}$. 
We already know that $R_0^{\mathsf{out}} = R_0^{\mathsf{out}}$. Moreover, $R_1^{\mathsf{out}} = S_1^{\mathsf{out}}$ since
\begin{align}
P(F_{1}, S_{n{-}i{-}1}^{\mathsf{mid}}) =& \ 
T_{\zeta_i}(P(F_{n{-}i+1}, P(\tE_1, \cdots, \tE_{n{-}i} ) )) \\
=& \ T_{\zeta_i}(P(\tE_1, \cdots, \tE_{n{-}i{-}1}, P(F_{n{-}i+1}, \tE_{n{-}i}))), 
\end{align} 
and $P(F_{n{-}i+1}, \tE_{n{-}i}) = 0$ by Lemma \ref{lem: FE vanishing}. 
We also claim that $[F_{k}, S_{k{-}2}^{\mathsf{out}}] = 0$ for each $2 \leq k \leq i{-}1$. Since $S_{k{-}2}^{\mathsf{out}} = P(\tE_{k{-}2}, \cdots, \tE_{1}, S_{n{-}i{-}1}^{\mathsf{mid}})$, it suffices to show that $[F_{k}, S_{n{-}i{-}1}^{\mathsf{mid}}]~=~0$. This is indeed true because, by Lemma \ref{lem: action on sroots 1}, 
\[
[F_{k}, S_{n{-}i{-}1}^{\mathsf{mid}}] = T_{\zeta_i} ( [F_{n+k{-}i}, P(\tE_1, \cdots, \tE_{n{-}i})] ) = 0. 
\]
It follows that  
\[
R_k^{\mathsf{out}} - S_k^{\mathsf{out}} = P(F_{k}, S_{k{-}1}^{\mathsf{out}}) = [[F_{k}, \tE_{k{-}1}]_q, S_{k{-}2}^{\mathsf{out}}]_q = 0. 
\]
By induction, we conclude that $R_{i{-}1}^{\mathsf{out}} = S_{i{-}1}^{\mathsf{out}}$. 
Altogether, we have proven that 
\[
P(B_{i-1}, \cdots B_1, P(B_{i+1}, \cdots, B_{2n-1},\tE_0)) =
P(\tE_{i-1}, \cdots \tE_1, P(\tE_{i+1}, \cdots, \tE_{2n-1},\tE_0)). 
\]
Finally, one checks that 
\begin{align*}
P(\tE_{i-1}, \cdots \tE_1, P(\tE_{i+1}, \cdots, \tE_{2n-1},\tE_0)) = q^{2n} T_{\omega'_{\tau(i)}}(E_{\tau(i)}) T_{\omega_i'}(\Kb_i). 
\end{align*}

Next, assume that $i=n$. Let 
\begin{align}
P_k &= {\bf P}( B_{2n-k}, B_{k}, {\bf P}( \cdots {\bf P}( B_{2n-2}, B_{2}, {\bf P}( B_{2n-1}, B_{1}, F_{0} )  ) )), \\ 
Q_k &= 
P(F_{k}, \cdots, F_1, F_{2n-k}, \cdots, F_{2n-1}, F_0), \\
R_k &= {\bf P}( B_{2n-k}, B_{k}, {\bf P}( \cdots {\bf P}( B_{2n-2}, B_{2}, {\bf P}( B_{2n-1}, B_{1}, \tE_{0} )  ) )), \\ 
S_k &= P(E_{k}, \cdots, E_1, E_{2n-k}, \cdots, E_{2n-1}, E_0) \Kb_0 \Kb_1 \Kb_{2n-1} \cdots \Kb_k \Kb_{2n-k}, 
\end{align}
for $1 \leq k \leq n-1$. We first prove, by induction, that $P_k \equiv Q_k$ modulo $\Uu_{d_{i} \geq 1,+}$. The case $k=1$ follows from Lemma \ref{lem: PBBE}. Assuming the claim holds for $k-1$, we get 
\begin{align*}
P_{k} = {\bf P}(B_{2n-k}, B_k, P_{k-1}) \equiv {\bf P}(B_{2n-k}, B_k, Q_{k-1}) \equiv  [F_{2n-k}, [F_k, Q_{k-1}]_q]_q = Q_k. 
\end{align*}

Let us now prove that $R_k = S_k$. As before, the case $k=1$ follows from Lemma \ref{lem: PBBE}. 
Assuming that the claim holds for $k-1$, we get 
\begin{align*}
R_{k} = {\bf P}(B_{2n-k}, B_k, R_{k-1}) \equiv {\bf P}(B_{2n-k}, B_k, S_{k-1}), 
\end{align*} 
and an easy calculation shows that 
\begin{align*}
{\bf P}(B_{2n-k}, B_k, S_{k-1}) &= [F_{2n-k}, [E_{2n-k} \Kb_k, S_{k-1}]_q]_q + [E_{k}\Kb_{2n-k}, [F_{k},  S_{k-1}]_q]_q \\ 
&\quad + [E_{k} \Kb_{2n-k}, [E_{2n-k} \Kb_k, S_{k-1}]_q]_q - q \Ka_{2n-k} \Kb_k S_{k-1} = S_k.  
\end{align*} 
It follows that 
\[ 
\Tbr_{\varpi'_n}(B_n) = P_{n-1} + R_{n-1} \equiv Q_{n-1} + S_{n-1} = T_{\omega_n'}(F_n) + q^{2n-2}T_{\omega_n'}(E_n) T_{\omega_n'}(\Kb_n),
\]
completing the proof. 
\enproof 

\subsection{Type ${\sf AIII}_{2n}^{(\tau)}$}

We will now prove analogous results for $N=2n$. 

\Lem \label{lem: T B pol C even case} 
Let $n \geq 2$. 
The expression $\Tbr_{\varpi'_i}(B_i)$ is equal to
\eq \label{eq: C good pol even case}  
\begin{cases}
P(B_{i-1}, \ldots, B_1, B_{i+1},B_{i+2}, \ldots, B_{n-1}, {\bf P}( B_n,B_{n+1}, P(B_{n+2}, \ldots, B_{2n},B_0 ) )  ) & \\
 \mbox{ if } 1 \leq i < n, \\
P( B_{i+1},\ldots,B_{2n},B_{i-1},B_{i-2},\ldots,B_{n+2}, {\bf P}( B_{n+1},B_n,P(B_{n-1},\ldots, B_1,B_0) ) ) & \\
 \mbox{ if } n+1 < i \le 2n, \\
P( B_i, {\bf P}(B_{n-1}, B_{n+2}, {\bf P}(B_{n-2}, B_{n+3}, {\bf P}( \ldots {\bf P}(B_1, B_{2n}, B_0) )) )) & \\ 
 \mbox{ if } i \in \{n,n+1\}.
\end{cases}  
\eneq 
\enlem

\begin{proof}
Let $1 \leq i<n$. By Proposition \ref{prop:reduced expression of omegai in extended affine Weyl groups type AIIItau}, one has $\varpi_i' = \zeta_i \tau_i$, where $\zeta_i=(r_0 [1,n])^i$ and $\tau_i=[n{-}i, n{-}1] \cdots [2,i{+}1][1,i{-}1]$. 
By \cite[Lemma 5.3]{LP25}, 
\begin{align*}
\Tbr_{\tau_i}(B_i) = P(B_{n-1}, B_{n-2}, \ldots, B_{n-i+1}, P(B_1, \ldots, B_{n-i})),
\end{align*}
and $\Tbr_{\zeta_i}(B_j) = B_{(j+i) \: {\rm mod} \: n}$ if $j \ne n-i$. Therefore
\begin{align*}
\Tbr_{\varpi'_i}(B_i) & = P(B_{i-1}, B_{i-2}, \ldots, B_{1}, P(B_{i+1}, B_{i+2}, \ldots, B_{n-1}, \Tbr_{\zeta_i}(B_{n-i}))) \\
& = P(B_{i-1}, B_{i-2}, \ldots, B_{1}, B_{i+1}, B_{i+2}, \ldots, B_{n-1}, \Tbr_{\zeta_i}(B_{n-i})).
\end{align*}
Since, by Lemma \ref{lem: TjTiBj formula 2}, 
\begin{align*}
\Tbr_{\zeta_i}( B_{n-i} ) & = \Tbr_{r_0[1,n]}( B_{n-1} ) = {\bf P}( B_n,B_{n+1}, \Tbr_{r_0[1,n-2]}(B_{n+2}) ) \\
& =  {\bf P}(B_n, B_{n+1}, P(B_{n+2}, B_{n+3}, \ldots, B_{2n}, B_0),
\end{align*}
the result follows. The case of $i \in [n+2, 2n]$ is proved similarly.

Let $i=n$. Then $\varpi_i' = \zeta'_n \tau_n$, where $\zeta_n'=(r_0[1,n])^{n-1}$, $\tau_n=r_1 \cdots r_{n-1}$. We have
\begin{align*}
{\bf T}_{\tau_n}(B_n) & = P(B_n, B_{n-1}, \ldots, B_1, B_0) \\
& = P(P( B_n, B_{n-1},\ldots,B_1 ), B_0) = P( {\bf T}_{[1,n-1]}(B_n), B_0 ),
\end{align*}
where we have used the fact that $B_0$ commutes with $B_2,\ldots,B_n$. Therefore
\begin{align*}
{\bf T}_{\zeta_n'\tau_n}(B_n) & = {\bf T}_{(r_0[1,n])^{n-1}}(P( {\bf T}_{[1,n-1]}(B_n), B_0 )) \\
& = P(B_n, {\bf T}_{(r_0[1,n])^{n-1}}( B_0 ) ) \\
& = P(B_n, {\bf P}( B_{n-1}, B_{n+2}, {\bf T}_{(r_0[1,n])^{n-2}}(B_0) ) ) \\
& = P( B_n, {\bf P}(B_{n-1}, B_{n+2}, {\bf P}(B_{n-2}, B_{n+3}, {\bf P}( \ldots {\bf P}(B_1, B_{2n}, B_0) )) )).
\end{align*}
The case of $i=n+1$ is proved similarly. 
\end{proof}

\Pro \label{pro: Aeve igood}
For $1 \le i \le 2n$, 
\begin{align*}
{\bf T}_{\varpi'_i}(B_i) \equiv \begin{cases} T_{\omega'_i}(F_i) - q^{2n-1}T_{\omega'_{\tau(i)}}(E_{\tau(i)})T_{\omega'_i}(\widetilde{K}'_i) \ \ \ \mod \Uu_{d_{i} \geq 1,+},   & i \not\in \{n,n+1\}, \\[3pt]
T_{\omega'_{\tau(i)}}(F_{\tau(i)}) - q^{2n}T_{\omega'_{i}}(E_{i})T_{\omega'_{\tau(i)}}(\widetilde{K}'_{\tau(i)}) \ \mod \Uu_{d_{\tau(i)} \geq 1,+},  & i \in \{n,n+1\}.
\end{cases}
\end{align*}
\enpro

\Proof
The proof is analogous to that of Proposition \ref{pro: Aodd igood}. 
\enproof

\section{Factorization formula} 
\label{sec: factor formula}

We consider $\Ui$ as a subalgebra of $\Uu$ under the embedding $\eta$ from \eqref{eq: Kolb emb}. Most of the time, we will omit $\eta$ from notation. In this section we prove the main result of the paper, i.e., a factorization formula in arbitrary rank. 

\Thm \label{thm: main overall} 
The `Drinfeld--Cartan' series $\thvar{i}$ admits, for each $i \in \indx_0$, the following factorization: 
\begin{align*}
\thvar{i} \equiv K_i K_{\tau(i)}\mi   \boldsymbol\phi^-_i(z\mi)\boldsymbol{\phi}^+_{\tau(i)}(\CCC z)  \mod \fext{\Uq_+}{z}. 
\end{align*}
\enthm

There are three cases to be considered, depending on whether $a_{\tau(i)} = 0, 2$ or $-1$. Even though the factorization formula above holds uniformly in all the three cases, the proofs vary slightly from case to case. 

\subsection{First estimate} 

It is easily verified that 
\[
\Tbr_{\varpi'_i}(\KK_i) = 
\begin{cases}
q^{-2} \CCC \KK_i\mi & \mbox{if } N \mbox{ is odd } \& \ a_{i,\tau(i)} = 0, \\ 
\CCC \KK_i\mi & \mbox{if } N \mbox{ is odd } \& \ a_{i,\tau(i)} = 2, \\ 
-q^{-2} \CCC \KK_i\mi & \mbox{if } N \mbox{ is even } \& \ a_{i,\tau(i)} = 0, \\ 
-q^{-1} \CCC \KK_{\tau(i)}\mi & \mbox{if } N \mbox{ is even } \& \ a_{i,\tau(i)} = -1, 
\end{cases}
\]
and
\[
T_{\omega'_i}(\Ka_i) = c_\delta^{1/2} \Ka_i\mi, \qquad T_{\omega'_i}(\Kb_i) = c_\delta^{1/2} (\Kb_i)\mi. 
\]

The statement of the following lemma should also be seen as the definition of the remainder term $Q_i$. 

\Lem \label{lem: A-1 gen}
For each $i \in \indx_0$: 
\[
A_{i,-1} = \big( \CCC\mi c_\delta^{1/2} \Ka_{i} x^+_{\tau(i),-1} + x^-_{i,1} \big) + Q_i, \qquad Q_i \in \Uu_{d_{\tau(i)} \geq 1,+}.   
\] 
In the special case $N=1$, the remainder term $Q_i$ is zero. 
\enlem

\Proof
By \cite[Proposition 4.7]{beck-94}, 
\begin{align} 
\label{eq: x+beck}
\ \ \ \ \  o(\tau(i)) x_{\tau(i),-1}^+ =& \ T_{\omega_{\tau(i)}}(E_{\tau(i)}) = T_{\omega'_{\tau(i)}}(- F_{\tau(i)}(\Kb_{\tau(i)})\mi) = -c_\delta^{-1/2} T_{\omega'_{\tau(i)}}(F_{\tau(i)})  \Kb_{\tau(i)}, \\
\label{eq: x-beck}
o(i) x_{i,1}^- =& \ T_{\omega_i}(F_i) = T_{\omega'_i}(- \Ka_i\mi E_i) = - c_\delta^{-1/2} \Ka_i T_{\omega'_i}(E_i). 
\end{align} 
Moreover, by \cite[Theorem 4.5]{LWZ-quasi} and \cite[Theorem 4.7]{LPWZ}, $o(i) A_{i,-1} = \Tbr_{\varpi_i}(B_i)$. 
Let us consider the four possible cases separately. 

If $N$ is odd and $a_{i,\tau(i)} = 0$, then 
\[
o(i) A_{i,-1} = \Tbr_{\varpi'_i}(-\KK_i\mi B_{\tau(i)}) = -q^2\CCC\mi \KK_i \Tbr_{\varpi'_{\tau(i)}}(B_{\tau(i)}). 
\]
Hence, by Proposition \ref{pro: Aodd igood}, as well as \eqref{eq: x+beck}--\eqref{eq: x-beck}, we get
\begin{align}
A_{i,-1} &= - o(i) q^2\CCC\mi \KK_i \big(  T_{\omega_{\tau(i)}'}(F_{\tau(i)}) + q^{2n} T_{\omega'_{i}}(E_{i}) T_{\omega_{\tau(i)}'}(\Kb_{\tau(i)})\big) + Q_i \\ 
&= - q^{2}o(i) \CCC\mi \Ka_i \Kb_{\tau(i)} \big( T_{\omega_{\tau(i)}'}(F_{\tau(i)}) + q^{2n} c_\delta^{1/2} T_{\omega'_{i}}(E_{i}) (\Kb_{\tau(i)})\mi \big) + Q_i \\ 
 &= \big( \CCC\mi c_\delta^{1/2} \Ka_{i} x^+_{\tau(i),-1}  + x^-_{i,1} \big) + Q_i, 
\end{align} 
for some $Q_i \in \Uu_{d_{\tau(i)} \geq 1,+}$. 
If $N$ is even and $a_{i,\tau(i)} = 0$, then a similar calculation using Proposition \ref{pro: Aeve igood} also shows that 
\[
A_{i,-1} = \big( \CCC\mi c_\delta^{1/2} \Ka_{i} x^+_{\tau(i),-1}  + x^-_{i,1} \big) + Q_i, 
\]
for some $Q_i \in \Uu_{d_{\tau(i)} \geq 1,+}$. 

If $N$ is odd and $a_{i,\tau(i)} = 2$, then 
\[
o(i) A_{i,-1} = \Tbr_{\varpi'_i}(\KK_i\mi B_{i}) = \CCC\mi \KK_i \Tbr_{\varpi'_{i}}(B_{i}). 
\]
Hence, again by Proposition \ref{pro: Aodd igood}, we get 
\begin{align}
A_{i,-1} &= o(i) \CCC\mi \KK_i \big(  T_{\omega_{{i}}'}(F_{{i}}) + q^{2n-2} T_{\omega'_{i}}(E_{i}) T_{\omega_{{i}}'}(\Kb_{{i}})\big) + Q_i \\ 
&= - o(i) q^{2} \Ka_i \Kb_{i} \CCC\mi \big( T_{\omega_{{i}}'}(F_{{i}}) + q^{2n-2} c_\delta^{1/2} T_{\omega'_{i}}(E_{i}) (\Kb_{{i}})\mi \big) + Q_i \\ 
 &= \big(\CCC\mi c_\delta^{1/2} \Ka_{i} x^+_{{i},-1}  +  x^-_{i,1} \big) + Q_i, 
\end{align} 
for some $Q_i \in \Uu_{d_{i} \geq 1,+}$. 

Finally, let $N$ be even and $a_{i,\tau(i)} = -1$. Then 
\[
o(i) A_{i,-1} = \Tbr_{\varpi'_i}(-q^{-2} B_{i}\KK_{\tau(i)}\mi) = q\mi \CCC\mi  \Tbr_{\varpi'_{i}}(B_{i}) \KK_i. 
\]
Therefore, by Proposition \ref{pro: Aeve igood}, we get 
\begin{align}
A_{i,-1} &= o(i) q\mi \CCC\mi  \big(  T_{\omega_{\tau(i)}'}(F_{\tau(i)}) - q^{2n} T_{\omega'_{i}}(E_{i}) T_{\omega_{\tau(i)}'}(\Kb_{\tau(i)})\big) \KK_i + Q_i \\ 
&= o(i) q^{-1} \CCC\mi  \big( T_{\omega_{\tau(i)}'}(F_{\tau(i)}) - q^{2n} c_\delta^{1/2} T_{\omega'_{i}}(E_{i}) (\Kb_{\tau(i)})\mi \big) \Ka_i \Kb_{\tau(i)} + Q_i \\ 
 &= \big( \CCC\mi c_\delta^{1/2} \Ka_{i} x^+_{\tau(i),-1} + x^-_{i,1} \big) + Q_i, 
\end{align} 
for some $Q_i \in \Uu_{d_{\tau(i)} \geq 1,+}$. 
\enproof

\subsection{The case: $a_{i, \tau(i)} = 0$} 

We first need an estimate of $H_{i,1}$. 

\Lem \label{lem: Hi1 aij0}
We have 
\[
H_{i,1} \equiv \CCC h_{\tau(i), 1} - h_{i,-1} - c_\delta^{1/2} (q^2 - q^{-2})x_{\tau(i),0}^+ x_{i,-1}^+ +  \CCC \KK_{\tau(i)}\mi [ x_{i,0}^- , \ Q_{\tau(i)} ] 
\quad \mod \Uu_{d_i, d_{\tau(i)} \geq 1,+}. 
\]
\enlem 

\Proof
Using Lemma \ref{lem: A-1 gen}, one checks that 
\begin{align}
H_{i,1} &=  \CCC \KK_{\tau(i)}\mi [ A_{i,0}, A_{\tau(i),-1}] \\ 
&= \CCC \KK_{\tau(i)}\mi [ x_{i,0}^- + x_{\tau(i),0}^+ \Kb_i , \ x^-_{\tau(i),1} + \CCC\mi c_\delta^{1/2}\Ka_{\tau(i)} x^+_{i,-1}  + Q_{\tau(i)} ] \\ 
&\equiv \big( \CCC h_{\tau(i), 1} - h_{i,-1} \big) 
- c_\delta^{1/2} (q^2 - q^{-2})x_{\tau(i),0}^+ x_{i,-1}^+  + \CCC \KK_{\tau(i)}\mi [ x_{i,0}^- , \ Q_{\tau(i)} ] 
\end{align}
modulo $\Uu_{d_i, d_{\tau(i)} \geq 1,+}$. 
\enproof

Let $\avarm{i}$ and $\avarp{i}$ denote the homogeneous components of $\avar{i}$ of degrees $- \alpha_i$ and $\alpha_{\tau(i)}$, respectively. 

\Pro \label{pro: A_i aij0}
We have 
\[
\avar{i} \equiv \avarm{i} + \avarp{i} \quad \mod \fext{\Uu_{d_{\tau(i)} \geq 1,+}}{z},
\]
with
\begin{align*}
\avarm{i} = \mathbf{x}^-_{i,-}(z\mi), \qquad \avarp{i} = c_i^{1/2} \mathbf{x}_{\tau(i)}^+(\CCC c_\delta^{-1/2} z) \boldsymbol\phi^-_{i, -}(z\mi). 
\end{align*} 
\enpro 

\Proof
According to Lemma \ref{lem: A Th as series}, 
\[ 
x^-_{i,0} + x^+_{\tau(i),0} \Kb_i = A_{i,0} = (1 - \Ad_{\overline{H}_{i,1}} z) \avar{i}. 
\] 
Applying Lemma \ref{lem: Hi1 aij0}, we deduce that, modulo $\fext{\Uu_{d_{\tau(i)} \geq 1,+}}{z}$, we have $\avar{i} \equiv \avarm{i} + \avarp{i} + \mathbf{R}_i(z)$, with 
\begin{align} 
x^-_{i,0} &= (1 + \ad_{\bar{h}_{i,-1}} z) \avarm{i}, \\
x^+_{\tau(i),0} \Kb_i &= (1 - \CCC \Ad_{\bar{h}_{\tau(i),1}} z) \avarp{i} + \rho c_\delta^{1/2} x_{\tau(i),0}^+ [x_{i,-1}^+, \avarm{i}] z, \\ \label{eq: Rz a0}
\mathbf{R}_i(z) &= [2]\mi \CCC [  \KK_{\tau(i)}\mi [ x_{i,0}^- , \ Q_{\tau(i)} ], \avarm{i}] z. 
\end{align} 
The first identity above implies immediately that $\avarm{i} = \mathbf{x}^-_{i,-}(z\mi)$. Therefore, the second identity can be rearranged to
\begin{align*}
\avarp{i} &= \left(1 - \CCC \Ad_{\bar{h}_{\tau(i),1}} z \right)\mi \left( x^+_{\tau(i),0} \Kb_i 
-  \rho c_\delta^{1/2}  x_{\tau(i),0}^+ [x_{i,-1}^+, \mathbf{x}^-_{i,-}(z\mi)] z \right) \\ 
&= \left(1 - \CCC \Ad_{\bar{h}_{\tau(i),1}} z \right)\mi \left( x^+_{\tau(i),0} \Kb_i 
+  c_i^{1/2} x_{\tau(i),0}^+ \boldsymbol\phi^-_{i, \leq -1}(z\mi) \right) \\ 
&= c_i^{1/2} \left(1 - \CCC \Ad_{\bar{h}_{\tau(i),1}} z \right)\mi x_{\tau(i),0}^+ \boldsymbol\phi^-_{i, -}(z\mi) \\ 
&= c_i^{1/2} \mathbf{x}_{\tau(i)}^+(\CCC c_\delta^{-1/2} z) \boldsymbol\phi^-_{i, -}(z\mi). 
\end{align*}
To complete the proof, one needs to show that $\mathbf{R}_i(z) \in \fext{\Uu_{d_{\tau(i)} \geq 1,+}}{z}$. This is sketched in Appendix \ref{sec: R(z) est} below. 
\enproof

\Thm
The `Drinfeld--Cartan' series $\thvar{i}$ admits, for each $i \in \indx_0$, the following factorization: 
\begin{align*}
\thvar{i} \equiv K_i K_{\tau(i)}\mi   \boldsymbol\phi^-_i(z\mi)\boldsymbol{\phi}^+_{\tau(i)}(\CCC z)  \mod \fext{\Uq_+}{z}. 
\end{align*}
\enthm

\Proof
According to Lemma \ref{lem: A Th as series} and Proposition \ref{pro: A_i aij0},  
\begin{align*}
\thvar{i} &= \KK_i \KK_{\tau(i)}\mi + \rho\KK_{\tau(i)}\mi [\avar{i}, A_{\tau(i),0}] \\ 
&\equiv  \KK_i \KK_{\tau(i)}\mi + \rho\KK_{\tau(i)}\mi [\mathbf{x}^-_{i,-}(z\mi) + c_i^{1/2} \mathbf{x}_{\tau(i)}^+(\CCC c_\delta^{-1/2} z) \boldsymbol\phi^-_{i, -}(z\mi), \ x^-_{\tau(i),0} + x^+_{i,0}\Kb_{\tau(i)}] \\ 
&\equiv \KK_i \KK_{\tau(i)}\mi + \rho\KK_{\tau(i)}\mi \left(\Kb_{\tau(i)} [\mathbf{x}^-_{i,-}(z\mi), x^+_{i,0}] + c_i^{1/2}\boldsymbol\phi^-_{i, -}(z\mi)[ \mathbf{x}_{\tau(i)}^+(\CCC c_\delta^{-1/2} z)  ,  x^-_{\tau(i),0}] \right) \\
&= c_i\KK_{\tau(i)}\mi\left(  K_i K_{\tau(i)}\mi  + K_{\tau(i)}\mi (\boldsymbol\phi^-_{i, -}(z\mi) - K_i) + \boldsymbol\phi^-_{i, -}(z\mi) (\boldsymbol\phi^+_{\tau(i)}(\CCC z) - K_{\tau(i)}\mi) \right) \\ 
&= K_i K_{\tau(i)}\mi   \boldsymbol\phi^-_i(z\mi)\boldsymbol{\phi}^+_{\tau(i)}(\CCC z), 
\end{align*}
as claimed. 
\enproof

\subsection{The case: $a_{i, \tau(i)} = 2$} 

This case only occurs when $N$ is odd and $i=n$. 
Again, we first need an estimate of $H_{n,1}$. 

\Lem \label{lem: Hi1 aij2} 
We have 
\[
H_{n,1} \equiv \big( \CCC h_{n, 1} - h_{n,-1} \big) + q^{-2} c_\delta^{1/2}[x_{n,0}^+, x_{n,-1}^+]_{q^6} - \CCC \KK_{n}\mi [ x_{n,0}^- , \ Q_{n} ]_{q^2} 
\quad \mod \Uu_{d_n \geq 2,+}. 
\]
\enlem 

\Proof
Using Lemma \ref{lem: A-1 gen}, one checks that 
\begin{align}
H_{n,1} &= -\CCC \KK_n\mi  [A_{n,0}, A_{n,-1}]_{q^{2}} \\ 
&= -\CCC \KK_{n}\mi [ x_{n,0}^- + x_{n,0}^+ \Kb_n , \ x^-_{n,1} +  \CCC\mi c_\delta^{1/2} \Ka_{n} x^+_{n,-1}  + Q_{n} ]_{q^2} \\ 
&\equiv \big( \CCC h_{n, 1} - h_{n,-1} \big) + q^{-2} c_\delta^{1/2}[x_{n,0}^+, x_{n,-1}^+]_{q^6} - \CCC \KK_{n}\mi [ x_{n,0}^- , \ Q_{n} ]_{q^2} 
\end{align}
modulo $\Uu_{d_n \geq 2,+}$. 
\enproof

Consider the \emph{non-commutative} diagram

\eq \label{eq: comm diagram compatible}
\begin{tikzcd}
\widetilde{\mathcal{O}}_q \arrow[r, "\eta", hookrightarrow] \arrow[d, "\wr"'] & \Uu(\widehat{\mathfrak{sl}}_2) \arrow[d,  "\iota_n", "\wr"']   \\
\Ui_{[n]} \arrow[d, hookrightarrow]  & \Uu_{[n]} \arrow[d, hookrightarrow]  \\
\Ui \arrow[r, hookrightarrow] & \Uu. 
\end{tikzcd} 
\eneq 
Let $\avarm{n}$ and $\avarp{n}$ denote the homogeneous components of $\avar{n}$ of degrees $- \alpha_n$ and $\alpha_{n}$, respectively. Here $\avar{n}$ should be considered as an element of $\fext{\Ui}{z}$ and $\fext{\Uu}{z}$ (in the bottom row of the diagram). On the other hand, series without subscripts, e.g., $\mathbf{A}_+(z)$ and $\boldsymbol{\grave{\Theta}}_+(z)$, are elements of $\fext{\widetilde{\mathcal{O}}_q}{z}$ (the upper left corner). 

\Thm \label{thm: main aij2}
We have 
\begin{align*}
\avar{n} &\equiv \iota_n \circ \eta (\mathbf{A}_+(z)) \quad \mod  \fext{\Uu_{d_n \geq 1,+}}{z}, \\ 
\thvar{n} &\equiv \iota_n \circ \eta (\boldsymbol{\grave{\Theta}}_+(z)) \quad \mod  \fext{\Uq_{+}}{z}. 
\end{align*}
Hence 
\begin{align*}
\thvar{n} \equiv \boldsymbol\phi^-_n(z\mi)\boldsymbol{\phi}^+_{n}(\CCC z)  \mod \fext{\Uq_+}{z}. 
\end{align*}
\enthm 

\Proof
According to Lemma \ref{lem: A Th as series}, 
\[
A_{n,0} + \CCC z A_{n,-1} = 
(1 - \Ad_{\overline{H}_{n,1}} z - \CCC z^2) \avar{n}.
\]
Applying Lemmas \ref{lem: A-1 gen} and \ref{lem: Hi1 aij2}, we deduce that,  
modulo $\fext{\Uu_{d_{n} \geq 1,+}}{z}$, we have $\avar{n} \equiv \avarm{n} + \avarp{n} + \mathbf{R}_n(z)$, with 
\begin{align} 
x^-_{n,0} + \CCC x_{n,1}^-z &= (1 - \CCC \ad_{\bar{h}_{n,1}} z)(1 + \ad_{\bar{h}_{n,-1}} z) \avarm{n}, \\
x^+_{n,0} \Kb_n + c_\delta^{1/2} \Ka_n x^+_{n,-1} z &=  (1 - \CCC \ad_{\bar{h}_{n,1}} z)(1 + \ad_{\bar{h}_{n,-1}} z) \avarp{n} \\
&\quad - q^{-2}[2]\mi c_\delta^{1/2}[[x_{n,0}^+, x_{n,-1}^+]_{q^6}, \avarm{n}] z, \\ \label{eq: Rz a2}
(1-\CCC z^2) \mathbf{R}_n(z) &= [2]\mi \CCC \KK_n\mi [ [ x_{n,0}^- , \ Q_{n} ]_{q^2}  , \avarm{n}] z. 
\end{align} 
The first identity above implies immediately that $\avarm{n} = \mathbf{x}^-_{n,-}(z\mi)$. Moreover, since the first two identities only involve generators indexed by $n$, we conclude that $\avarm{n} = \iota_n(\mathbf{A}_+^{(-)})$ and $\avarp{n} = \iota_n(\mathbf{A}_+^{(+)})$. Therefore, it suffices to show that $\mathbf{R}_n(z)~=~0$. 
This is sketched in Appendix \ref{sec: R(z) est} below. 

Moving on to the second equivalence, from Lemma \ref{lem: A Th as series} we get 
\[
\thvar{n} = 1 + \frac{q^2\rho\CCC \KK_n\mi z^2}{1 - \CCC z^2}\left(z\mi[A_{n,-1}, \avar{n}]_{q^{-2}} - q^{-2}[A_{n,0}, \avar{n}]_{q^2} \right).  
\]
We can write $\avar{n} = \iota_n \circ \eta (\mathbf{A}_+(z)) + \mathbf{S}(z)$ for some $\mathbf{S}(z) \in \fext{\Uq_{d_n \geq 1,+}}{z}$. It is clear from \eqref{eq: Kolb emb} and Lemma \ref{lem: A-1 gen} that $[Q_i, \avar{n}]_{q^{-2}}$, $[A_{n,-1}, \mathbf{S}(z)]_{q^{-2}}$ and $[A_{n,0}, \mathbf{S}(z)]_{q^2}$ lie in $\fext{\Uq_+}{z}$. This implies $\thvar{n} \equiv \iota_n \circ \eta (\boldsymbol{\grave{\Theta}}_+(z))$. Since the Drinfeld positive half of $\Uq(\widehat{\mathfrak{sl}}_2)$ 
is mapped by $\iota_n$ into the Drinfeld positive half of $\Uq$, the equivalence $\thvar{n} \equiv \boldsymbol\phi^-_n(z\mi)\boldsymbol{\phi}^+_{n}(\CCC z)$ 
now follows from Theorem \ref{thm: rank 1 factorization}. 
\enproof

\subsection{The case: $a_{i, \tau(i)} = -1$} 

This case only occurs when $N$ is even and $i \in \{ n, n+1 \}$. 
Again, we first need an estimate of $H_{i,1}$. The calculation is somewhat harder than in the preceding two cases. 
In particular, we will need a finer decomposition of the remainder term $Q_{\tau(i)} = Q_{\tau(i)}' + Q_{\tau(i)}'' + Q_{\tau(i)}'''$ from Lemma \ref{lem: Qi a-1}. 

\Lem \label{lem: Hi1 aij-1} 
We have 
\[
H_{i,1} \equiv \big( \CCC h_{\tau(i), 1} - h_{i,-1} \big) + \mathcal{Y}_i + q\CCC \KK_{\tau(i)}\mi [ x_{i,0}^- , \ Q''_{\tau(i)} ]_{q\mi}
\quad \mod \Uu_{d_i, d_{\tau(i)} \geq 1,+}. 
\]
with 
\[
\mathcal{Y}_i = c_\delta^{1/2} \left( (1-q^{-2}) x^+_{\tau(i), -1} x_{i,0}^+ + (q^{-2} - q)  x_{\tau(i),0}^+x^+_{i, -1} +  (1-q) x^+_{i, -1} x_{\tau(i),0}^+ \right). 
\]
\enlem 

\Proof
Let us first estimate the term $\mathbf{T}\mi_{\boldsymbol\theta_n}(B_0)$. As shown in, e.g., the proof of \cite[Proposition 3.4]{LPWZ}, 
\[
\mathbf{T}\mi_{\boldsymbol\theta_n}(B_0) = \mathbf{P}(P(B_{n+2}, \cdots B_{N}), P(B_{n-1}, \cdots, B_1), B_0). 
\]
An argument similar to those in \S \ref{sec: compat of br} shows that 
\[
\mathbf{T}\mi_{\boldsymbol\theta_n}(B_0) \equiv [ P(F_{n-1}, \cdots, F_1), P(F_{n+2}, \cdots, F_N, F_0)   ]_q \quad \mod \Uu_{d_i, d_{\tau(i)} \geq 1,+}. 
\]
Let us abbreviate the RHS above as $Y$. 
On the other hand, by, e.g., \cite[Lemma 5.3]{LP25}, 
\[
x_{i,-1}^+ = o(i) T_{\omega_i}(E_i) = -o(i) c_\delta^{-1/2} [ P(F_{i-1}, \cdots, F_1), P(F_{i+1}, \cdots, F_N, F_0)   ]_q  \Kb_i. 
\]
Therefore, 
\[
[x_{i,-1}^+, x_{\tau(i),0}^+]_q = o(i)q c_\delta^{-1/2} [ P(F_{n-1}, \cdots, F_1), P(F_{n+2}, \cdots, F_N, F_0)   ]_q \Kb_i \Kb_{\tau(i)}.
\]
By Lemma \ref{lem: Qi a-1}, we have $Q'_{\tau(i)} = -o(i) q\mi \CCC\mi [E_i \Kb_{\tau(i)}, Y]_q \KK_{\tau(i)}$. Then
\begin{align*}
 q \CCC \KK_{\tau(i)}\mi [F_i, Q'_{\tau(i)}]_{q\mi} &= -o(i) [F_i, [E_i , Y ]_q ]_{q^2} \\
 &= \rho o(i) [F_i, E_i Y ]_{q} \Kb_{\tau(i)} \\ 
 &=  \rho o(i) \big( [F_i, E_i]Y + E_i [F_i, Y]_q \big) \Kb_{\tau(i)} \\ 
 &= -o(i)Y(q \Ka_i - q\mi\Kb_i)\Kb_{\tau(i)} + \rho c_\delta^{1/2} x_{i,0}^+ x_{\tau(i),-1}^+ \\ 
 &= -o(i) q \mathbf{T}\mi_{\boldsymbol\theta_n}(B_0) \KK_i + q^{-2} c_\delta^{1/2} [x_{i,-1}^+, x_{\tau(i),0}^+]_q  + \rho c_\delta^{1/2} x_{i,0}^+ x_{\tau(i),-1}^+. 
\end{align*} 
Using relation \eqref{eq: rel sln xx}, we get
\[
q \CCC \KK_{\tau(i)}\mi [F_i, Q'_{\tau(i)}]_{q\mi} + o(i) q \mathbf{T}\mi_{\boldsymbol\theta_n}(B_0) \KK_i 
= (1 - q^{-2}) x^+_{\tau(i),-1} x^+_{i,0} + [x^+_{i,-1}, x^+_{\tau(i),0}]_q. 
\]

Using the definition of $H_{i,1}$ (see \cite[(4.7)]{LPWZ}) and the above calculations, one now checks that 
\begin{align*}
H_{i,1} &= q\CCC \KK_{\tau(i)}\mi  [A_{i,0}, A_{\tau(i),-1}]_{q\mi}  + o(i) q \mathbf{T}\mi_{\boldsymbol\theta_n}(B_0) \KK_i \\ 
&\equiv q \CCC \KK_{\tau(i)}\mi \left[ x_{i,0}^- + x_{\tau(i),0}^+ \Kb_i , \ x^-_{\tau(i),1} + \CCC\mi c_\delta^{1/2}\Ka_{\tau(i)} x^+_{i,-1}  + Q'_{\tau(i)} + Q''_{\tau(i)} \right]_{q\mi} \\
&\quad  + o(i) q \mathbf{T}\mi_{\boldsymbol\theta_n}(B_0) \KK_i \\ 
&\equiv \big( \CCC h_{\tau(i), 1} - h_{i,-1} \big) + \mathcal{Y}_i + q\CCC \KK_{\tau(i)}\mi [ x_{i,0}^- , \ Q''_{\tau(i)} ]_{q\mi}.  
\end{align*}
modulo $\Uu_{d_i, d_{\tau(i)} \geq 1,+}$. 
\enproof

Consider the \emph{non-commutative} diagram

\eq \label{eq: comm diagram compatible}
\begin{tikzcd}
\Ui(\widehat{\mathfrak{sl}}_3, \tau) \arrow[r, "\eta", hookrightarrow] \arrow[d, "\wr"'] & \Uu(\widehat{\mathfrak{sl}}_3) \arrow[d,  "\iota_{n,n+1}", "\wr"']   \\
\Ui_{[n]} \arrow[d, hookrightarrow]  & \Uu_{[n,n+1]} \arrow[d, hookrightarrow]  \\
\Ui \arrow[r, hookrightarrow] & \Uu. 
\end{tikzcd} 
\eneq 
Given $i \in \{n,n+1\}$, let $\bar{i} = i-n+1$. 
Let $\avarm{i}$ and $\avarp{i}$ denote the homogeneous components of $\avar{i}$ of degrees $- \alpha_i$ and $\alpha_{\tau(i)}$, respectively. 
Here $\avar{i}$ should be considered as an element of $\fext{\Ui}{z}$ and $\fext{\Uu}{z}$ (in the bottom row of the diagram). On the other hand, series with barred subscripts, e.g., $\mathbf{A}_{\bar{i},+}(z)$ and $\boldsymbol{\grave{\Theta}}_{\bar{i},+}(z)$, are elements of $\fext{\Ui(\widehat{\mathfrak{sl}}_3, \tau)}{z}$ (the upper left corner).

\Thm \label{thm: main aij2}
For $i \in \{n,n+1\}$, we have 
\begin{align*}
\avar{i} &\equiv \iota_{n,n+1} \circ \eta (\mathbf{A}_{\bar{i},+}(z)) \quad \mod  \fext{\Uu_{d_{\tau(i)} \geq 1,+}}{z}, \\ 
\thvar{i} &\equiv \iota_{n,n+1} \circ \eta (\boldsymbol{\grave{\Theta}}_{\bar{i},+}(z)) \quad \mod  \fext{\Uq_{+}}{z}. 
\end{align*}
Hence 
\begin{align*}
\thvar{i} \equiv K_i K_{\tau(i)}\mi  \boldsymbol\phi^-_i(z\mi)\boldsymbol{\phi}^+_{\tau(i)}(\CCC z)  \mod \fext{\Uq_+}{z}. 
\end{align*}
\enthm 

\Proof
According to Lemma \ref{lem: A Th as series}, 
\[
[2] A_{i,0} - \CCC z A_{i,-1} = 
([2] - \Ad_{H_{i,1}} z + \CCC z^2) \avar{i}.
\]
Applying Lemmas \ref{lem: A-1 gen} and \ref{lem: Hi1 aij-1}, we deduce that,  
modulo $\fext{\Uu_{d_{\tau(i)} \geq 1,+}}{z}$, we have $\avar{i} \equiv \avarm{i} + \avarp{i} + \mathbf{R}_i(z)$, with 
\begin{align} 
[2]x^-_{i,0} - \CCC x_{i,1}^-z &= [2] (1 + \Ad_{\bar{h}_{i,-1}} z) (1 - \Ad_{\bar{h}_{\tau(i),1}} \CCC z) \avarm{i}, \\
[2]x^+_{\tau(i),0} \Kb_{i} - c_\delta^{1/2} \Ka_i x^+_{\tau(i),-1} z &=  [2] (1 + \Ad_{\bar{h}_{i,-1}} z) (1 - \Ad_{\bar{h}_{\tau(i),1}} \CCC z) \avarp{i} \\
&\quad - [\mathcal{Y}_i, \avarm{n}] z, \\ \label{eq: Rz a-1}
([2]+\CCC z^2) \mathbf{R}_i(z) &= \CCC [ \KK_{\tau(i)}\mi [ x_{i,0}^- , \ Q''_{\tau(i)} ]_{q\mi}  , \avarm{i}] z. 
\end{align} 
The first identity above implies immediately that $\avarm{i} = \mathbf{x}^-_{i,-}(z\mi)$. Moreover, since the first two identities only involve generators indexed by $i$ and $\tau(i)$, we conclude that $\avarm{i} = \iota_{n,n+1}(\mathbf{A}_{\bar{i},+}^{(-)})$ and $\avarp{i} = \iota_{n,n+1}(\mathbf{A}_{\bar{i},+}^{(+)})$. Therefore, it suffices to show that $\mathbf{R}_i(z) =0$. This is sketched in Appendix \ref{sec: R(z) est} below. 

Moving on to the second equivalence, 
let us write $\avar{i} = \iota_{n,n+1} \circ \eta (\mathbf{A}_{\bar{i},+}(z)) + \mathbf{S}(z)$ for some $\mathbf{S}(z) \in \fext{\Uq_{d_{\tau(i)} \geq 1,+}}{z}$. It is clear from \eqref{eq: Kolb emb} and Lemma \ref{lem: A-1 gen} that $[Q_{\tau(i)}, \avar{i}]_{q}$, $[A_{\tau(i),-1}, \mathbf{S}(z)]_{q}$ and $[A_{\tau(i),0}, \mathbf{S}(z)]_{q\mi}$ lie in $\fext{\Uq_+}{z}$. 
Using Lemma \ref{lem: A Th as series}, 
this implies that $\thvar{i} \equiv \iota_{n,n+1} \circ \eta (\boldsymbol{\grave{\Theta}}_{\bar{i},+}(z))$. Since the Drinfeld positive half of $\Uq(\widehat{\mathfrak{sl}}_3)$ 
is mapped by $\iota_{n,n+1}$ into the Drinfeld positive half of $\Uq$, the desired equivalence \linebreak $\thvar{i} \equiv K_i K_{\tau(i)}\mi \boldsymbol\phi^-_i(z\mi)\boldsymbol{\phi}^+_{\tau(i)}(\CCC z)$ 
now follows from Theorem \ref{thm: rank 1 factorization aij-1}. 
\enproof

\section{Coproduct formula}
\label{sec: coproduct}

In this short section, we apply the factorization formula (Theorem \ref{thm: main overall}) to show that the series $\thvar{i}$ are ``approximately group-like''.  

\Thm \label{thm: coproduct main}
For all $i \in \indx_0$, the following coproduct formula holds: 
\begin{align*}
\Delta(\thvar{i}) &\equiv \thvar{i} \otimes K_i K_{\tau(i)}\mi  \boldsymbol\phi^-_i(z\mi)\boldsymbol{\phi}^+_{\tau(i)}(\CCC z) \\
&\equiv \thvar{i} \otimes \thvar{i} 
\end{align*}
modulo $\fext{(\mathbf{U}^\imath_{\mathbf{c}} \otimes \Uq_+)}{z}$. 
\enthm

\Proof
The second equivalence clearly follows from the first in light of Theorem \ref{thm: main overall}. So let us prove the first equivalence. By Theorem \ref{thm: main overall} and the `group-like' property of the series $\boldsymbol\phi_i^{\pm}(z)$ (see \cite[Proposition 7.1]{damiani-r}), we get 
\begin{align*}
\Delta(\thvar{i}) &= \Delta \big(K_i K_{\tau(i)}\mi  \boldsymbol\phi^-_i(z\mi)\boldsymbol{\phi}^+_{\tau(i)}(\CCC z) + R \big) \\
&= (K_i K_{\tau(i)}\mi \otimes K_i K_{\tau(i)}\mi) \Delta(\boldsymbol\phi^-_i(z\mi)) \Delta(\boldsymbol{\phi}^+_{\tau(i)}(\CCC z)) + \Delta(R)  \\
&= \big(K_i K_{\tau(i)}\mi \otimes K_i K_{\tau(i)}\mi \big) \big(\boldsymbol\phi^-_i(z\mi) \otimes \boldsymbol\phi^-_i(z\mi) + R_1 \big) \big(\boldsymbol{\phi}^+_{\tau(i)}(\CCC z) \otimes \boldsymbol{\phi}^+_{\tau(i)}(\CCC z) + R_2 \big) \\
&\quad + \Delta(R)  \\ 
&= (\thvar{i} - R) \otimes \big( K_i K_{\tau(i)}\mi  \boldsymbol\phi^-_i(z\mi)\boldsymbol{\phi}^+_{\tau(i)}(\CCC z) \big) + \Delta(R) + Z, 
\end{align*}
for some $R \in \Uq_+$; $R_1, R_2 \in \Uq_- \otimes \Uq_+$, and   
\[
Z = \big(K_i K_{\tau(i)}\mi \otimes K_i K_{\tau(i)}\mi \big) \left[ \left( \boldsymbol\phi^-_i(z\mi) \otimes \boldsymbol\phi^-_i(z\mi) \right) R_2 + R_1 \left( \boldsymbol{\phi}^+_{\tau(i)}(\CCC z) \otimes \boldsymbol{\phi}^+_{\tau(i)}(\CCC z) \right) \right]. 
\]
Since $\mathbf{U}^\imath_{\mathbf{c}}$ is a coideal subalgebra, we must have $\Delta(\thvar{i}) \in \fext{\big(\mathbf{U}^\imath_{\mathbf{c}} \otimes \Uq \big)}{z}$. In particular, it follows that 
\begin{align} \label{eq: Z-term degrees}
-R \otimes \big( K_i K_{\tau(i)}\mi  \boldsymbol\phi^-_i(z\mi)\boldsymbol{\phi}^+_{\tau(i)}(\CCC z) \big) + \Delta(R) + Z  
\in \fext{\big(\mathbf{U}^\imath_{\mathbf{c}} \otimes \Uq \big)}{z}. 
\end{align}
Since $R \in \fext{\Uq_+}{z}$, we can write 
\[
\Delta(R) = a + b, \qquad a \in \fext{(\Uq \otimes \Uq_+)}{z}, \qquad b \in \bigcup_{j \in \indx_0} \fext{(\Uq_{d_j \geq 1} \otimes \Uq_{d_j \leq 0})}{z}. 
\] 
Decomposing \eqref{eq: Z-term degrees} as a sum of homogeneous terms with respect to the right tensor factor, we deduce that 
\begin{align*}
Z + a &\in \fext{( \mathbf{U}^\imath_{\mathbf{c}} \otimes \Uq_+ )}{z}, \\ 
b -R \otimes \big( K_i K_{\tau(i)}\mi  \boldsymbol\phi^-_i(z\mi)\boldsymbol{\phi}^+_{\tau(i)}(\CCC z) \big)   &\in \bigcup_{j \in \indx_0} \fext{( (\mathbf{U}^\imath_{\mathbf{c}} \cap \Uq_{d_j \geq 1}) \otimes \Uq_{d_j \leq 0})}{z}. 
\end{align*} 
Using a similar argument to that in the proof of \cite[Theorem 7.5]{Przez-23}, one can show that $\mathbf{U}^\imath_{\mathbf{c}} \cap \Uq_{d_j \geq 1} = \varnothing$, completing the proof. 
\enproof

\section{Boundary $q$-characters} 
\label{sec: bound q-char}
\nc{\Uh}{U_q(\widetilde{\mathfrak{h}})}
\nc{\Uic}{\mathbf{U}^\imath_{\mathbf{c}}}

In this section, we show that, as an application of Theorems \ref{thm: main overall} and \ref{thm: coproduct main}, boundary $q$-characters for quantum affine symmetric pairs are compatible, in an appropriate sense, with Frenkel and Reshetikhin's $q$-characters for quantum affine algebras. 

Fix $\mathbf{c} = (c_0, \cdots, c_N) \in (\C^\times)^{N+1}$ satisfying $c_i = c_{\tau(i)}$, and let $C$ be the image of $\mathfrak{C}$ in~$\Uic$. 
Let $\Rep \Uq$ be the monoidal category of finite dimensional representations of~$\Uq$, acting on $\Rep \Uic$, the category of finite dimensional representations of $\Uic$, via the coproduct \eqref{eq: coproduct explicit}.

\subsection{Reminder on $q$-characters for quantum affine algebras} 

The $q$-characters for quantum affine algebras were introduced by Frenkel and Reshetikhin \cite{FrenRes}.  
Below we briefly recall their two definitions and their equivalence. 

\subsubsection{$R$-matrix formulation} 

Let $\Uh$ be the subalgebra of $\Uq$ generated by $h_{i,k}$ ($i \in \indx_0$, $k < 0$). Following \cite{FrenRes, FrenMuk-comb}, set 
\begin{align} \label{eq: defi of Yia}
Y_{i,a} =& \  K_{\omega_i}^{-1} \exp \left( -(q-q\mi) \sum_{k > 0} \tilde{h}_{i,-k} a^k z^k  \right) \in \Uhz \qquad (a \in \C^\times),
\end{align}
where 
\eq
\tilde{h}_{i,-k} = \sum_{j \in \indx_0} \widetilde{C}_{ji}(q^k) \hg_{j,-k},
\eneq 
and $\widetilde{C}(q)$ is the inverse of the $q$-Cartan matrix. 
Let $\Yring = \Z[Y_{i,a}^{\pm 1}]_{i \in \indx_0, a \in \C^{\times}}$. 
By \cite[Theorem 3]{FrenRes}, there exists an injective ring homomorphism 
\[
\chmap \colon [\Rep \Uq] \to \Yring \subset \Uhz,
\]
called the \emph{$q$-character map}, given by 
\eq \label{eq: LSS dec}
[V] \mapsto \Tr_V \left[ \exp \left( -(q-q\mi) \sum_{i \in \indx_0} \sum_{k >0} \frac{k}{[k]} \pi_V(h_{i,k}) \otimes \tilde{h}_{i,-k} z^k \right) \cdot (\pi_V \otimes 1)(T) \right], 
\eneq
where $T$ is as in \cite[(3.8)]{FrenRes}, and $\pi_V \colon \Uq \to \End(V)$ is the representation. 
The expression \eqref{eq: LSS dec} is derived from the factorization of the universal $R$-matrix due to Damiani, 
Khoroshkin--Tolstoy and Levendorsky--Soibelman--Stukopkin \cite{damiani-r, khor-tol, lev-sob-st}.

\subsubsection{Algebraic formulation} 
\label{subsec: alg form}

By \cite[Proposition 2.4]{FrenMuk-comb}, the $q$-character map can equivalently be defined more explicitly in terms of the joint spectrum of the Drinfeld--Cartan operators series $\boldsymbol{\phi}^{\pm}_i(z)$. More precisely, there is a one-to-one correspondence between the monomials occurring in $\chmap(V)$ and the common eigenvalues of $\boldsymbol{\phi}^{\pm}_i(z)$ on a finite-dimensional representation $V$. 
Let 
$V = \bigoplus_{\gamma} V_{\gamma}$, with $\gamma = (\gamma^\pm_{i,\pm m})_{i \in \indx_0, m \in \Z_{\ge 0}}$, where 
\begin{align*}
V_{\gamma} = \{ v \in V \mid \exists p \ \forall i \in \indx_0 \ \forall m \in \Z_{\ge 0}:  (\phi^\pm_{i,\pm m} - \gamma^\pm_{i,\pm m})^p \cdot v = 0 \}, 
\end{align*}
be the generalized eigenspace decomposition of $V$. 
Collect the eigenvalues into generating series $\gamma^\pm_i(z) = \sum_{m \ge 0} \gamma^\pm_{i,\pm m} z^{\pm m}$. 
By \cite[Proposition 2.4]{FrenMuk-comb}, the series $\gamma^\pm_i(z)$ are expansions (at $0$ and $\infty$, respectively) of  the same rational function of the form
\eq \label{eq: FR eig}
\gamma^\pm_i(z) = q^{\deg Q_i - \deg R_i} \frac{Q_i(q^{-1}z)R_i(qz)}{Q_i(qz)R_i(q^{-1}z)}, 
\eneq
for some polynomials $Q_i(z), R_i(z)$ with constant term $1$. Writing  
\eq \label{eq: QR notation}
Q_i(z) = \prod_{r=1}^{k_i} (1 - z a_{i,r}), \qquad R_i(z) = \prod_{s=1}^{l_i} (1 - z b_{i,s}), 
\eneq
the $q$-character of $V$ can now be expressed as 
\eq \label{eq: FR qchar alg}
\chi_q(V) = \sum_{\gamma} \dim(V_\gamma) M_\gamma, \qquad 
M_\gamma = \prod_{i \in \indx_0} \prod_{r=1}^{k_i} Y_{i,a_{i,r}} \prod_{s=1}^{l_i} Y_{i,b_{i,s}}\mi. 
\eneq

\subsection{Boundary $q$-characters} 

In \cite{Przez-23, LP25, LP25b}, we proposed to imitate the algebraic approach recalled in \S \ref{subsec: alg form} to define $q$-characters for quantum affine symmetric pairs. In reference to the role quantum symmetric pairs play in the theory of integrable systems with boundary, we call them `\emph{boundary $q$-characters}'. Our approach does not require any information about the universal $K$-matrix. Although the theory of the universal $K$-matrix for quantum affine symmetric pairs has recently seen significant progress, especially in a series of works by Appel and Vlaar \cite{AppelVlaar, AppelVlaar4, AppelVlaar3, AppelVlaar2}, it appears that still more work is required to obtain an explicit $q$-character theory via this route. 

A key difficulty of the algebraic approach is that it is highly non-trivial to show that boundary $q$-characters are compatible with the usual $q$-characters for quantum affine algebras. The desired compatibility can, roughly speaking, be expressed as the boundary $q$-character map being a \emph{module homomorphism}. As explained below, this follows as an application of Theorems \ref{thm: main overall} and \ref{thm: coproduct main}.

\subsubsection{Eigenvalues on restricted representations} 

First, we need some notation. Given a polynomial $P(z) \in \C[z]$ with constant term $1$, let $P^\dag(z)$ be the polynomial with constant term $1$ whose roots are obtained from those of $P(z)$ via the transformation $a \mapsto \CC\mi a\mi$; and let $P^*(z)$ be the polynomial with constant term $1$ whose roots are the inverses of the roots of $P(z)$. Observe that
\eq \label{eq: P* formula} 
P(z) = \xi_P z^{\deg P} P^\dag(\CC\mi z\mi), \qquad 
P(z) = \xi_{P} z^{\deg P} P^*(z^{-1}),
\eneq
where $\xi_P$ is the product of the roots of $P^*(z)$. 

\Thm
\label{cor: FR thm Oq} 
Let $W \in \Rep \Uq$. 
Then the generalized eigenvalues of $\thvar{i}$ on the restricted representation $W|_{\Uic}$ are of the form: 
\eq \label{eq: FR formula Oq}
\boldsymbol\gamma_i(z) = q^{\deg \mathbf{Q}_i - \deg \mathbf{Q}_{\tau(i)}}  \kappa_i \kappa_{\tau(i)}\mi \frac{\ourQ_i(q^{-1}z)}{\ourQ_i(qz)} \frac{\ourQ_{\tau(i)}^\dag(q z)}{\ourQ_{\tau(i)}^\dag(q^{-1}z)}, 
\eneq
where $\ourQ_i(z)$ is a polynomial with constant term $1$, and $\kappa_i = q^{(\alpha_i, \lambda)}$, for $\lambda$ the weight of the associated eigenvector. 
Explicitly, 
\[
\ourQ_i(z) = {Q_{\tau(i)}(C z)R_i^*(z)}, \quad  \ourQ_i^\dag(z) = Q_{\tau(i)}^*(z)R_i(C z). 
\]
\enthm

\Proof 
Theorem \ref{thm: main overall} implies that the action of $\thvar{i}$ on $W$ is the sum of the action of $K_i K_{\tau(i)}\mi   \boldsymbol\phi^-_i(z\mi)\boldsymbol{\phi}^+_{\tau(i)}(C z)$ and a nilpotent operator. Hence the eigenvalues of $\thvar{i}$ coincide with those of $K_i K_{\tau(i)}\mi   \boldsymbol\phi^-_i(z\mi)\boldsymbol{\phi}^+_{\tau(i)}(C z)$. Therefore, \eqref{eq: FR eig} and \eqref{eq: P* formula} imply that 
\begin{align*}
\boldsymbol\gamma_i(z) &= \kappa_i \kappa_{\tau(i)}\mi \gamma_i^-(z\mi) \gamma_{\tau(i)}^+(\CC  z) \\ 
&= q^{s_i+s_{\tau(i)}-r_i-r_{\tau(i)}}  \kappa_i \kappa_{\tau(i)}\mi \frac{Q_i(q^{-1}{z\mi})R_i(q{z\mi})}{Q_i(q{z\mi})R_i(q^{-1}{z\mi})} \frac{Q_{\tau(i)}(q^{-1}Cz)R_{\tau(i)}(qCz)}{Q_{\tau(i)}(qCz)R_{\tau(i)}(q^{-1}Cz)} \\ 
&= q^{-s_i+s_{\tau(i)}+r_i-r_{\tau(i)}}  \kappa_i \kappa_{\tau(i)}\mi \frac{Q_i^*(qz)R_i^*(q\mi{z})}{Q^*_i(q\mi{z})R^*_i(q{z})} \frac{Q_{\tau(i)}(q^{-1}Cz)R_{\tau(i)}(qCz)}{Q_{\tau(i)}(qCz)R_{\tau(i)}(q^{-1}Cz)} \\ 
&= q^{-s_i+s_{\tau(i)}+r_i-r_{\tau(i)}}  \kappa_i \kappa_{\tau(i)}\mi \left( \frac{Q_{\tau(i)}( q^{-1}Cz)R^*_i(q\mi z)}{Q_{\tau(i)}(q Cz)R_i^*(q z)} \right) \left( \frac{Q_i^*(q z)R_{\tau(i)}(q Cz)}{Q_i^*(q\mi z) R_{\tau(i)}( q^{-1}Cz)}\right) \\
&= q^{\deg \mathbf{Q}_i - \deg \mathbf{Q}_{\tau(i)}}  \kappa_i \kappa_{\tau(i)}\mi \frac{\ourQ_i(q^{-1}z)}{\ourQ_i(qz)} \frac{\ourQ_{\tau(i)}^\dag(q z)}{\ourQ_{\tau(i)}^\dag(q^{-1}z)}, 
\end{align*} 
where $s_i = \deg Q_i$ and $r_i = \deg R_i$. 
\enproof 

\subsubsection{Definition}

Recall that $\rho = q-q\mi$. 
Define the following series 
\[
\mathbf{H}_i(z) = \sum_{r \geq 1} H_{i,r} z^r, \qquad \mathbf{\grave{H}}_i(z) = \sum_{r \geq 1} \grave{H}_{i,r} z^r = \rho\mi \on{log} \left( \thvar{i} \right). 
\] 
In particular, if $a_{i,\tau(i)} = 0$, then $\mathbf{\grave{H}}_i(z) = \mathbf{H}_i(z)$. 
Let\footnote{The definition of $\mathcal{K}^0$ in \cite[Definition 9.1]{LP25} was not correct. One should replace it with the version presented here. Note that, in the setting of that paper, $i = \tau(i)$ holds always, so $TT_\tau\mi = 1$.}
\eq \label{eq: K0}
\mathcal{K}^0 = \exp \left( -\rho \sum_{i \in \indx_0} \sum_{k >0} \frac{k}{[k]} \grave{H}_{i,k} \otimes \tilde{h}_{i,-k} z^k \right) \cdot TT_\tau\mi \in \fext{\Uic \otimes \Uh}{z}, 
\eneq
where 
\[
T = q^{\sum_{i,j \in\indx_0} d_{ij} h_{\alpha_i} \otimes h_{\alpha_j}}, \qquad T_\tau = q^{\sum_{i,j\in\indx_0} d_{ij} h_{\alpha_i} \otimes h_{\alpha_{\tau(j)}}}, 
\]
with coroots $h_{\alpha_i}$ and inverse Cartan matrix $(d_{ij})$. 
We define the \emph{boundary $q$-character map} to be 
\[
\ichmap \colon [\Rep \Uic] \ \to \ \Uhz, \qquad [V] \mapsto \Tr_V((\pi_{V(z)} \otimes 1)(\mathcal{K}^0)). 
\]

\subsubsection{Compatibility} 

Consider $\Uhz$ as a $\Yring$-module via the ring homomorphism
\eq \label{eq: Yring module}
\Yring \to \Yring \hookrightarrow \Uhz, \qquad Y_{i,a} \mapsto \mathbf{Y}_{i,a} := Y_{\tau(i),Ca}Y_{i,a\mi}\mi, 
\eneq
i.e., $Y_{i,a}^{\pm1}$ acts on $\Uhz$ via multiplication by  $\mathbf{Y}_{i,a}^{\pm1}$.

\Thm \label{cor: comm diagram qchar actions}
Let $(\Uq, \Uic)$ be a quasi-split quantum affine symmetric pair of type~$\mathsf{AIII}_N^{(\tau)}$. Then the following diagram commutes: 
\[
\begin{tikzcd}[ row sep = 0.2cm]
{[\Rep \Uq]}  \arrow[r, "\chmap"] & \Yring  \\
 \curvearrowright  & \curvearrowright  \\
{[\Rep \Uic]} \arrow[r, "\ichmap"] & \Uhz. 
\end{tikzcd}
\]
\enthm

\Proof 
The proof is a modification of the proof of \cite[Proposition 2.4]{FrenMuk-comb}. 
Let $W \in \Rep \Uq$ and $V \in \Rep \Uic$. We need to calculate $\ichmap(V \otimes W)$. Theorem \ref{thm: coproduct main} implies that the eigenvalues of $\thvar{i}$ on $V \otimes W$ are a product of the eigenvalues of $\thvar{i}$ on $V$ and $W$, respectively. Hence $\Tr_{V\otimes W}((\pi_{V \otimes W(z)} \otimes 1)(\mathcal{K}^0)) = \Tr_V((\pi_{V(z)} \otimes 1)(\mathcal{K}^0)) \cdot \Tr_W((\pi_{W(z)} \otimes 1)(\mathcal{K}^0))$. Moreover, we \emph{claim} that 
the eigenvalues of $\grave{H}_{i,m}$ on $W$ are of the form 
\eq \label{eq: eig with twist}
\frac{[m]}{m}\left( \sum_{r=1}^{k_{\tau(i)}} (\CC a_{\tau(i),r})^m - \sum_{r=1}^{k_{i}} a_{i,r}^{-m} -
\sum_{s=1}^{l_{\tau(i)}} (\CC b_{\tau(i),s})^m + \sum_{s=1}^{l_i} b_{i,s}^{-m} \right),
\eneq
keeping the notation from \eqref{eq: QR notation}. Plugging \eqref{eq: eig with twist} into the definition of $\ichmap(W)$, and comparing with \eqref{eq: defi of Yia}, we conclude that $\ichmap(W)$ is a linear combination of monomials of the form 
\begin{align*}
\prod_{i \in \indx_0} \prod_{r=1}^{k_{\tau(i)}} Y_{i,C a_{\tau(i),r}}  \prod_{r=1}^{k_i} &Y_{i,a_{i,r}\mi}\mi \prod_{s=1}^{l_i} Y_{i, b_{i,s}\mi}  \prod_{s=1}^{l_{\tau(i)}} Y_{i,C b_{\tau(i),s}}\mi = \\ 
&= \prod_{i \in \indx_0}  \prod_{r=1}^{k_{\tau(i)}} Y_{i,C a_{\tau(i),r}} Y_{\tau(i) ,a_{\tau(i),r}\mi}\mi \prod_{s=1}^{l_{\tau(i)}} Y_{i, Cb_{\tau(i),s}}\mi Y_{\tau(i), b_{\tau(i),s}\mi} \\
&=  \prod_{i \in \indx_0}  \prod_{r=1}^{k_{\tau(i)}} \mathbf{Y}_{\tau(i),a_{\tau(i),r}} \prod_{s=1}^{l_{\tau(i)}} \mathbf{Y}_{\tau(i),b_{\tau(i),s}}\mi = \prod_{i \in \indx_0}  \prod_{r=1}^{k_{i}} \mathbf{Y}_{i,a_{i,r}} \prod_{s=1}^{l_{i}} \mathbf{Y}_{i,b_{i,s}}\mi. 
\end{align*} 
Comparing with \eqref{eq: FR qchar alg}, it follows that $\ichmap(W)$ is equal to the image of $\chmap(W)$ under \eqref{eq: Yring module}. 

To complete the proof, let us prove the \emph{claim}. By Theorem \ref{cor: FR thm Oq}, the eigenvalues of $\thvar{i}$ on $W$ are equal to $\boldsymbol\gamma_i(z)$. Let $\boldsymbol\xi_i(z)$ denote the series of eigenvalues of $\grave{\mathbf{H}}_i(z)$. Then 
\begin{align*}
\boldsymbol\xi_i(z) &= \rho\mi \on{log} \left( \boldsymbol\gamma_i(z) \right). 
\end{align*} 
Using the formula $\on{log}(1-x) = - \sum_{k \geq 1} \frac{x^k}{k}$, the $m$-th coefficient of $\boldsymbol\xi_i(z)$ is easily seen to equal \eqref{eq: eig with twist}, for $m \geq 1$. 
\enproof

\appendix
\section{Auxiliary calculations} 

\subsection{The case: ${\sf AIII}_{2n-1}^{(\tau)}$} 

\Lem \label{lem: FE vanishing}
If $a_{ji} = -1$, then $[F_j,\tE_i]_{q} = 0$. 
\enlem 

\Proof
The equality
\eq 
F_jE_{\tau(i)}\Kb_i - q E_{\tau(i)}\Kb_i F_j = q E_{\tau(i)}\Kb_i F_j- q E_{\tau(i)}\Kb_i F_j = 0 
\eneq 
yields the statement. 
\enproof 

\Lem \label{lem: TjTiBj formula}
Let $a_{ij} = -1$. Then 
\[
\Tbr_j \Tbr_i (B_j) = 
\left\{
\begin{array}{r l}
\mathbf{P}(B_i, B_{\tau(i)}, B_j) & \quad \mbox{if } \ \tau(j) = j \ \& \ \tau(i) \neq i, \\ \\
{[B_i, B_{\tau(j)}]}_{q} & \quad \mbox{if }\ \tau(j) \neq j \ \& \ \tau(i) = i.
\end{array}
\right. 
\]
\enlem 

\begin{proof}
We begin by checking the first case. Assume that $\tau(j)=j$, $\tau(i) \ne i$. We have 
\begin{align*}
\Tbr_i(B_j) = B_j B_i B_{\tau(i)}-q B_i B_j B_{\tau(i)}-q B_{\tau(i)} B_j B_i+q^2 B_{\tau(i)} B_i B_j-q B_j \KK_i,
\end{align*}
and $\Tbr_j \Tbr_i(B_j)$ is
\begin{align*}
& \KK^{-1}_j B_j B_i B_j B_{\tau(i)} B_j-q \KK^{-1}_j B_j B_i B_j B_j B_{\tau(i)}-q \KK^{-1}_j B_j B_j B_i B_{\tau(i)} B_j \\
& + q^2 \KK^{-1}_j B_j B_j B_i B_j B_{\tau(i)} -q B_i B_j \KK^{-1}_j B_j B_{\tau(i)} B_j + q^2 B_i B_j \KK^{-1}_j B_j B_j B_{\tau(i)} \\
& + q^2 B_j B_i \KK^{-1}_j B_j B_{\tau(i)} B_j -q^3 B_j B_i \KK^{-1}_j B_j B_j B_{\tau(i)} -q B_{\tau(i)} B_j \KK^{-1}_j B_j B_i B_j \\
& + q^2 B_{\tau(i)} B_j \KK^{-1}_j B_j B_j B_i + q^2 B_j B_{\tau(i)} \KK^{-1}_j B_j B_i B_j -q^3 B_j B_{\tau(i)} \KK^{-1}_j B_j B_j B_i \\
& + q^2 B_{\tau(i)} B_j B_i B_j \KK^{-1}_j B_j -q^3 B_{\tau(i)} B_j B_j B_i \KK^{-1}_j B_j -q^3 B_j B_{\tau(i)} B_i B_j \KK^{-1}_j B_j \\
& + q^4 B_j B_{\tau(i)} B_j B_i \KK^{-1}_j B_j -q \KK^{-1}_j B_j \KK_i \KK_j. 
\end{align*}
Applying the relations 
\begin{align*}
& B_{\tau(i)}B_{j}B_j =  [2] B_j B_{\tau(i)} B_{j} -  B_j B_j B_{\tau(i)} - q^{-1} B_{\tau(i)} \KK_j, \\
& B_{j}B_j B_{i} =  [2] B_j B_{i} B_{j} -  B_i B_j B_{j} - q^{-1} B_{i} \KK_j, \\
& B_{\tau(i)}B_i = B_i B_{\tau(i)} + \frac{1}{q-q^{-1}} \KK_i - \frac{1}{q-q^{-1}} \KK_{\tau(i)}, \\
& B_j \KK_i = q^{-a_{ij} + a_{\tau(i),j}} \KK_i B_j,
\end{align*}
we obtain 
\begin{align*}
& B_i B_{\tau(i)} B_j-q B_i B_j B_{\tau(i)} - \frac{1}{q-q^{-1}} (\KK_{\tau(i)}q^2 - \KK_i) B_j-q B_{\tau(i)} B_j B_i + q^2 B_j B_i B_{\tau(i)} \\
& = B_i [ B_{\tau(i)}, B_j ]_q - q \KK_i B_j - q [ B_{\tau(i)}, B_j ]_q B_i = [ B_i,  [ B_{\tau(i)}, B_j ]_q ]_q - q \KK_i B_{j}.
\end{align*} 

Now we check the second equation. Assume that $\tau(j)\ne j$, $\tau(i)=i$. Then 
\begin{align*}
\Tbr_i(B_j) = B_j B_i - q B_i B_j,
\end{align*}
and $\Tbr_j \Tbr_i(B_j)$ is equal to
\begin{align*}
&  ( - \KK_j^{-1} B_{\tau(j)} ) (  B_i B_j B_{\tau(j)} -q B_j B_i B_{\tau(j)} -q B_{\tau(j)} B_i B_j + q^2 B_{\tau(j)} B_j B_i-q B_i \KK_j )  \\
& - q (   B_i B_j B_{\tau(j)} -q B_j B_i B_{\tau(j)} -q B_{\tau(j)} B_i B_j + q^2 B_{\tau(j)} B_j B_i-q B_i \KK_j  ) (  - \KK_j^{-1} B_{\tau(j)}  ).
\end{align*}
By using  relations $ B_j \KK_i = q^{-a_{ij} + a_{\tau(i),j}} \KK_i B_j$, this is equal to
\begin{align*}
& - q[2] \KK^{-1}_j B_{\tau(j)} B_i B_j B_{\tau(j)} + [2] \KK^{-1}_j B_{\tau(j)} B_j B_i B_{\tau(j)} + q\KK^{-1}_j B_{\tau(j)} B_{\tau(j)} B_i B_j + q^{-1} B_{\tau(j)} B_i  \\
& -q^2\KK^{-1}_j B_{\tau(j)} B_{\tau(j)} B_j B_i + q\KK^{-1}_j B_i B_j B_{\tau(j)} B_{\tau(j)}-q^2\KK^{-1}_j B_j B_i B_{\tau(j)} B_{\tau(j)}-q^2 B_i B_{\tau(j)}.
\end{align*}
Applying the relations
\begin{align*}
& B_{\tau(j)}B_{\tau(j)}B_i = - B_i B_{\tau(j)} B_{\tau(j)} + [2] B_{\tau(j)} B_i B_{\tau(j)}, \\
& B_{\tau(j)}B_j = B_j B_{\tau(j)} + \frac{1}{q-q^{-1}} \KK_j - \frac{1}{q-q^{-1}} \KK_{\tau(j)}, \\
&  B_j \KK_i = q^{-a_{ij} + a_{\tau(i),j}} \KK_i B_j,
\end{align*}
 we obtain
\begin{align*}
q^2 [2] \KK_j^{-1} B_j B_{\tau(j)} B_i B_{\tau(j)} + [B_i, B_{\tau(j)}]_q - q^2 \KK_j^{-1} B_j B_{\tau(j)} B_{\tau(j)} B_i - q^2 \KK_j^{-1} B_j B_{i} B_{\tau(j)} B_{\tau(j)}.
\end{align*}
Applying the Serre relation again, we finally obtain $\Tbr_j \Tbr_i (B_j) = [B_i, B_{\tau(j)}]_q$.  
\end{proof}

\Lem \label{lem: PBBE}
We have 
\[
\mathbf{P}(B_{2n-1}, B_1, \tE_0) =  [E_{2n-1},[E_1, E_0]_q]_q  \Kb_0\Kb_{1}\Kb_{2n-1}. 
\]
\enlem 

\Proof 
It suffices to prove the result in the case of $n=3$. By definition, ${\bf P}(B_5, B_1, E_0 \tilde{K}'_0)$ is equal to
\begin{align*}
& F_5 F_1 E_0 \widetilde{K}'_0-q F_1 E_0 \widetilde{K}'_0 F_5-q F_5 E_0 \widetilde{K}'_0 F_1+q^2 E_0 \widetilde{K}'_0 F_1 F_5+ F_5 E_5 \widetilde{K}'_1 E_0 \widetilde{K}'_0
-q E_5 \widetilde{K}'_1 E_0 \widetilde{K}'_0 F_5 \\
& -q F_5 E_0 \widetilde{K}'_0 E_5 \widetilde{K}'_1+q^2 E_0 \widetilde{K}'_0 E_5 \widetilde{K}'_1 F_5+E_1 \widetilde{K}'_5 F_1 E_0 \widetilde{K}'_0
-q F_1 E_0 \widetilde{K}'_0 E_1 \widetilde{K}'_5-q E_1 \widetilde{K}'_5 E_0 \widetilde{K}'_0 F_1 \\
& +q^2 E_0 \widetilde{K}'_0 F_1 E_1 \widetilde{K}'_5+E_1 \widetilde{K}'_5 E_5 \widetilde{K}'_1 E_0 \widetilde{K}'_0
-q E_5 \widetilde{K}'_1 E_0 \widetilde{K}'_0 E_1 \widetilde{K}'_5-q E_1 \widetilde{K}'_5 E_0 \widetilde{K}'_0 E_5 \widetilde{K}'_1 \\
& +q^2 E_0 \widetilde{K}'_0 E_5 \widetilde{K}'_1 E_1 \widetilde{K}'_5-q K_5 \widetilde{K}'_1 E_0 \widetilde{K}'_0.
\end{align*}
Using the relations
\begin{equation}  \label{eq:relation-set-first-step-bPB5B1E0tilde}
\Ka_i\Eg_j^{\pm} = q^{\pm a_{ij}} \Eg_j^\pm\Ka_i, \quad 
\Kb_i\Eg_j^{\pm} = q^{\mp a_{ij}} \Eg_j^\pm\Kb_i, \quad [\Ka_i,\Ka_j] = [\Kb_j, \Kb_i] = [\Ka_i, \Kb_j] = 0, 
\end{equation}
the above is equal to 
\begin{align*}
& F_5 F_1 E_0 \widetilde{K}'_0- F_1 E_0 F_5 \widetilde{K}'_0- F_5 E_0 F_1 \widetilde{K}'_0+E_0 F_1 F_5 \widetilde{K}'_0 +q F_5 E_5 E_0 \widetilde{K}'_1 \widetilde{K}'_0-q E_5 E_0 F_5 \widetilde{K}'_1 \widetilde{K}'_0 \\
&+q^2 E_0 E_5 F_5 \widetilde{K}'_1 \widetilde{K}'_0+q E_1 F_1 E_0 \widetilde{K}'_5 \widetilde{K}'_0-q^2 F_1 E_0 E_1 \widetilde{K}'_5 \widetilde{K}'_0 -q E_1 E_0 F_1 \widetilde{K}'_5 \widetilde{K}'_0 \\
& +q^2 E_0 F_1 E_1 \widetilde{K}'_5 \widetilde{K}'_0+E_1 E_5 E_0 \widetilde{K}'_5 \widetilde{K}'_1 \widetilde{K}'_0-q E_5 E_0 E_1 \widetilde{K}'_5 \widetilde{K}'_1 \widetilde{K}'_0 -q E_1 E_0 E_5 \widetilde{K}'_5 \widetilde{K}'_1 \widetilde{K}'_0 \\
& +q^2 E_0 E_5 E_1 \widetilde{K}'_5 \widetilde{K}'_1\widetilde{K}'_0 -q E_0 K_5 \widetilde{K}'_1 \widetilde{K}'_0-q^2 F_5 E_0 E_5 \widetilde{K}'_1.  \widetilde{K}'_0.
\end{align*}
Using relations $[\Eg_i^{+}, \Eg_j^{-}] = \delta_{ij} \frac{\Ka_i - {\Kb_i}}{q - q^{-1}}$ and \eqref{eq:relation-set-first-step-bPB5B1E0tilde}, one checks that 
the above is equal to
\begin{align*}
(E_5 E_1 E_0 -q E_5 E_0 E_1 -q E_1 E_0 E_5 +q^2 E_0 E_5 E_1 ) \widetilde{K}'_5 \widetilde{K}'_1 \widetilde{K}'_0,
\end{align*}
completing the proof. 
\enproof

\Lem \label{lem: P4 van} 
We have 
\[
P(F_{n-1}, \tE_{n}, \tE_{n+1}, \tE_{n+2}) = 0. 
\]
\enlem 

\Proof 

It suffices to consider the case of $n=3$. We have 
\begin{align*}
P(& F_{n-1}, \widetilde{E}_n, \widetilde{E}_{n+1}, \widetilde{E}_{n+2} )  =[F_2,  [ E_3 \widetilde{K}'_3,  [ E_2 \widetilde{K}'_4,  E_1 \widetilde{K}'_5]_q ]_q ]_q = \\
& =  F_2  E_3 \widetilde{K}'_3 E_2 \widetilde{K}'_4 E_1 \widetilde{K}'_5 - q E_3 \widetilde{K}'_3 E_2 \widetilde{K}'_4 E_1 \widetilde{K}'_5 F_2 
-q F_2 E_2 \widetilde{K}'_4 E_1 \widetilde{K}'_5 E_3 \widetilde{K}'_3 \\
&\quad  + q^2 E_2 \widetilde{K}'_4 E_1 \widetilde{K}'_5 E_3 \widetilde{K}'_3 F_2 -q F_2 E_3 \widetilde{K}'_3 E_1 \widetilde{K}'_5 E_2 \widetilde{K}'_4 + q^2 E_3 \widetilde{K}'_3 E_1 \widetilde{K}'_5 E_2 \widetilde{K}'_4 F_2 \\
&\quad  + q^2 F_2 E_1 \widetilde{K}'_5 E_2 \widetilde{K}'_4 E_3 \widetilde{K}'_3 -q^3 E_1 \widetilde{K}'_5 E_2 \widetilde{K}'_4 E_3 \widetilde{K}'_3 F_2 \\
 & = \big(q F_2 E_3 E_2 E_1 -q E_3 E_2 F_2 E_1-q^2 F_2 E_2 E_1 E_3 + q^2 E_2 F_2 E_1 E_3 \\
&\quad -q^2 F_2 E_3 E_1 E_2 + q^2 E_3 E_1 E_2 F_2 + q^3 F_2 E_1 E_2 E_3 -q^3 E_1 E_2 F_2 E_3 \big) \widetilde{K}'_5 \widetilde{K}'_4 \widetilde{K}'_3.
\end{align*}
Since $F_2 E_1 = E_1 F_2$ and $F_2 E_3 = E_3 F_2$, we get 
\begin{align*}
& q F_2 E_3 E_2 E_1 -q E_3 E_2 F_2 E_1-q^2 F_2 E_2 E_1 E_3 + q^2 E_2 F_2 E_1 E_3 \\
& -q^2 F_2 E_3 E_1 E_2 + q^2 E_3 E_1 E_2 F_2 + q^3 F_2 E_1 E_2 E_3 -q^3 E_1 E_2 F_2 E_3 = \\
&\qquad  = q  E_3 [F_2, E_2] E_1 -q^2 [F_2, E_2] E_1 E_3   -q^2  E_3 E_1 [F_2, E_2] + q^3  E_1 [F_2, E_2] E_3.
\end{align*}
The above is equal to $\frac{q^2}{q^2-1}$ times 
\begin{align*}
& E_3 ( \widetilde{K}'_2 - \widetilde{K}_2 )E_1  + q ( \widetilde{K}_2 - \widetilde{K}'_2) E_1 E_3  + q E_3 E_1 (\widetilde{K}_2 -  \widetilde{K}'_2) - q^2 E_1 ( \widetilde{K}_2 -  \widetilde{K}'_2) E_3, 
\end{align*}
which is $0$. 
\enproof

\subsection{The case: ${\sf AIII}_{2n}^{(\tau)}$} 

\Lem \label{lem: TjTiBj formula 2}
In type ${\sf AIII}_{2n}^{(\tau)}$, let $a_{ij} = -1$, $a_{j, \tau(j)} = 0$, $a_{i,\tau(i)} = -1$. Then 
\[
\Tbr_j \Tbr_i (B_j) = 
\mathbf{P}(B_i, B_{\tau(i)}, B_{\tau(j)}).
\]
\enlem 

\Proof
Let $a_{i,j}=-1$, $a_{j,\tau(j)}=0$, and $a_{i,\tau(i)}=-1$. Then, by the definition of the braid group action, 
\begin{align*}
{\bf T}_i(B_j) = B_j B_i B_{\tau(i)} -q B_i B_j B_{\tau(i)} -q B_{\tau(i)} B_j B_i + q^2 B_{\tau(i)} B_i B_j -  \mathbb{K}_i B_j,
\end{align*}
and 
\begin{align*}
& {\bf T}_j {\bf T}_i ( B_j) = - \mathbb{K}_j^{-1} B_{\tau(j)} B_i B_j B_{\tau(i)} B_{\tau(j)} + 
 q \mathbb{K}_j^{-1} B_{\tau(j)} B_i B_j B_{\tau(j)} B_{\tau(i)} \\
 & 
 -q^2 \mathbb{K}_j^{-1} B_{\tau(j)} B_j B_i B_{\tau(j)} B_{\tau(i)}+
 q B_i B_j \mathbb{K}_j^{-1} B_{\tau(j)} B_{\tau(i)} B_{\tau(j)}
 -q^2 B_i B_j \mathbb{K}_j^{-1} B_{\tau(j)} B_{\tau(j)} B_{\tau(i)}
\\
& -q^2 B_j B_i \mathbb{K}_j^{-1} B_{\tau(j)} B_{\tau(i)} B_{\tau(j)}+
 q^3 B_j B_i \mathbb{K}_j^{-1} B_{\tau(j)} B_{\tau(j)} B_{\tau(i)}+
 q B_{\tau(i)} B_{\tau(j)} \mathbb{K}_j^{-1} B_{\tau(j)} B_i B_j
\\
& -q^2 B_{\tau(i)} B_{\tau(j)} \mathbb{K}_j^{-1} B_{\tau(j)} B_j B_i
 -q^2 B_{\tau(j)} B_{\tau(i)} \mathbb{K}_j^{-1} B_{\tau(j)} B_i B_j+
 q^3 B_{\tau(j)} B_{\tau(i)} \mathbb{K}_j^{-1} B_{\tau(j)} B_j B_i
\\
& -q^2 B_{\tau(i)} B_{\tau(j)} B_i B_j \mathbb{K}_j^{-1} B_{\tau(j)}+
 q^3 B_{\tau(i)} B_{\tau(j)} B_j B_i \mathbb{K}_j^{-1} B_{\tau(j)}+
 q^3 B_{\tau(j)} B_{\tau(i)} B_i B_j \mathbb{K}_j^{-1} B_{\tau(j)}
\\
& -q^4 B_{\tau(j)} B_{\tau(i)} B_j B_i \mathbb{K}_j^{-1} B_{\tau(j)}
 -q \mathbb{K}_i \mathbb{K}_j \mathbb{K}_j^{-1} B_{\tau(j)} +
 q \mathbb{K}_j^{-1} B_{\tau(j)} B_j B_i B_{\tau(i)} B_{\tau(j)}. 
\end{align*}
We can reorder the above monomials as follows: 
\begin{align*}
& -\mathbb{K}_j^{-1} B_{\tau(j)} B_i B_j B_{\tau(i)} B_{\tau(j)}+
 q\mathbb{K}_j^{-1} B_{\tau(j)} B_i B_j B_{\tau(j)} B_{\tau(i)}+
 q\mathbb{K}_j^{-1} B_{\tau(j)} B_j B_i B_{\tau(i)} B_{\tau(j)}
\\
&  -q^2\mathbb{K}_j^{-1} B_{\tau(j)} B_j B_i B_{\tau(j)} B_{\tau(i)}+
 \mathbb{K}_j^{-1} B_i B_j B_{\tau(j)} B_{\tau(i)} B_{\tau(j)}
 -q\mathbb{K}_j^{-1} B_i B_j B_{\tau(j)} B_{\tau(j)} B_{\tau(i)}
\\
& -q\mathbb{K}_j^{-1} B_j B_i B_{\tau(j)} B_{\tau(i)} B_{\tau(j)}+
 q^2\mathbb{K}_j^{-1} B_j B_i B_{\tau(j)} B_{\tau(j)} B_{\tau(i)}+
 q^2\mathbb{K}_j^{-1} B_{\tau(i)} B_{\tau(j)} B_{\tau(j)} B_i B_j
\\
& -q^3\mathbb{K}_j^{-1} B_{\tau(i)} B_{\tau(j)} B_{\tau(j)} B_j B_i
 -q^3\mathbb{K}_j^{-1} B_{\tau(j)} B_{\tau(i)} B_{\tau(j)} B_i B_j+
 q^4 \mathbb{K}_j^{-1} B_{\tau(j)} B_{\tau(i)} B_{\tau(j)} B_j B_i
\\
& -q^2\mathbb{K}_j^{-1} B_{\tau(i)} B_{\tau(j)} B_i B_j B_{\tau(j)}+
 q^3\mathbb{K}_j^{-1} B_{\tau(i)} B_{\tau(j)} B_j B_i B_{\tau(j)}+
 q^3\mathbb{K}_j^{-1} B_{\tau(j)} B_{\tau(i)} B_i B_j B_{\tau(j)} \\
 &
 -q^4 \mathbb{K}_j^{-1} B_{\tau(j)} B_{\tau(i)} B_j B_i B_{\tau(j)}
 - q \mathbb{K}_i B_{\tau(j)}.
\end{align*}
Repeatedly applying the defining relations 
\begin{align*}
B_{\tau(j)} B_{\tau(j)} B_{\tau(i)} &= - B_{\tau(i)} B_{\tau(j)} B_{\tau(j)} + [2] B_{\tau(j)} B_{\tau(i)} B_{\tau(j)}, \\
B_{\tau(j)} B_j &= B_j B_{\tau(j)} + (q - q^{-1})^{-1} (\mathbb{K}_j - \mathbb{K}_{\tau(j)}),
\end{align*} 
we obtain the desired result. 
\enproof

\section{Estimation of $\mathbf{R}_i(z)$}
\label{sec: R(z) est} 

 Below we outline the necessary steps needed to show that $\mathbf{R}_i(z) \in \fext{\Uu_{d_{\tau(i)} \geq 1,+}}{z}$ (see \eqref{eq: Rz a0}, \eqref{eq: Rz a2}, and \eqref{eq: Rz a-1} for the definition). Most importantly, the estimation of $\mathbf{R}_i(z)$ requires finer control over the remainder term $Q_{\tau(i)}$ than we needed in \S \ref{sec: compat of br} and Lemma~\ref{lem: A-1 gen}. 
 
\Rem
In the split case \cite[\S 9]{Przez-23}, $\mathbf{R}_i(z)$ always vanishes. This is generally \emph{not} the case here.  
\enrem

\subsection{Case: $a_{i,\tau(i)} = 0$}

We need the following lemma. 

\Lem
If $a_{i,\tau(i)} = 0$ then 
\eq \label{eq: sum to vanish}
Q_{\tau(i)} = -q^2\CCC\mi \KK_{\tau(i)}  \sum_{\substack{\varepsilon_k \leq \varepsilon_{k+1},\\ k=1,...,i-1}} \sum_{\substack{\varepsilon_l \geq \varepsilon_{l+1},\\ l=i+1,...,N-1}} P(\tilde{e}_{i-1}^{\varepsilon_{i-1}}, \cdots, \tilde{e}_{1}^{\varepsilon_{1}}, P(\tilde{e}_{i+1}^{\varepsilon_{i+1}}, \cdots, \tilde{e}_{N}^{\varepsilon_{N}}, F_{0})), 
\eneq
where $(\varepsilon_{i-1}, \varepsilon_{i+1}) \in \{(+,-), (-,+), (+,+)\}$. 
\enlem 

\Proof
The proof is a straightforward generalization of the proof of Proposition \ref{pro: Aodd igood}. 
\enproof

Observe that commutation with $x_{i,r}^-$ is only affected by the superscripts $\varepsilon_{i-1}, \varepsilon_{i+1}, \varepsilon_{\theta(i)}$ in \eqref{eq: sum to vanish}. Let 
$P_{\varepsilon_{i-1}, \varepsilon_{i+1}, \varepsilon_{\theta(i)}}$ be the subsum of all the terms in \eqref{eq: sum to vanish} (including the factor $-q^2\CCC\mi \KK_{\tau(i)}$) with specified $\varepsilon_{i-1}, \varepsilon_{i+1}, \varepsilon_{\theta(i)}$ and other superscripts arbitrary. Let $P_{*,*,+} = P_{+,-,+} + P_{-,+,+} + P_{-,-,+}$. Then $P_{*,*,+} \in \Uu_{d_i = 2, d_{\tau(i)} = 1, +}$, and so 
\[
\CCC [  \KK_{\tau(i)}\mi [ x_{i,0}^- , \ P_{*,*,+}  ], \avarm{i}] \in \fext{\Uu_{d_{\tau(i)} \geq 1,+}}{z}. 
\] 
It is also easy to check that $[F_i, P_{+,+,-}] = 0$. Therefore, it remains to be shown that 
\[
\CCC [  \KK_{\tau(i)}\mi [ x_{i,0}^- , \ P_{-,+,-} + P_{+,-,-}  ], \avarm{i}] = 0. 
\] 
To do that, one shows that $\KK_{\tau(i)}\mi [ x_{i,0}^- , \ P_{-,+,-} + P_{+,-,-}  ]$ can be expressed as a polynomial in $x_{j,r}^+$ ($j \neq i$) with coefficients in $\C$. 
This can be done in the same way as in \cite[Lemma 9.14]{Przez-23}. 

\subsection{Case: $a_{i,\tau(i)} = 2$} 

This case occurs only when $N = 2n-1$ is odd and $i=n$. 

\Lem \label{lem: app a2}
We have 
\eq \label{eq: app a2}
[F_n, Q_n]_{q^2} =  \gamma \sum_{\mp} \sum_{k=0}^{n-2} P(x_{n\mp 1,-1}^+, x^+_{n\mp2,0}, \cdots, x^+_{n\mp k \mp1,0}) P(x^+_{n\pm1,0}, \cdots, x^+_{n\pm k\pm1,0}), 
\eneq
for $\gamma = -o(n{-}1) \rho q^2 c_n c_\delta^{1/2}\CCC\mi  \in \C$. 
\enlem 

\Proof
One can show by induction that 
\begin{align*}
[F_n, Q_n]_{q^2} &= \gamma \sum_{\on{id},\tau} \sum_{k=0}^{n-2} P(F_{n-k-2}, \cdots, F_1, F_n, \cdots,F_{2n-1},F_0)  \times \\
&\quad \times \Kb_{n-1} \cdots \Kb_{n-k-1} P(E_{n+1}, \cdots, E_{n+ k+1}), 
\end{align*} 
where $\sum_{\on{id},\tau}$ indicates the symmetrization of the subscripts of all the generators with respect to $\tau$. The result now follows by applying the formula 
\[
P(F_{n-k-2}, \cdots, F_1, F_n, \cdots,F_{2n-1},F_0) \Kb_{n-1} \cdots \Kb_{n-k-1} = P(x_{n- 1,-1}^+, x^+_{n-2,0}, \cdots, x^+_{n-k-1,0}), 
\]
whose derivation can be found inside the proof of \cite[Lemma 9.14]{Przez-23}. 
\enproof

It follows from Lemma \ref{lem: app a2} that $[F_n, Q_n]_{q^2}$ can be expressed as a polynomial in $x_{j,r}^+$  ($j \neq i$) with coefficients in $\C$ and, therefore, commutes with $\avarm{n}$. Hence $\mathbf{R}_n(z) = 0$. 

\subsection{Case: $a_{i,\tau(i)} = -1$}

This case is slightly more nuanced than the previous two. 

\Lem \label{lem: Qi a-1} 
If $a_{i,\tau(i)} = -1$ then 
\[
Q_{\tau(i)} = Q_{\tau(i)}' + Q_{\tau(i)}'' + Q_{\tau(i)}''', 
\]
for 
\begin{align*}
Q'_{\tau(i)} &= o(\tau(i)) q\mi \CCC\mi P(E_{i} \Kb_{\tau(i)}, P(P(F_{n-1}, \cdots, F_1), P(F_{n+2}, \cdots, F_N, F_0))) \KK_{\tau(i)}, 
\end{align*}
element $Q_{\tau(i)}''$ satisfying
\[
\KK_{\tau(i)}\mi [F_i, Q''_{\tau(i)}]_{q\mi} = \gamma \sum_{k=0}^{n-2} P(x_{n\mp 1,-1}^+, x^+_{n\mp2,0}, \cdots, x^+_{n\mp k \mp1,0}) P(x^+_{n\pm2,0}, \cdots, x^+_{n\pm k\pm2,0}), 
\]
with $\gamma = o(i)c_\delta^{1/2} \rho q^{-2} \CCC\mi \in \C$, the upper (resp.\ lower) sign in $\mp, \pm$ corresponding to $i = n$ (resp.\ $i = n+1$),  
and some $Q'''_{\tau(i)} \in \Uu_{d_{i} \geq 2, d_{\tau(i)} \geq 1, +}$. 
\enlem 

\Proof
One proceeds in the same way as in the proof of Lemma \ref{lem: app a2}, using the formulae from Lemma \ref{lem: T B pol C even case}. 
\enproof

Lemma \ref{lem: Qi a-1} implies that $\KK_{\tau(i)}\mi [F_i, Q''_{\tau(i)}]_{q\mi}$ can be expressed as a polynomial in $x_{j,r}^+$  ($j \neq i$) with coefficients in $\C$ and, therefore, commutes with $\avarm{i}$. Hence $\mathbf{R}_i(z) = 0$. 

\addtocontents{toc}{\SkipTocEntry}
\section*{Conflict of interest statement} 

The authors have no competing interests to declare that are relevant to the content of this
article. 

\addtocontents{toc}{\SkipTocEntry}
\section*{Data availability statement} 

The authors declare that all the data supporting the findings of this article are available within the paper.


\newcommand{\etalchar}[1]{$^{#1}$}
\providecommand{\MR}{\relax\ifhmode\unskip\space\fi MR }
\providecommand{\MRhref}[2]{%
  \href{http://www.ams.org/mathscinet-getitem?mr=#1}{#2}
}
\providecommand{\href}[2]{#2}

\end{document}